\RequirePackage{fix-cm}
\documentclass[smallextended]{svjour3}       %
\usepackage{booktabs}

\usepackage{amsthm}

\makeatletter
\def\thmt@innercounters{section,equation,theorem}
\makeatother

\usepackage{springerpaper}
\usepackage[hyperpageref]{backref}
\renewcommand*{\backref}[1]{}%
\renewcommand*{\backrefalt}[4]{%
\ifcase #1 %
No citations.%
\or
\marginpar{\tiny cited on p. #2}%
\else
\marginpar{\tiny cited on pp. #2}%
\fi
}
\makeatletter
\newcommand*{\eqdef}{\mathrel{\vcenter{\baselineskip0.5ex \lineskiplimit0pt
                     \hbox{\scriptsize.}\hbox{\scriptsize.}}}%
                     =}
\makeatother

\newcommand{\skcp}{s_{k, \textup{cp}}}
\newcommand{\varphicp}{\varphi_{\textup{cp}}}
\newcommand{\xicp}{\xi_{\textup{cp}}}
\newcommand{\xikcp}{\xi_{k, \textup{cp}}}
\newcommand{\alphakCG}{\alpha_k^{\text{\tiny CG}}}
\newcommand{\alphakCR}{\alpha_k^{\text{\tiny CR}}}
\newcommand{\alphakcp}{\alpha_{k, \textup{cp}}}
\newcommand{\skCG}{s_k^{\text{\tiny CG}}}
\newcommand{\skCR}{s_k^{\text{\tiny CR}}}

\newcommand{\calS}{\mathcal{S}}
\newcommand{\calSpos}{\calS_{> 0}}
\newcommand{\calSneg}{\calS_{\leq 0}}
\newcommand{\calU}{\mathcal{U}}
\newcommand{\cupdot}{\mathbin{\mathaccent\cdot\cup}}
\newcommand{\cahiernumber}{81}  %
\def\cahieryear{2026}  %
\gitline{}             %
\usepackage[T1]{fontenc}   %
\usepackage{mathptmx}      %

\usepackage{tikz}
\usepackage{pgfplots}
\pgfplotsset{compat=1.18}
\usepgfplotslibrary{groupplots}
\usepackage{tikzscale}
\usetikzlibrary{external}
\graphicspath{%
    {figs/generated/}%
    {figs/}%
    {figs/R2N/}%
    {figs/R2NLS/}%
}

\newcommand*{\includetikzgraphics}[2][]{%
      \includegraphics[#1]{#2}%
}

\makeatletter
\renewcommand{\todo}[2][]{\tikzexternaldisable\@todo[#1]{#2}\tikzexternalenable}
\makeatother
\renewcommand{\theenumi}{(\arabic{enumi})}
\spdefaulttheorem{problemassumption}{Problem Assumption}{\bf}{\it}%
\spdefaulttheorem{modelassumption}{Model Assumption}{\bf}{\it}%
\spdefaulttheorem{stepassumption}{Step Assumption}{\bf}{\it}%
\spdefaulttheorem{assumption}{Assumption}{\bf}{\it}%
\crefname{problemassumption}{Problem Assumption}{Problem Assumptions}
\Crefname{problemassumption}{Problem Assumption}{Problem Assumptions}
\crefname{modelassumption}{Model Assumption}{Model Assumptions}
\Crefname{modelassumption}{Model Assumption}{Model Assumptions}
\crefname{stepassumption}{Step Assumption}{Step Assumptions}
\Crefname{stepassumption}{Step Assumption}{Step Assumptions}
\Crefname{assumption}{Assumption}{Assumptions}
\Crefname{section}{Section}{Sections}
\crefname{section}{\S}{\S}
\Crefname{subsection}{Section}{Sections}
\crefname{subsection}{\S}{\S}
\journalname{Mathematical Programming Computation}
\def\thistitle{A Modified Quasi-Newton Method for Unconstrained Optimization with Nonpositive Curvature Detection and Structured Complexity}
\def\authorFR{Farhad Rahbarnia}
\def\authorDM{Dominique Monnet}
\def\authorDPO{Dominique Orban}
\newcommand\nsercACK[1]{was partially supported by NSERC Discovery Grant RGPIN–#1}
\hypersetup{
  pdftitle={\thistitle},
  pdfauthor={\authorFR{}, \authorDM{}, and \authorDPO{}},
  pdfsubject={Unconstrained Optimization},
  pdfkeywords={Unconstrained optimization, Modified quasi-Newton methods, Regularization methods, Nonpositive curvature, Worst-case complexity, Julia},
}
\begin{document}

\title{%
  \thistitle
  \thanks{%
    \authorFR{} was partially supported by Alexander Graham Bell Canada Graduate Scholarships - Doctoral number 568744-2022.%
    \authorDPO{} \nsercACK{2020–06535}.
  }
}

\titlerunning{Modified Quasi-Newton Method with Nonpositive Curvature Detection}        %

\author{\authorFR{} \and
  \authorDM{} \and
  \authorDPO{}}

\institute{%
  \authorFR{} \at
  GERAD and Department of Mathematics and Industrial Engineering, Polytechnique Montr\'eal.
  \email{farhad.rahbarnia@polymtl.ca}
  \and
  \authorDM{} \at
  Université de Rennes, INSA Rennes, CNRS, IRMAR-UMR 6625, F-35000 Rennes, France.
  \email{dominique.monnet@insa-rennes.fr}
  \and
  \authorDPO{} \at
  GERAD and Department of Mathematics and Industrial Engineering, Polytechnique Montr\'eal.
  \email{dominique.orban@gerad.ca}           %
}

\date{Received: date / Accepted: date}

\pagestyle{myheadings}

\maketitle
\thispagestyle{mytitlepage}

\begin{abstract}
  We consider the convergence, worst-case complexity, and implementation of R2N, a modified quasi-Newton method for large-scale unconstrained smooth optimization.
  At each iteration, we compute a step that approximately minimizes the sum of a quadratic model of the objective with an adaptive quadratic regularization, and that results in sufficient decrease.
  The objective is only assumed to be \(C^1\), although second derivatives may be used if they are available.
  The implementation features a broad range of step-computation strategies, from sparse factorization to matrix-free Krylov methods, that satisfy global convergence requirements.
  When detected, directions of nonpositive curvature are exploited.
  Our analysis allows unbounded model Hessians, provided they satisfy a mild growth condition.
  We derive a worst-case evaluation complexity bound that separates the progress made along positive- and nonpositive-curvature directions.
  By forcing the regularization parameter to vanish near isolated local minimizers, we obtain superlinear local convergence.
  The implementation specializes naturally to nonlinear least-squares problems.
  We provide a modular Julia implementation, and use it in a numerical study on the CUTEst test set that assesses the treatment of nonpositive curvature, the subproblem solver, the update of the regularization parameter, and direct versus iterative subproblem solves.
  An Armijo-Goldstein line search along nonpositive-curvature directions, MINRES-QLP as subproblem solver, and a regularization parameter driven to zero near a solution form our recommended configuration.
  With exact second derivatives, R2N is competitive with a trust-region method, solving slightly more problems with fewer objective and gradient evaluations.
  With convex models, the trust-region method performs as well or better.
\end{abstract}

\section{Introduction}

We consider the solution of
\begin{equation}%
  \label{eq:nlp}
  \minimize{x \in \R^n} \ f(x),
\end{equation}
where \(f: \R^n \to \R\) is \(C^1\) by way of a modified Newton or quasi-Newton method with quadratic regularization.
At each iteration, a step is computed by approximately minimizing the sum of a quadratic model of \(f\) about the current iterate and a multiple of the squared step norm.
The multiple in question is a regularization parameter that is updated according to the decrease in \(f\) induced by the step.
Such methods recently received renewed attention, e.g., from \citet{kanzow-steck-2023}, who study direct solvers based on limited-memory quasi-Newton approximations, and \citet{diouane-habiboullah-orban-2024b}, who derive convergence and complexity properties under weak assumptions, even when the objective includes a nonsmooth regularizer.
Quadratic regularization is related to line search methods and is a well-known alternative to trust-region methods \citep{conn-gould-toint-2000}.
One of its appealing features is that subproblems are unconstrained.
It is famously used in the method of \citet{levenberg-1944} and \citet{marquardt-1963} for nonlinear least squares.

We propose solver R2N, whose salient features are
\begin{itemize}
  \item a global convergence analysis that follows from \citep{diouane-habiboullah-orban-2024b} under the sole assumption that \(f \in C^1\), with model Hessians allowed to grow unbounded at a controlled rate, and under a lower bound on Armijo-Goldstein step lengths when directions of nonpositive curvature are exploited;
  \item a structured worst-case evaluation complexity analysis that separates positive- and nonpositive-curvature steps under Lipschitz continuity of \(\nabla f\) and \(\nabla^2 f\) at the iterates along the steps computed by the algorithm that corrects misleading bounds from \citep{liu-roosta-2022b}:
    \begin{itemize}
      \item if the sequence of model Hessians \(\{B_k\}\) remains bounded, we show that accepted steps along which curvature is positive contribute a term proportional to \(\epsilon^{-2}\) to the total number of evaluations to bring \(\|\nabla f(x_k)\|\) below \(\epsilon > 0\);
      \item when \(\{B_k\}\) is allowed to grow unbounded, this first bound deteriorates---see \Cref{prop:succ-pc};
      \item when \(B_k = \nabla^2 f(x_k)\), accepted steps along which curvature is nonpositive contribute a term proportional to \(\epsilon^{-3/2}\);
      \item this second bound may further deteriorate when the model curvature is inexact along the step---see \Cref{thm:decrease-npc};
    \end{itemize}
  \item superlinear or quadratic local convergence near an isolated local minimizer for a variant in which the regularization parameter is driven to zero;
  \item a modular implementation that is suited for large-scale problems, allows multiple iterative and direct subproblem solvers, and accommodates exact, limited-memory quasi-Newton, or other Hessian approximations;
  \item a variant for nonlinear least squares named R2NLS that may be viewed as an implementation of the method of \citet{levenberg-1944} and \citet{marquardt-1963}, to which our global convergence and complexity results apply, and our local convergence analysis applies when \(\sum_i f_i(x_k) \nabla^2 f_i(x_k) \to 0\), in particular for zero-residual problems or when the residuals are linear;
  \item a thorough numerical study, on unconstrained and nonlinear least-squares problems from the CUTEst collection, of the choices that our analysis leaves open: the treatment of nonpositive curvature, the subproblem solver, the update of the regularization parameter, and direct versus iterative subproblem solves, together with a comparison with a trust-region method.
\end{itemize}

Our current implementation supports the following iterative solvers: the truncated conjugate-gradient method of \citet{dembo-steihaug-1983}, the truncated conjugate residuals method of \citet{dahito-orban-2019}, the truncated minimum residual method of \citet{liu-roosta-2022}, where \emph{truncated} refers to the identification of directions of nonpositive curvature.
R2N currently supports the direct solvers MA57 \citep{duff-2004} and MA97 \citep{hogg-scott-2011}.
When using limited-memory BFGS model Hessians, we also support the direct solver of \citet{erway-jain-marcia-2014}.

One of our objectives is to settle, through a thorough numerical study, the choices that our analysis leaves open.
We use exact second derivatives, limited-memory BFGS approximations and Gauss-Newton models, and compare R2N with the trust-region solver Trunk using the same model, subproblem solver and stopping test.
Although R2N also supports MA97, only MA57 results are reported.
Exploiting nonpositive curvature with the Armijo-Goldstein line search is more efficient and more robust than falling back on the Cauchy step, and comparable to the other strategies we consider.
MINRES-QLP is more efficient than truncated CG, and driving \(\sigma_k\) to zero near a solution improves efficiency without affecting robustness.
Exact subproblem solves, whether by sparse factorization or by the solver of \citet{erway-jain-marcia-2014}, bring no consistent advantage over MINRES-QLP\@.
With exact second derivatives, R2N solves slightly more problems than Trunk and requires fewer objective and gradient evaluations, while Trunk requires fewer Hessian-vector products.
With convex models, Trunk performs as well or better.

\subsection*{Related Research}

\citet{li-fukushima-qi-yamashita-2004}, \citet{li-li-2009} and \citet{polyak-2009} study the properties of the regularized Newton method under the assumption that \(f\) is convex.
\citet{ueda-yamashita-2010} further study the properties of a regularized Newton method with Armijo line search for nonconvex \(f\) under the assumption that \(f\) is \(C^2\).
At each iteration, they set the regularization parameter based on the smallest eigenvalue of the Hessian and the gradient norm, thereby ensuring that subproblems are convex.
They establish a \(O(\epsilon^{-2})\) complexity bound if \(\nabla^2 f\) is Lipschitz continuous and superlinear convergence under a local error bound condition.
Under similar assumptions, \citet{ueda-yamashita-2014} replace the Armijo line search with a parameter update that plays a role in the definition of the regularization parameter.
They establish a complexity bound in \(O(\epsilon^{-2})\), but do not show that it is sharp.
\citet{zhang-ni-2018} combine a similar strategy with a non-monotone Armijo line search.

\citet{kanzow-steck-2023} study a method similar to ours but rely on limited-memory quasi-Newton methods and variants of the direct method developed by \citet{erway-jain-marcia-2014} to solve the first-order stationarity conditions of the model exactly.
They assume that model Hessians remain bounded, and do not exploit nonpositive-curvature directions, but merely increase the regularization parameter when they occur.

\citet{gratton-jerad-toint-2025} propose a regularized Newton method for nonconvex \(f\) in which the shift has the form \(\sqrt{\sigma_k \|\nabla f(x_k)\|}\), and \(\sigma_k\) is updated as in R2N.
When the smallest eigenvalue of the Hessian, possibly restricted to a Krylov subspace, is sufficiently negative, they take a step of prescribed length along an approximate eigenvector, and otherwise absorb the negative curvature into the shift.
They establish a bound in \(O(|\log \epsilon| \, \epsilon^{-3/2})\) under Lipschitz continuity of the exact Hessian.
By contrast, R2N uses \(\sigma_k\) itself as the shift, requires no eigenvalue estimate, exploits any direction of nonpositive curvature detected by the subproblem solver in a line search, and allows quasi-Newton models.

Sharp complexity bounds for the regularized (quasi-)Newton method were first established by \citet{diouane-habiboullah-orban-2024b}, even when model Hessians, whether they are exact or not, are allowed to grow unbounded. %
Because their focus is on nonsmooth problems, they solve subproblems with a first-order method and have no concept of curvature.
We build upon their convergence analysis, which is conducted under the sole assumption that \(f \in C^1\), and discuss multiple subproblem solvers.
We also benefit from their complexity analysis, and supplement it to emphasize the contribution of nonpositive-curvature steps to the overall complexity bound.
Indeed, contrary to the references above, we allow models to be nonconvex, and we exploit nonpositive curvature wherever possible, in a manner similar to trust-region methods.
Our results correct several results by \citet{liu-roosta-2022b}, who study a regularized Newton method with line search and MINRES as a subproblem solver, and claim an overall \(O(\epsilon^{-3/2})\) complexity bound.

\subsection*{Organization of the Paper}

The paper is organized as follows.
In \cref{sec:background}, we review quasi-Newton and regularization methods.
In \cref{sec:models}, we introduce the models and the Cauchy decrease condition, and state Algorithm R2N in \cref{sec:algorithm}. 
In \cref{sec:step}, we describe subproblem solvers that satisfy the Cauchy decrease condition, and our treatment of nonpositive curvature.
We establish global convergence in \cref{sec:convergence}.
In \cref{sec:complexity}, we derive the structured complexity bound, and show in \cref{sec:slow,sec:conv-npc} that its two terms are attained, at least when exact derivatives are used.
In \cref{sec:local}, we present a variant with fast local convergence, and specialize R2N to nonlinear least squares \cref{sec:nlls}. 
Finally, we describe our implementation and reports numerical experiments \cref{sec:implementation}. 

\subsection*{Notation}

We use uppercase and lowercase Latin letters for matrices and vectors, respectively.
Lowercase Greek letters denote scalars.
The Euclidean norm of vector \(x\) is denoted \(\|x\|\).
Accordingly, the spectral norm of matrix \(A\) is \(\|A\|\).
Sets are denoted by uppercase calligraphic letters such as \(\mathcal{S}\), and their cardinality by \(|\mathcal{S}|\).
A subscript \(k\) denotes an index or iteration number.
For positive sequences \(\{\alpha_k\}\) and \(\{\beta_k\}\), \(\alpha_k = O(\beta_k)\) means that there exists a constant \(C > 0\) such that \(\alpha_k \leq C \beta_k\) for all \(k\) large enough, and \(\alpha_k = o(\beta_k)\) means that \(\lim_{k \to \infty} \alpha_k / \beta_k = 0\).
If \(A\) is symmetric, \(A \succ 0\) and \(A \succeq 0\) mean that it is positive definite and positive semi-definite, respectively.

\section{Background}%
\label{sec:background}

Traditional quasi-Newton methods for unconstrained optimization compute steps \(s_k\) that approximately minimize a local quadratic model of the objective about the current iterate \(x_k\) \citep{dennis-more-1977}.
The model Hessian \(B_k = B_k^T\) is typically updated using first-order information to satisfy the secant equation \(B_{k} s_{k-1} = y_{k-1}\), where \(y_{k-1} \eqdef \nabla f(x_k) - \nabla f(x_{k-1})\), and \(x_k = x_{k-1} + s_{k-1}\) \citep{dennis-more-1977}.
In a standard line-search framework, the step \(s_k\) is computed by solving \(B_k s_k = -\nabla f(x_k)\), either exactly or inexactly.
A step length \(\alpha_k > 0\) is then computed to satisfy the Armijo condition, which requires \(s_k\) to be a direction of sufficient descent.
While it is not strictly necessary for \(B_k\) to be positive definite for \(s_k\) to be a descent direction, as is the case, for example, with the Symmetric Rank-1 (SR1) update, imposing \(B_k \succ 0\) is a standard way to guarantee descent.
However, when the underlying objective is highly non-convex, maintaining positive-definite approximations can lead to poor performance.

Trust-region methods are able to exploit directions of nonpositive curvature naturally by imposing a bound \(\|s\| \leq \Delta_k\) \citep{conn-gould-toint-2000}.
Rather than imposing a hard constraint on the step norm, quadratic regularization methods augment the quadratic model with an adaptive penalty term \(\frac{1}{2} \sigma_k \|s\|^2\), where \(\sigma_k > 0\).
The step \(s_k\) is then computed as an approximate solution to the shifted system
\begin{equation}%
  \label{eq:shifted-system-background}%
  (B_k + \sigma_k I)
  s_k = -\nabla f(x_k).
\end{equation}
Curiously, this regularized quasi-Newton approach is not mentioned in \citet{cartis-gould-toint-2022}, although it has received recent attention in the literature, for smooth \citep{li-fukushima-qi-yamashita-2004,li-li-2009,ueda-yamashita-2010,ueda-yamashita-2014,kanzow-steck-2023}, and nonsmooth optimization \citep{diouane-habiboullah-orban-2024b}.
If \(\sigma_k\) is large enough, \(B_k + \sigma_k I \succ 0\), and~\eqref{eq:shifted-system-background} admits a unique minimizer.
However, maintaining a uniformly large \(\sigma_k\) artificially restricts the step size and impedes fast local convergence.
The algorithm must therefore dynamically adjust \(\sigma_k\) based on the predictive quality of the model, while gracefully handling iterations where \(B_k + \sigma_k I \not\succ 0\).
That can be achieved by detecting and exploiting directions of nonpositive curvature---see the definition below.

\section{Models and Cauchy Decrease}%
\label{sec:models}

We consider two types of models of \(f\) about the current iterate \(x_k\).
The first, linear, model is defined as the first-order expansion of \(f\) about \(x_k\):
\begin{equation}%
  \label{eq:def-phi-cp}
  \varphicp(s; x_k) \eqdef f(x_k) + \nabla f(x_k)^T s.
\end{equation}
The \emph{Cauchy point} at iteration \(k\), denoted \(\skcp(\alphakcp)\), is defined as an appropriate multiple \(\alphakcp > 0\) of \(-\nabla f(x_k)\).
We intentionally leave some freedom in the choice of \(\alphakcp\).
In our context, inspired from \citep{diouane-habiboullah-orban-2024b}, it is convenient to view \(\skcp(\alphakcp) \eqdef -\alphakcp \nabla f(x_k)\) as the unique minimizer of
\begin{equation}%
  \label{eq:def-m-cp}
  m(s; x_k, \alphakcp) \eqdef \varphicp(s; x_k) + \tfrac{1}{2} \alphakcp^{-1} \|s\|^2.
\end{equation}
We denote the decrease in \(\varphicp(\cdot; x_k)\) achieved by \(\skcp(\alphakcp)\) as
\begin{equation}%
  \label{eq:def-xi-cp}
  \begin{aligned}
    \xicp(x_k, \alphakcp) & \eqdef \varphicp(0; x_k) - \varphicp(\skcp; x_k)
    \\
                          & = \alphakcp \|\nabla f(x_k)\|^2 = \alphakcp^{-1} \|\skcp(\alphakcp)\|^2 > 0,
  \end{aligned}
\end{equation}
where the final inequality holds as long as \(\nabla f(x_k) \neq 0\).

The second, quadratic, model is
\begin{equation}
  \label{eq:def-phi}
  \varphi(s; x_k) \eqdef \varphicp(s; x_k) + \tfrac{1}{2} s^T B_k s,
\end{equation}
where \(B_k = B_k^T \in \R^{n \times n}\) may be \(\nabla^2 f(x_k)\) if the latter is well-defined, available, and can be obtained at reasonable expense, or an approximation to it.
Model~\eqref{eq:def-phi} serves to compute a step \(s_k\) at iteration \(k\) as an approximate minimizer of
\begin{equation}%
  \label{eq:def-m}
  m(s; x_k, \sigma_k) \eqdef \varphi(s; x_k) + \tfrac{1}{2} \sigma_k \|s\|^2,
\end{equation}
where \(\sigma_k > 0\) is a regularization parameter adjusted at each iteration, and whose role is to drive convergence.
The step \(s_k\) need not minimize \(m(\cdot; x_k, \sigma_k)\) accurately, but should achieve \emph{sufficient decrease}, in the sense that
\begin{equation}%
  \label{eq:def-xi}
  \xi(x_k, \sigma_k) \eqdef \varphi(0; x_k) - \varphi(s_k; x_k)
\end{equation}
should satisfy
\begin{equation}%
  \label{eq:def-cauchy-decrease}
  \xi(x_k, \sigma_k) \geq \tau \, \xicp(x_k, \nu_k)
  \quad \text{for all } k \in \N,
  \quad \text{where } 0 < \tau \leq 1,
\end{equation}
where
\begin{equation}%
  \label{eq:def-nuk}
  \nu_k \eqdef \frac{ \theta_1 }{ \|B_k\| + \sigma_k },
  \quad (0 < \theta_1 < 1).
\end{equation}
\citet{diouane-habiboullah-orban-2024b} use the special value \(\alphakcp = \nu_k\) to define the Cauchy point because \(\skcp(\nu_k)\) is guaranteed to induce a decrease in~\eqref{eq:def-m} as well.
Our algorithm below allows to do the same, although \(\nu_k\) can be overly restrictive in the smooth context.
Thus, we provide other options.

Importantly, \(\skcp(\alphakcp)\) need not necessarily be computed during the iterations.
Instead, we will establish that the procedures used to compute \(s_k\) satisfy~\eqref{eq:def-cauchy-decrease} under certain assumptions, and we only compute \(\skcp(\alphakcp)\) when those assumptions fail.

The next result parallels \citep[Proposition~\(3.7\)]{diouane-habiboullah-orban-2024b} but allows for a step \(s_k\) that achieves a fraction of the decrease in \(m(\cdot; x_k, \sigma_k)\) induced by \(\skcp(\nu_k)\).

\begin{proposition}%
  \label{prop:sufficient-decrease}
  At a general iterate \(x_k \in \R^n\), let \(\sigma_k > 0\) be the regularization parameter, \(B_k = B_k^T\) be the model Hessian, and \(\skcp(\nu_k) \eqdef -\nu_k \nabla f(x_k)\).
  Assume that there exists \(\tau_0 > 0\) such that \(s_k\) satisfies
  \begin{equation}%
    \label{eq:sufficient-decrease}
    m(0; x_k, \sigma_k) - m(s_k; x_k, \sigma_k) \geq \tau_0 \, (m(0; x_k, \sigma_k) - m(\skcp(\nu_k); x_k, \sigma_k)).
  \end{equation}
  Then,~\eqref{eq:def-cauchy-decrease} holds with \(\tau \eqdef \min(1, \, \tau_0 (1 - \tfrac{1}{2} \theta_1))\).
\end{proposition}

\begin{proof}
  We combine~\eqref{eq:def-xi-cp},~\eqref{eq:def-m},~\eqref{eq:def-xi},~\eqref{eq:def-nuk}, and our assumptions to obtain
  \begin{align*}
    m(0; x_k, \sigma_k) - m(s_k; x_k, \sigma_k) & = \xi(x_k, \sigma_k) - \tfrac{1}{2} \sigma_k \|s_k\|^2
    \\
                                                & \geq \tau_0 \, (m(0; x_k, \sigma_k) - m(\skcp(\nu_k); x_k, \sigma_k))
    \\
                                                & = \tau_0 \, (\xicp(x_k, \nu_k) - \tfrac{1}{2} \skcp(\nu_k)^T (B_k + \sigma_k I) \skcp(\nu_k))
    \\
                                                & \geq \tau_0 \, (\xicp(x_k, \nu_k) - \tfrac{1}{2} (\|B_k\| + \sigma_k) \|\skcp(\nu_k)\|^2)
    \\
                                                & = \tau_0 \, (\xicp(x_k, \nu_k) - \tfrac{1}{2} (\|B_k\| + \sigma_k) \nu_k \xicp(x_k, \nu_k))
    \\
                                                & = \tau_0 \, (1 - \tfrac{1}{2} \theta_1) \xicp(x_k, \nu_k)
    \\
                                                & \geq \tau \, \xicp(x_k, \nu_k),
  \end{align*}
  where \(\tau\) is as defined in the statement of the present result.
  We combine the first and last lines above to obtain \(\xi(x_k, \sigma_k) \geq \tau \, \xicp(x_k, \nu_k) + \tfrac{1}{2} \sigma_k \|s_k\|^2 \geq \tau \, \xicp(x_k, \nu_k)\).
\end{proof}

In \Cref{sec:step}, we review several strategies to compute a step that satisfies the assumption of \Cref{prop:sufficient-decrease}.
First, we state a general algorithmic framework in the next section.

\section{Algorithm Statement}%
\label{sec:algorithm}

Our algorithm follows that of \citep{diouane-habiboullah-orban-2024b} with the differences that we need not systematically compute a Cauchy point in the outer loop, and that we have explicit treatment of directions of nonpositive curvature.
The latter is also a major difference with the algorithms of \citep{liu-roosta-2022b} and \citep{kanzow-steck-2023}.
Specifically, when such a direction is detected during the step computation, we perform an Armijo-Goldstein line search in order to ensure controlled decrease along said direction.
The Armijo-Goldstein line search is the key mechanism that allows us to quantify the decrease achieved along a nonpositive-curvature direction in terms of the user-defined stopping tolerance---see \Cref{thm:decrease-npc}.
First, we assume that steps are generated along \emph{gradient-related} directions.

\begin{definition}[Gradient-Related Direction]%
  \label{def:gradient-related-dir}
  We call \(d \in \R^n\) a \emph{gradient-related direction} at \(x_k\) if there exist \(\kappa_1 > 0\) and \(\kappa_2 > 0\), independent of \(k\), such that
  \begin{enumerate}
    \item \(\nabla f(x_k)^T d \leq -\kappa_1 \|\nabla f(x_k)\|^2\), and
    \item \(\|d\| \leq \kappa_2 \|\nabla f(x_k)\|\).
  \end{enumerate}
\end{definition}

If \(d_k \in \R^n\) is a gradient-related direction, the first part of \Cref{def:gradient-related-dir} and the Cauchy-Schwarz inequality yield \(\|\nabla f(x_k)\| \, \|d_k\| \geq |\nabla f(x_k)^T d_k| = -\nabla f(x_k)^T d_k \geq \kappa_1 \|\nabla f(x_k)\|^2\).
Thus,
\begin{equation}%
  \label{eq:gradient-related-dir-lower-bound}
  \|d_k\| \geq \kappa_1 \|\nabla f(x_k)\|.
\end{equation}
Similarly, the second part of \Cref{def:gradient-related-dir} and the Cauchy-Schwarz inequality yield \(-\nabla f(x_k)^T d_k \leq \|\nabla f(x_k)\| \, \|d_k\| \leq \kappa_2 \|\nabla f(x_k)\|^2\).
Thus, we also have
\begin{equation}%
  \label{eq:slope-lower-bound}
  \nabla f(x_k)^T d_k \geq -\kappa_2 \|\nabla f(x_k)\|^2.
\end{equation}

\begin{definition}[Non-Positive Curvature Direction]%
  \label{def:npc-dir}
  We call \(d_k \in \R^n\) a \emph{direction of nonpositive curvature} at \(x_k\) if it is gradient-related and \(d_k^T (B_k + \sigma_k I) d_k \leq 0\).
\end{definition}

We formally state R2N as \Cref{alg:r2n}.
Our formulation uses a constant \(\sigma_{\min}\) that may be chosen to be zero, or a positive value.
This freedom is used to highlight the convergence guarantees that \(\sigma_{\min} > 0\) affords and the extra assumptions required to provide similar guarantees when \(\sigma_{\min} = 0\).
In \Cref{sec:local}, we require \(\sigma_{\min} = 0\) to establish superlinear convergence near an isolated local minimizer.

\begin{algorithm}[ht]%
  \caption[caption]{%
    \label{alg:r2n}
    R2N\@: A modified Quasi-Newton method.
  }
  \begin{algorithmic}[1]%
    \State Choose constants \(\sigma_{\min} \geq 0\), \(\theta_2 > 1\), \(0 < \rho_L < \rho_U < 1\), \(0 < \eta_1 \leq \eta_2 < 1\) and \(0 < \gamma_3 \leq 1 < \gamma_1 \leq \gamma_2\).
    \State Choose \(x_0 \in \R^n\), \(\sigma_0 \geq \sigma_{\min}\).
    \For{\(k = 0, 1, \dots\)}
    \State%
    \label{alg:r2n:Bk}%
    Choose \(B_k = B_k^T \in \R^{n \times n}\).
    \State%
    \label{alg:r2n:step-nuk}%
    Compute Cauchy step \(\skcp\) using \Cref{alg:cauchy-step} %
    \State%
    \label{alg:r2n:step-computation}%
    Compute a step \(s_k\) that satisfies~\eqref{eq:sufficient-decrease} or identify a direction of nonpositive curvature \(d_k\) at \(x_k\).
    \If{a direction of nonpositive curvature \(d_k\) at \(x_k\) was identified}
    \State Compute \(\alpha_k > 0\) such that    \Comment{perform Armijo-Goldstein line search}
    \begin{equation}%
      \label{eq:ag}
      f(x_k) + \rho_U \alpha_k \nabla f(x_k)^T d_k \leq f(x_k + \alpha_k d_k) \leq f(x_k) + \rho_L \alpha_k \nabla f(x_k)^T d_k.
    \end{equation}
    \State%
    \label{alg:r2n:armijo-goldstein-step}%
    Reset \(s_k = \alpha_k d_k\).
    \EndIf
    \If{\(\|s_k\| > \theta_2 \; \|\skcp\|\)}%
    \label{alg:r2n:step-comparison}%
    \State reset \(s_k = \skcp\).
    \EndIf
    \State%
    \label{alg:r2n:step-rhok}%
    Compute the ratio
    \[
      \rho_k \eqdef
      \frac{
        f(x_k) - f(x_k + s_k)
      }{
        \varphi(0; x_k) - \varphi(s_k; x_k)
      }.
    \]
    \If{\(\rho_k \geq \eta_1\)}%
    \label{alg:r2n:step-accept}%
    \State set \(x_{k+1} = x_k + s_k\).
    \Else
    \State set \(x_{k+1} = x_k\).
    \EndIf
    \State%
    \label{alg:r2n:step-update}%
    Update the regularization parameter according to
    \[
      \sigma_{k+1} \in
      \begin{cases}
        \begin{aligned}
           & [\gamma_3 \sigma_k, \, \sigma_k]          &  & \text{ if } \rho_k \geq \eta_2,          &  & \quad \text{very successful iteration} \\
           & [\sigma_k, \, \gamma_1 \sigma_k]          &  & \text{ if } \eta_1 \leq \rho_k < \eta_2, &  & \quad \text{successful iteration}      \\
           & [\gamma_1 \sigma_k, \, \gamma_2 \sigma_k] &  & \text{ if } \rho_k < \eta_1.             &  & \quad \text{unsuccessful iteration}
        \end{aligned}
      \end{cases}
    \]
    \State Reset \(\sigma_{k+1} = \max(\sigma_{k+1}, \, \sigma_{\min})\).
    \EndFor
  \end{algorithmic}
\end{algorithm}

\begin{algorithm}[ht]%
  \caption[caption]{%
    \label{alg:cauchy-step}
    R2N\@: Calculation of Cauchy Step.
  }
  \begin{algorithmic}[1]%
    \Require Gradient $\nabla f(x_k)$, model Hessian \(B_k = B_k^T\), regularization parameter \(\sigma_k\).
    \Require Parameters $0 < \theta_1 < 1$.
    \State Compute \(\gamma_k\) according to~\eqref{eq:def-gamma}.
    \If{$\gamma_k \le 0$}
    \State%
    \label{alg:cauchy-step:nuk}%
    Set $\alphakcp = \theta_1 / (\|B_k\| + \sigma_k)$.
    \Else
    \State Set $\alphakcp = \bar{\alpha}_{k, \textup{cp}} = 1 / \gamma_k$.
    \EndIf
    \State \Return $\skcp \eqdef -\alphakcp \nabla f(x_k)$.
  \end{algorithmic}
\end{algorithm}

In \Cref{alg:cauchy-step}, the computation of \(\|B_k\|\) at \Cref{alg:cauchy-step:nuk} may be costly if operator-vector products with \(B_k\) are expensive.
Firstly, although we use the spectral norm by convention, any other norm may be used in its place, as that will only affect constants in the analysis.
Secondly, \citet{diouane-habiboullah-orban-2024b} explain that computing \(\beta_k \geq \mu_{\mathrm{B}} \|B_k\|\) for \(0 < \mu_{\mathrm{B}} \leq 1\) and setting \(\nu_k \eqdef \theta_1 / (\beta_k + \sigma_k)\) with \(\theta_1 < \mu_{\mathrm{B}}\) continues to ensure that our convergence and complexity results hold unchanged.
Estimating \(\|B_k\|\) with an iterative method only requires operator-vector products, and is thus feasible when the latter are cheap, e.g., when \(B_k\) is a limited-memory quasi-Newton approximation.

When operator-vector products are expensive, an alternative to \(\skcp(\nu_k)\) is to define the Cauchy step as \(-\bar{\alpha}_k \nabla f(x_k)\) where \(\bar{\alpha}_k\) minimizes \(\psi_k(\alpha) \eqdef m(-\alpha \nabla f(x_k); x_k, \sigma_k)\) when it is bounded below.
Specifically,
\begin{equation}
  \label{eq:cauchy-min}
  \bar{\alpha}_{k, \textup{cp}} \eqdef 1 / \gamma_k,
\end{equation}
whenever \(\gamma_k > 0\), where
\begin{equation}%
  \label{eq:def-gamma}
  \gamma_k \eqdef \frac{ \nabla f(x_k)^T (B_k + \sigma_k I) \nabla f(x_k) }{ \|\nabla f(x_k)\|^2 },
\end{equation}
is the model curvature along the gradient direction.

When \(\gamma_k > 0\), \(\psi_k\) is strictly convex, \(\psi_k(0) = f(x_k) = \psi_k(\alpha_{k, \max})\), where \(\alpha_{k, \max} \eqdef 2 / \gamma_k = 2 \bar{\alpha}_{k, \textup{cp}}\), and \(\psi_k(\alpha) < f(x_k)\) for all \(0 < \alpha < \alpha_{k, \max}\).
Hence,~\eqref{eq:sufficient-decrease} holds for \(s_k = \skcp(\alpha_k)\) whenever \(\gamma_k > 0\) and there is \(0 < \delta < 1\) independent of \(k\) such that \((1 - \delta) \bar{\alpha}_{k, \textup{cp}} \leq \alpha_k \leq (1 + \delta) \bar{\alpha}_{k, \textup{cp}}\).
In that situation, \(\alpha_k\) may be used in place of \(\nu_k\) everywhere, or when convenient.
Indeed, such \(\alpha_k\) is guaranteed to result in a fraction of the optimal decrease in \(\psi\), hence a fraction of \(\psi(0) - \psi(\nu_k)\).
Continuing further, any \(s_k\) that results in a fraction of the decrease \(\psi(0) - \psi(\alpha_k)\) also satisfies~\eqref{eq:sufficient-decrease}.

Treatment of directions of nonpositive curvature does not appear in \citep{diouane-habiboullah-orban-2024b}.
Such directions indicate that \(m(\cdot; x_k, \sigma_k)\) is unbounded below.
The algorithm of \citep{diouane-habiboullah-orban-2024b} would simply increase \(\sigma_k\) and compute a new step, or fall back on \(\skcp(\nu_k)\).
Instead, \Cref{alg:r2n} exploits such directions by performing an Armijo-Goldstein line search.

We now justify that, under reasonable assumptions, a step length that satisfies the Armijo-Goldstein conditions along a nonpositive-curvature direction yields a step that satisfies~\eqref{eq:def-cauchy-decrease}.

Our first result establishes that there exists an interval of positive step lengths that satisfy the Armijo-Goldstein line search conditions, whether the search direction  is one of nonpositive curvature or not.
We state this result as we were not able to find it in the literature, although it is very likely known.

\begin{proposition}
  Let \(f: \R^n \to \R\) with \(f \in C^1\) be bounded below, \(x \in \R^n\), and \(d\) be a descent direction for \(f\) from \(x\), i.e., \(\nabla f(x)^T d < 0\).
  Let also \(0 < \rho_L < \rho_U < 1\).
  Then, there exist \(0 < \alpha_U < \alpha_L\) (possibly dependent on \(x\)) such that~\eqref{eq:ag} holds for all \(\alpha_U \leq \alpha \leq \alpha_L\).
\end{proposition}

\begin{proof}
  Denote \(\phi(\alpha) = f(x + \alpha d)\).
  Since \(f\) is bounded below, so is \(\phi\).
  Because \(\phi'(0) = \nabla f(x)^T d < 0\) and \(\phi'\) is continuous, there exists \(\delta > 0\) such that \(\phi'(\alpha) < 0\), and hence \(\phi(\alpha) < \phi(0) = f(x)\), for all \(0 < \alpha < \delta\).
  Therefore, because \(0 < \rho_L < 1\) and \(d\) is a descent direction, the line \(\ell_L(\alpha) = f(x) + \rho_L \alpha \nabla f(x)^T d\) is unbounded below, and necessarily intersects the graph of \(\phi\) for some \(0 < \alpha < \delta\).
  Let  \(\alpha_L > 0\) be the smallest \(\alpha > 0\) such that \(\ell_L(\alpha) = \phi(\alpha)\), i.e.,
  \begin{equation}%
    \label{eq:def-alpha_L}
    f(x + \alpha_L d) = f(x) + \rho_L \alpha_L \nabla f(x)^T d.
  \end{equation}
  By construction,
  \begin{equation}%
    \label{eq:armijo}
    f(x + \alpha d) \leq f(x) + \rho_L \alpha \nabla f(x)^T d \quad \text{for all } 0 \leq \alpha \leq \alpha_L.
  \end{equation}
  As a consequence, any such \(\alpha\) satisfies the Armijo condition.

  Define now \(\ell_U(\alpha) = f(x) + \rho_U \alpha \nabla f(x)^T d\).
  By the same reasoning as above, there exists a smallest \(\hat{\alpha}_U > 0\) such that \(\phi(\hat{\alpha}_U) = \ell_U(\hat{\alpha}_U)\), i.e., \(f(x + \hat{\alpha}_U d) = f(x) + \rho_U \hat{\alpha}_U \nabla f(x)^T d\).
  By construction,
  \begin{equation}%
    \label{eq:armijo-U}
    f(x + \alpha d) \leq f(x) + \rho_U \alpha \nabla f(x)^T d \quad \text{for all } 0 \leq \alpha \leq \hat{\alpha}_U.
  \end{equation}
  If, by contradiction, it were that \(\alpha_L \leq \hat{\alpha}_U\),~\eqref{eq:def-alpha_L} and~\eqref{eq:armijo-U} would imply
  \[
    f(x) + \rho_U \alpha_L \nabla f(x)^T d \geq f(x + \alpha_L d) = f(x) + \rho_L \alpha_L \nabla f(x)^T d.
  \]
  But because \(\nabla f(x)^T d < 0\), this would imply \(\rho_U \leq \rho_L\), which is a contradiction.
  Hence, \(\alpha_L > \hat{\alpha}_U\), and \([\hat{\alpha}_U, \alpha_L] \neq \varnothing\) and is not a singleton.

  There may exist other values \(\hat{\alpha}_U < \alpha < \alpha_L\) such that \(\phi(\alpha) = \ell_U(\alpha)\).
  If so, let \(\alpha_U\) be the largest such \(\alpha\)---note that \(\alpha_U < \alpha_L\) because \(\ell_U(\alpha) < \ell_L(\alpha)\) for all \(\alpha > 0\).
  If not, set \(\alpha_U = \hat{\alpha}_U\).
  By continuity,
  \begin{equation}%
    \label{eq:goldstein}
    f(x + \alpha d) \geq f(x) + \rho_U \alpha \nabla f(x)^T d \quad \text{for all } \alpha_U \leq \alpha \leq \alpha_L,
  \end{equation}
  i.e., the Goldstein condition is satisfied.
\end{proof}

The next result states that if \(d_k\) is a nonpositive-curvature direction at \(x_k\) and \(\alpha_k > 0\) satisfies~\eqref{eq:ag}, \(s_k \eqdef \alpha_k d_k\) satisfies~\eqref{eq:def-cauchy-decrease} under mild assumptions.

\begin{proposition}%
  \label{prop:npc-dir-cauchy}
  Let \(f: \R^n \to \R\) with \(f \in C^1\) be bounded below.
  At iteration \(k\), let \(x_k \in \R^n\) be the current iterate, \(d_k \in \R^n\) satisfy \Cref{def:npc-dir}, and \(\alpha_k > 0\) satisfy~\eqref{eq:ag}.
  Assume that either (i) \(\sigma_{\min} > 0\) in \Cref{alg:r2n}, or that (ii) there exist \(\kappa_3 > 0\) independent of \(k\) such that \(d_k^T B_k d_k \leq -\kappa_3 \|d_k\|^2\) and \(\kappa_b > 0\) such that \(\|B_k\| \geq \kappa_b\) for all \(k\).
  Assume also that there exists \(\alpha_{\min} > 0\) independent of \(k\) such that \(\alpha_k \geq \alpha_{\min}\).
  There exists \(\tau > 0\) sufficiently small such that \(s_k \eqdef \alpha_k d_k\) satisfies~\eqref{eq:def-cauchy-decrease}.
\end{proposition}

\begin{proof}
  If \(\sigma_k \geq \sigma_{\min} > 0\), \(d_k^T B_k d_k \leq -\sigma_k \|d_k\|^2 \leq -\sigma_{\min} \kappa_1^2 \|\nabla f(x_k)\|^2\).
  Otherwise, our assumption ensures that \(d_k^T B_k d_k \leq -\kappa_3 \kappa_1^2 \|\nabla f(x_k)\|^2\).
  We continue the proof under the first assumption, but the second case is similar.
  \Cref{def:npc-dir} yields
  \begin{align*}
    \varphi(0; x_k) - \varphi(s_k; x_k) & = -\alpha_k \nabla f(x_k)^T d_k - \tfrac{1}{2} \alpha_k^2 d_k^T B_k d_k
    \\
                                        & \geq (\kappa_1 \alpha_k + \tfrac{1}{2} \sigma_{\min} \kappa_1^2 \alpha_k^2) \|\nabla f(x_k)\|^2.
  \end{align*}
  Thus,~\eqref{eq:def-cauchy-decrease} holds if there exists \(\tau > 0\) such that \(\kappa_1 \alpha_k + \tfrac{1}{2} \sigma_{\min} \kappa_1^2 \alpha_k^2 \geq \tau \nu_k\).
  Under assumption (i), \(\nu_k \leq \nu_{\max} \eqdef 1 / \sigma_{\min}\) for all \(k\), which under assumption (ii), the same holds with \(\nu_{\max} \eqdef 1 / \kappa_b\).
  Thus, the above holds if \(\kappa_1 \alpha_k + \tfrac{1}{2} \sigma_{\min} \kappa_1^2 \alpha_k^2 \geq \tau \nu_{\max}\).
  Since \(\alpha_k \geq \alpha_{\min}\), it suffices to choose \(\tau \leq (\kappa_1 \alpha_{\min} + \tfrac{1}{2} \sigma_{\min} \kappa_1^2 \alpha_{\min}^2) / \nu_{\max}\).
\end{proof}

In \Cref{prop:npc-dir-cauchy}, Assumptions~(i) and~(ii) are assumptions on the algorithm.
In particular, \(B_k\) must be chosen so nonpositive-curvature directions are directions of sufficient negative curvature for \(B_k\).
The assumption \(\|B_k\| \geq \kappa_b\) simply requires that \(B_k\) not asymptotically vanish.
On the contrary, the assumption on \(\alpha_{\min}\) is an assumption on \(f\).
We now show that it holds, in particular, when \(\nabla f\) is Lipschitz continuous.

\begin{proposition}%
  \label{prop:ag-step-min}
  Let \(d_k \in \R^n\) satisfy \Cref{def:npc-dir} at \(x_k \in \R^n\), and let \(\alpha_k > 0\) satisfy~\eqref{eq:ag}.
  Assume that \(\nabla f\) is Lipschitz continuous with constant \(L > 0\).
  Finally, let \(\kappa_1 > \rho_U \kappa_2\).
  Then, there exists \(\alpha_{\min} > 0\) independent of \(k\) such that \(\alpha_k \geq \alpha_{\min}\).
\end{proposition}

\begin{proof}
  We combine~\eqref{eq:ag} with Lipschitz continuity of \(\nabla f\) to obtain
  \[
    \rho_U \alpha_k \nabla f(x_k)^T d_k \leq
    f(x_k + \alpha_k d_k) - f(x_k) \leq
    \alpha_k \nabla f(x_k)^T d_k + \tfrac{1}{2} L \alpha_k^2 \|d_k\|^2.
  \]
  The left- and right-hand sides combine with \Cref{def:gradient-related-dir} and~\eqref{eq:slope-lower-bound} to give
  \[
    -\rho_U \kappa_2 \alpha_k \|\nabla f(x_k)\|^2 \leq -\kappa_1 \alpha_k \|\nabla f(x_k)\|^2 + \tfrac{1}{2} L \kappa_2^2 \alpha_k^2 \|\nabla f(x_k)\|^2.
  \]
  After simplifying by \(\alpha_k \|\nabla f(x_k)\|^2\), there remains
  \begin{align*}
    \alpha_k \geq \alpha_{\min} \eqdef \frac{ 2 (\kappa_1 - \rho_U \kappa_2) }{ L \kappa_2^2 } > 0.
    \tag*{\qedhere}
  \end{align*}
\end{proof}

\section{Step Computation}%
\label{sec:step}

We turn to practical procedures to compute \(s_k\).
The only requirement on \(s_k\) to preserve the convergence properties of \citep{diouane-habiboullah-orban-2024b} is that it satisfy~\eqref{eq:def-cauchy-decrease}.
For reference, we note that the first-order optimality conditions of the subproblem are the symmetric linear system
\begin{equation}%
  \label{eq:kkt-model}%
  (B_k + \sigma_k I)
  s = -\nabla f(x_k).
\end{equation}
It is possible that an (inexact) solution of~\eqref{eq:kkt-model} be a descent direction for \(m(\cdot; x_k, \sigma_k)\) from \(s = 0\) even if \(B_k + \sigma_k I \not \succeq 0\).
We now outline several procedures to compute an exact or inexact solution of~\eqref{eq:kkt-model} that satisfies~\eqref{eq:def-cauchy-decrease}.

\subsection{Exact Minimizer with Limited-Memory Quasi-Newton Approximation}%
\label{sec:exact-min}

When \(B_k\) is a limited-memory quasi-Newton approximation, and in particular, an L-BFGS approximation, \citet{erway-jain-marcia-2014} describe a procedure to solve~\eqref{eq:kkt-model}, which yields the unique exact minimizer of~\eqref{eq:def-m}.
Clearly, this choice of \(s_k\) satisfies~\eqref{eq:sufficient-decrease}.
The procedure scans the pairs \(\{s_i, y_i\}\), \(i = 1, \ldots, q\), that compose the approximation, where \(y_i \eqdef \nabla f(x_i + s_i) - \nabla f(x_i)\), and \(q \in \N_0\) is the memory parameter.
Overall, it requires \(3 q^2 + 6 q + 3\) dot products and \(2 q^2 - 2 q + 1\) vector updates.
In addition, it requires \(2 q\) system solves with \(B_{k, 0} + \sigma_k I\) as coefficient, where \(B_{k, 0}\) is the initial L-BFGS approximation at iteration \(k\).
In our case, this coefficient is a multiple of the identity.
In practice, \(q\) is typically around \(5\) or \(10\), which makes this procedure relatively cheap.

The conditions of \Cref{def:gradient-related-dir} are satisfied provided \(B_k + \sigma_k I\) is uniformly bounded and uniformly positive definite.

\subsection{The Truncated Conjugate-Gradient Method}%
\label{sec:truncated-cg}

The conjugate-gradient (CG) method of \citet{hestenes-stiefel-1952} was designed to solve symmetric positive-definite linear systems.
In the absence of positive definiteness, \citet{steihaug-1983} argued that it remains relevant for optimization in the sense that it may be applied to~\eqref{eq:kkt-model} with appropriate stopping conditions, and yield a step that induces sufficient decrease in \(m(\cdot; x_k, \sigma_k)\).
The procedure proposed in \citep{steihaug-1983} is called the truncated conjugate-gradient method, which we refer to as \emph{truncated CG}.%
Even though we do not assume that \(B_k + \sigma_k I \succ 0\), in order for this approach to apply, we assume that the first CG iteration is well defined, i.e.,
\begin{equation}%
  \label{eq:cg-curvature-iter-1}
  \nabla f(x_k)^T (B_k + \sigma_k I) \nabla f(x_k) > 0.
\end{equation}
When~\eqref{eq:cg-curvature-iter-1} does not hold, \(d_k = -\nabla f(x_k)\) is a direction of nonpositive curvature, and we perform an Armijo-Goldstein line search along that direction.

CG generates iterates \(s_{k, j}\), \(j = 0, 1, \ldots\), starting with \(s_{k, 0} \eqdef 0\), by way of an exact line search procedure \(s_{k, j+1} = s_{k, j} + \alpha_{k, j} p_{k, j}\), where \(\alpha_{k, j} > 0\) is a step length and \(p_{k, j}\) is a search direction.
The \(s_{k, j}\) are increasingly better approximations of an exact solution of~\eqref{eq:kkt-model} in the sense that \(m(s_{k, j}; x_k, \sigma_k)\) decreases monotonically as \(j\) increases as long as \(p_{k, j}^T (B_k + \sigma_k I) p_{k, j} > 0\).
If \(B_k + \sigma_k I \succ 0\), CG will, in exact arithmetic, terminate after at most \(n\) iterations with \(s_{k, n} = s_k\), an exact solution of~\eqref{eq:kkt-model}.
In practice, it is not necessary to achieve optimality of the subproblem.
Indeed, by the monotonicity property just mentioned, \(m(s_{k, j}; x_k, \sigma_k) \leq m(s_{k, 1}; x_k, \sigma_k)\) for all \(j \geq 1\).
Below, we justify that \(s_{k, 1}\) satisfies~\eqref{eq:sufficient-decrease}, which continues to hold if we perform further CG iterations.

If \(B_k + \sigma_k I \not \succ 0\), CG will generate a search direction \(p_{k, j}\) such that \(p_{k, j}^T (B_k + \sigma_k I) p_{k, j} \leq 0\), at least in exact arithmetic.
We call such \(p_{k, j}\) a direction of nonpositive curvature.
The logic of truncated CG is that we apply the CG updates as long as no such \(p_{k, j}\) is generated, and terminate when the residual norm \(\|(B_k + \sigma_k I) s + \nabla f(x_k)\|\) is smaller than a prescribed tolerance.
If a direction of nonpositive curvature \(p_{k, j}\) is generated at iteration \(j > 0\) of CG, the iterate \(s_{k, j}\) is returned.
If instead \(j = 0\),~\eqref{eq:cg-curvature-iter-1} does not hold because \(p_{k, 0} \eqdef -\nabla f(x_k)\), and the direction \(-\nabla f(x_k)\) is returned as it is known to be a descent direction for \(m(\cdot; x_k, \sigma_k)\) from \(s = 0\).

\begin{proposition}%
  \label{prop:cg-cauchy}
  Let~\eqref{eq:cg-curvature-iter-1} hold.
  The iterate \(s_{k, j}\) returned by truncated CG satisfies~\eqref{eq:sufficient-decrease}.
\end{proposition}

\begin{proof}
  By assumption, truncated CG stops at an iteration \(j \geq 1\).
  At iteration \(j = 1\), \(s_{k, 1}\) is selected as the unique global minimizer of the convex quadratic \(\alpha \mapsto m(-\alpha \nabla f(x_k); x_k, \sigma_k)\).
  CG iterations \(2, \ldots, j-1\) further decrease \(m(\cdot; x_k, \sigma_k)\), and therefore,~\eqref{eq:sufficient-decrease} must be satisfied with \(\tau_0 = 1\).
\end{proof}

Recall that, when~\eqref{eq:cg-curvature-iter-1} holds, \(s_{k, 1} = \alpha_{k, 0} p_{k, 0}\) where \(p_{k, 0} \eqdef -\nabla f(x_k)\) and
\begin{equation}%
  \label{eq:alpha-cg}
  \alpha_{k, 0} = \alphakCG \eqdef \bar{\alpha}_{k, \textup{cp}},
\end{equation}
where \(\bar{\alpha}_{k, \textup{cp}}\) is defined in~\eqref{eq:cauchy-min}.
We also define \(\skCG \eqdef -\alphakCG \nabla f(x_k)\).

In CG, detecting directions of nonpositive curvature is straightforward as \(p_{k, j}^T (B_k + \sigma_k I) p_{k, j}\) is computed once per iteration.
When that quantity is nonpositive, CG terminates and returns
\begin{itemize}
  \item the search direction \(-\nabla f(x_k)\) if \(j = 0\), and flags that direction as a direction of nonpositive curvature;
  \item the iterate \(s_{k, j}\) as approximate solution if \(j > 0\), and flags \(p_{k, j}\) as a direction of nonpositive curvature.
\end{itemize}
In the second case, one may choose to simply use \(s_{k, j}\) as step \(s_k\) as it satisfies~\eqref{eq:sufficient-decrease}, or to exploit the direction of nonpositive curvature via the Armijo-Goldstein line search.

The conditions of \Cref{def:gradient-related-dir} are satisfied provided the projection of \(B_k + \sigma_k I\) in the subspace generated by the \(p_{k, j}\) during the CG iterations is uniformly bounded and uniformly positive definite.

\subsection{The Truncated Conjugate-Residual Method}%
\label{sec:truncated-cr}

The conjugate-residual (CR) method was also proposed by \citet{hestenes-stiefel-1952}.
Contrary to CG, it is designed to solve~\eqref{eq:kkt-model} by minimizing its residual norm.
Similarly to CG, CR generates iterates \(s_{k, j}\) from \(s_{k, 0} \eqdef 0\) by way of an exact line search procedure and is known to induce monotonic decrease in \(m(\cdot; x_k, \sigma_k)\) as long as no nonpositive curvature direction has been found \citep{dahito-orban-2019,fong-saunders-2012}.

\citet{dahito-orban-2019} formulate the truncated CR method, which is relevant for optimization and shares the main features of truncated CG\@.
As in \Cref{sec:truncated-cg}, we show that the first truncated CR iterate \(s_{k, 1}\) satisfies~\eqref{eq:sufficient-decrease}.
Further iterations only result in further decrease in \(m(\cdot; x_k, \sigma_k)\) as long as positive curvature is observed.
For the truncated CR method to apply, we again assume that~\eqref{eq:cg-curvature-iter-1} holds for \(s_{k, 1} \eqdef \alpha_{k, 0} p_{k, 0}\) with \(p_{k, 0} \eqdef -\nabla f(x_k)\) as in truncated CG, but this time,
\begin{equation}%
  \label{eq:alpha-cr}
  \alpha_{k, 0} = \alphakCR \eqdef \frac{ \nabla f(x_k)^T (B_k + \sigma_k I) \nabla f(x_k) }{ \nabla f(x_k)^T (B_k + \sigma_k I)^2 \nabla f(x_k)},
\end{equation}
and \(\skCR \eqdef -\alphakCR \nabla f(x_k)\).
We examine two cases.
In the first, \(B_k + \sigma_k I \succ 0\), and~\eqref{eq:sufficient-decrease} holds without further assumptions.
In the second, \(B_k + \sigma_k I \not \succ 0\), and~\eqref{eq:sufficient-decrease} follows provided the left-hand side of~\eqref{eq:cg-curvature-iter-1} remains sufficiently positive and \(\{B_k\}\) is bounded.

In both cases, the following result holds.

\begin{proposition}%
  \label{prop:alpha-cr-cg}
  Assume~\eqref{eq:cg-curvature-iter-1} holds.
  Then, \(\alphakCR \leq \alphakCG\).
\end{proposition}

\begin{proof}
  By the Cauchy-Schwarz inequality,
  \begin{align*}
    (\nabla f(x_k)^T (B_k + \sigma_k I) \nabla f(x_k))^2 & \leq \|\nabla f(x_k)\|^2 \, \|(B_k + \sigma_k I) \nabla f(x_k)\|^2
    \\
                                                         & = \|\nabla f(x_k)\|^2 \, (\nabla f(x_k)^T (B_k + \sigma_k I)^2 \nabla f(x_k)).
  \end{align*}
  Therefore,~\eqref{eq:alpha-cg} and~\eqref{eq:alpha-cr} yield \(\alphakCR \leq \alphakCG\).
\end{proof}

\begin{proposition}%
  \label{prop:cr-pos-def}
  Assume~\eqref{eq:cg-curvature-iter-1} holds.
  Let \(B_k + \sigma_k I \succ 0\).
  Then, \(\nu_k \leq \alphakCR \leq \alphakCG\) and~\eqref{eq:sufficient-decrease} holds with \(\tau_0 = 1\).
\end{proposition}

\begin{proof}
  Because \(B_k + \sigma_k I \succ 0\),
  \[
    \nabla f(x_k)^T (B_k + \sigma_k I)^2 \nabla f(x_k) \leq \|B_k + \sigma_k I\| \, \nabla f(x_k)^T (B_k + \sigma_k I) \nabla f(x_k).
  \]
  Thus,~\eqref{eq:alpha-cr} yields \(\alphakCR \geq 1 / \|B_k + \sigma_k I\| \geq 1 / (\|B_k\| + \sigma_k) \geq \nu_k\).
  \Cref{prop:alpha-cr-cg} provides \(\alphakCR \leq \alphakCG\).
  Because \(\alphakCG\) is the step length to the unique global minimizer of the convex quadratic \(\alpha \mapsto m(-\alpha \nabla f(x_k); x_k, \sigma_k)\) by definition,
  \begin{align*}
    m(0; x_k, \sigma_k) - m(\skcp; x_k, \sigma_k) & \leq m(0; x_k, \sigma_k) - m(\skCR; x_k, \sigma_k)
    \\
                                                  & \leq m(0; x_k, \sigma_k) - m(\skCG; x_k, \sigma_k),
  \end{align*}
  so that~\eqref{eq:sufficient-decrease} holds with \(\tau_0 = 1\).
\end{proof}

Whenever \(B_k \succeq 0\) for all \(k\), as is the case if it is an L-BFGS approximation, among other situations, \Cref{prop:cr-pos-def} shows that the truncated CR method always yields sufficient decrease.
When \(B_k + \sigma_k I \not \succ 0\), we require additional assumptions.

For reference, define
\begin{align}%
  \Delta_k(\alpha) & \eqdef m(0; x_k, \sigma_k) - m(-\alpha \nabla f(x_k); x_k, \sigma_k) \nonumber
  \\
                   & = \alpha (\|\nabla f(x_k)\|^2 - \tfrac{1}{2} \alpha \nabla f(x_k)^T (B_k + \sigma_k I) \nabla f(x_k)).
  \label{eq:gradient-model-decrease}
\end{align}

\begin{proposition}%
  \label{prop:cr-indef}
  Let \(s_k\) be computed by the truncated CR method.
  Assume that there exist \(\gamma_{\min} > 0\), \(\kappa_B > 0\) and \(\sigma_{\max} > 0\) such that
  \begin{equation}%
    \label{eq:sufficient-curvature}
    \nabla f(x_k)^T (B_k + \sigma_k I) \nabla f(x_k) \geq \gamma_{\min} \|\nabla f(x_k)\|^2,
  \end{equation}
  \(\|B_k\| \leq \kappa_B\), and \(\sigma_k \leq \sigma_{\max}\) for all \(k\).
  There exists \(\tau_0 > 0\) such that~\eqref{eq:sufficient-decrease} holds for all \(k\).
\end{proposition}

\begin{proof}
  By~\eqref{eq:gradient-model-decrease}, the model decrease at \(\skcp\) reads
  \begin{align*}
    \Delta_k(\nu_k) & = \frac{\theta_1}{\|B_k\| + \sigma_k} \left( \|\nabla f(x_k)\|^2 - \tfrac{1}{2} \theta_1 \frac{ \nabla f(x_k)^T (B_k + \sigma_k I) \nabla f(x_k) }{ \|B_k\| + \sigma_k } \right)
    \\
                    & \leq \theta_1 \left( 1 - \tfrac{1}{2} \theta_1 \frac{ \gamma_{\min} }{ \kappa_B + \sigma_{\max} } \right) \frac{ \|\nabla f(x_k)\|^2 }{ \|B_k\| + \sigma_k }.
  \end{align*}
  Without loss of generality, we may assume that \(\kappa_B\) and \(\sigma_{\max}\) are large enough that the parenthesized factor above is positive.

  Our assumption on \(\gamma_{\min}\) and~\eqref{eq:alpha-cr} yield \(\alphakCR \geq \gamma_{\min} / (\|B_k\| + \sigma_k)^2\).
  The same Cauchy-Schwarz inequality as in the proof of \Cref{prop:alpha-cr-cg} gives
  \begin{align*}
    \alphakCR \nabla f(x_k)^T (B_k + \sigma_k I) \nabla f(x_k) & = \frac{ (\nabla f(x_k)^T (B_k + \sigma_k I) \nabla f(x_k))^2 }{ \nabla f(x_k)^T (B_k + \sigma_k I)^2 \nabla f(x_k) }
    \\
                                                               & \leq \|\nabla f(x_k)\|^2.
  \end{align*}
  Combined with~\eqref{eq:gradient-model-decrease}, the above shows that the model decrease at \(\skCR\) is
  \[
    \Delta_k(\alphakCR) \geq \tfrac{1}{2} \frac{ \gamma_{\min} }{ (\|B_k\| + \sigma_k)^2 } \|\nabla f(x_k)\|^2 \geq \tfrac{1}{2} \frac{ \gamma_{\min} }{\kappa_B + \sigma_{\max} } \frac{ \|\nabla f(x_k)\|^2 }{ \|B_k\| + \sigma_k }.
  \]
  Finally, we obtain \(\Delta_k(\alphakCR) \geq  \tau_0 \Delta_k(\nu_k)\) with
  \begin{align*}
    \tau_0 \eqdef \tfrac{1}{2} \frac{ \gamma_{\min} }{\kappa_B + \sigma_{\max} } / \left( \theta_1 \left( 1 - \tfrac{1}{2} \theta_1 \frac{ \gamma_{\min} }{ \kappa_B + \sigma_{\max} } \right) \right).
    \tag*{\qedhere}
  \end{align*}
\end{proof}

In \Cref{alg:r2n}, there is no guarantee that \(\{\sigma_k\}\) remains bounded, even if \(\{B_k\}\) is bounded, unless \(\nabla f\) is Lipschitz continuous.
However, under the latter assumption, boundedness of \(\{B_k\}\) implies boundedness of \(\{\sigma_k\}\).
The proof of this assertion is identical to that of \citep[Theorem~\(4.1\)]{aravkin-baraldi-orban-2024}.

In CR, it is possible to detect two types of directions of nonpositive curvature.
As in CG, nonpositive curvature may occur along a search direction \(p_{k, j}\), but it may also occur along a residual \(r_{k, j} = -\nabla f(x_k) - (B_k + \sigma_k I) s_{k, j}\).
\citet{dahito-orban-2019} modify CR to instrument it with nonpositive curvature detection along both directions without incurring the cost of extra operator-vector products with \(B_k + \sigma_k I\).
When nonpositive curvature is detected along either, or both, directions, CR terminates and returns
\begin{itemize}
  \item the search direction \(-\nabla f(x_k)\) if \(j = 0\), and flags that direction as a direction of nonpositive curvature;
  \item the iterate \(s_{k, j}\) as approximate solution if \(j > 0\), and flags either or both of \(p_{k, j}\) and \(r_{k, j}\) as a direction of nonpositive curvature.
\end{itemize}
In the second case, one may simply use \(s_{k, j}\) as step \(s_k\) as it satisfies~\eqref{eq:sufficient-decrease}, or to exploit the direction(s) of nonpositive curvature via the Armijo-Goldstein line search.%

As in the case of truncated CG, the conditions of \Cref{def:gradient-related-dir} are satisfied provided the projection of \(B_k + \sigma_k I\) in the subspace generated by the \(p_{k, j}\) during the CR iterations is uniformly bounded and uniformly positive definite.

\subsection{Modified MINRES}%
\label{sec:modified-minres}

In exact arithmetic, MINRES \citep{paige-saunders-1975} is equivalent to CR as long as no direction of nonpositive curvature is detected.
In practice, MINRES can be significantly more stable than CR in the presence of rounding errors and ill conditioning.
In MINRES as in CR, nonpositive curvature may occur along a search direction \(p_{k, j}\) or along a residual \(r_{k, j} = -\nabla f(x_k) - (B_k + \sigma_k I) s_{k, j}\).
Based on the work of \citet{dahito-orban-2019}, \citet{liu-roosta-2022} modified MINRES to add nonpositive curvature detection.
The results on truncated CR apply to modified MINRES, and the same guarantees hold.

\subsection{Exact Solution via Sparse Factorization}%
\label{sec:sparse-factorization}

When \(\nabla^2 f(x_k)\) exists and is available at acceptable expense, it often occurs that it is large and sparse.
In our implementation of \Cref{alg:r2n}, it is possible to solve~\eqref{eq:kkt-model} via a sparse factorization method.
Such a method typically raises a warning when \(B_k + \sigma_k I\) is singular, and its solution may fail to be a descent direction.
In that event, we declare iteration \(k\) unsuccessful without evaluating \(f\), i.e., we set \(x_{k+1} = x_k\) and increase \(\sigma_k\).
If \(s_k\) is a positive-curvature direction, we now verify that it satisfies~\eqref{eq:def-cauchy-decrease}.

\begin{proposition}%
  \label{prop:sparse-factorization}
  Assume that either \(\sigma_{\min} > 0\) in \Cref{alg:r2n} or that there exists \(\kappa_b > 0\) such that \(\|B_k\| \geq \kappa_b\) for all \(k\).
  Let \(s_k\) be computed by solving the nonsingular system~\eqref{eq:kkt-model} via a sparse factorization, and satisfy \Cref{def:gradient-related-dir}.
  Then, there exists \(\tau > 0\) independent of \(k\) such that~\eqref{eq:def-cauchy-decrease} holds.
\end{proposition}

\begin{proof}
  If \(\sigma_{\min} > 0\), let \(\nu_{\max} \eqdef 1 / \sigma_{\min}\).
  Otherwise, if \(\|B_k\| \geq \kappa_b\) for all \(k\), let \(\nu_{\max} \eqdef 1 / \kappa_b\).
  In both cases, \(\nu_k \leq \nu_{\max}\) for all \(k\).
  Define also \(\tau \eqdef \tfrac{1}{2} \kappa_1 / \nu_{\max}\).
  From~\eqref{eq:kkt-model} and \Cref{def:gradient-related-dir},
  \begin{align*}
    \varphi(0; x_k) - \varphi(s_k; x_k) & = -\nabla f(x_k)^T s_k - \tfrac{1}{2} s_k^T B_k s_k
    \\
                                        & = -\tfrac{1}{2} \nabla f(x_k)^T s_k + \tfrac{1}{2} \sigma_k \|s_k\|^2
    \\
                                        & \geq -\tfrac{1}{2} \nabla f(x_k)^T s_k
    \\
                                        & \geq \tfrac{1}{2} \kappa_1 \|\nabla f(x_k)\|^2
    \\
                                        & = \tau \nu_{\max} \|\nabla f(x_k)\|^2
    \\
                                        & \geq \tau \nu_k \|\nabla f(x_k)\|^2
    \\
                                        & = \tau \xikcp \geq \min(1, \, \tau) \xikcp, \tag*{\qedhere}
  \end{align*}
\end{proof}

The conditions of \Cref{def:gradient-related-dir} are satisfied provided \(B_k + \sigma_k I\) is uniformly bounded and uniformly nonsingular.

\subsection{Treatment of Nonpositive Curvature Directions}%
\label{sec:npc-treatment}%

As indicated above, if \(d_k\) satisfying \Cref{def:npc-dir} is identified using any of the subproblem solvers, we perform an Armijo-Goldstein line search along that direction.
Other strategies are possible.
In truncated CG, CR and MINRES, one may choose to set \(s_k\) to the most recent subproblem solver iterate before nonpositive curvature was identified.
Or one may revert to the Cauchy step.
We compare such strategies in \Cref{sec:implementation}.

As a final remark, under reasonable assumptions, we may bypass the computation of \(\rho_k\) and the acceptance test \(\rho_k \geq \eta_1\) for steps \(s_k = \alpha_k d_k\) computed from a nonpositive-curvature direction and a step length satisfying the Armijo-Goldstein conditions.

\begin{proposition}%
  \label{prop:ag-accept}
  Let \(d_k \in \R^n\) be a direction of nonpositive curvature at \(x_k \in \R^n\), and assume that there exists \(\kappa_4 > 0\) such that \(d_k^T B_k d_k \geq -\kappa_4 \|d_k\|^2\).
  Let \(\alpha_k > 0\) be such that~\eqref{eq:ag} holds, and assume that there exists \(\alpha_{\max} > 0\) such that \(\alpha_k \leq \alpha_{\max}\) for all \(k\).
  Finally, let \(s_k \eqdef \alpha_k d_k\).
  Then, \(\rho_k \geq \tilde{\eta}_1\), where \(\tilde{\eta}_1 \eqdef \rho_L \kappa_1 / (\kappa_2 (1 + \tfrac{1}{2} \kappa_2 \kappa_4 \alpha_{\max})) > 0\).
\end{proposition}

\begin{proof}
  The Armijo condition in~\eqref{eq:ag} and~\eqref{eq:slope-lower-bound} yield
  \[
    f(x_k) - f(x_k + s_k) \geq -\rho_L \alpha_k \nabla f(x_k)^T d_k \geq \rho_L \kappa_1 \alpha_k \|\nabla f(x_k)\|^2.
  \]
  The definition of \(\rho_k\), \Cref{def:npc-dir},~\eqref{eq:slope-lower-bound} again, and our assumption give
  \begin{align*}
    \rho_k & =
    \frac{
      f(x_k) - f(x_k + s_k)
    }{
      -\alpha_k \nabla f(x_k)^T d_k - \tfrac{1}{2} \alpha_k^2 d_k^T B_k d_k
    }
    \\
           & \geq
    \frac{
      \rho_L \kappa_1 \alpha_k \|\nabla f(x_k)\|^2
    }{
      \kappa_2 \alpha_k \|\nabla f(x_k)\|^2 + \tfrac{1}{2} \kappa_4 \alpha_k^2 \|d_k\|^2
    }
    \\
           & \geq
    \frac{
      \rho_L \kappa_1 \alpha_k \|\nabla f(x_k)\|^2
    }{
      \kappa_2 \alpha_k \|\nabla f(x_k)\|^2 + \tfrac{1}{2} \kappa_2^2 \kappa_4 \alpha_k^2 \|\nabla f(x_k)\|^2
    }
    \\
           & =
    \frac{
      \rho_L \kappa_1
    }{
      \kappa_2 + \tfrac{1}{2} \kappa_2^2 \kappa_4 \alpha_k
    }
    \geq \tilde{\eta}_1. \tag*{\qedhere}
  \end{align*}
\end{proof}

The first assumption in \Cref{prop:ag-accept} asserts that nonpositive-curvature directions are not directions of arbitrarily-negative curvature for \(B_k\).
The second assumption requires that step lengths that satisfy the Armijo condition in~\eqref{eq:ag} not be arbitrarily large.
It holds, in particular, when \(f\) satisfies a Polyak-\L{}ojasiewicz-type inequality\footnote{The actual inequality uses the optimal value of \(f\) instead of a lower bound.}, i.e., there exists \(\mu_{\mathrm{PL}} > 0\) such that
\begin{equation}%
  \label{eq:pl-like}
  \|\nabla f(x_k)\|^2 \geq \mu_{\mathrm{PL}} (f(x_k) - f_{\mathrm{low}}),
\end{equation}
where \(f_{\mathrm{low}}\) is a lower bound on \(\{f(x_k)\}\), i.e., \(f(x_k) \geq f_{\mathrm{low}}\) for all \(k\).

\begin{proposition}
  Assume that \(f(x_0) > f_{\mathrm{low}}\), and that \(x_k\) is not stationary.
  Let~\eqref{eq:pl-like} hold.
  Then, \(\alpha_k \leq 1 / (\rho_L \kappa_1 \mu_{\mathrm{PL}})\).
\end{proposition}

\begin{proof}
  The rightmost inequality in~\eqref{eq:ag} and \Cref{def:gradient-related-dir} imply
  \[
    f(x_k) - f(x_k + s_k) \geq -\rho_L \alpha_k \nabla f(x_k)^T d_k \geq \rho_L \kappa_1 \alpha_k \|\nabla f(x_k)\|^2.
  \]
  Thus,~\eqref{eq:pl-like} yields
  \begin{align*}
    \alpha_k \leq
    \frac{
      f(x_k) - f(x_k + s_k)
    }{
      \rho_L \kappa_1 \|\nabla f(x_k)\|^2
    }
    \leq
    \frac{
      f(x_k) - f_{\mathrm{low}}
    }{
      \rho_L \kappa_1 \|\nabla f(x_k)\|^2
    }
    \leq
    \frac{
      1
    }{
      \rho_L \kappa_1 \mu_{\mathrm{PL}}
    }.
    \tag*{\qedhere}
  \end{align*}
\end{proof}

In view of \Cref{prop:ag-accept}, we may accept \(s_k\) without checking whether \(\rho_k \geq \eta_1\), and, in that instance, we leave \(\sigma_k\) unchanged.

\section{Convergence Analysis}%
\label{sec:convergence}

The main differences between \Cref{alg:r2n} and \citep[Algorithm~4.1]{diouane-habiboullah-orban-2024b} are that~\eqref{eq:def-xi-cp} holds as an equality, whereas only \(\xikcp \geq \tfrac{1}{2} \nu_k^{-1} \|\skcp\|^2\) holds in the presence of a nonsmooth term, and \Cref{prop:sufficient-decrease} is less demanding than \citep[Proposition~3.7]{diouane-habiboullah-orban-2024b}.
Because \(h = 0\) in~\eqref{eq:nlp}, the model \(\psi = 0\) satisfies \citep[Model Assumptions~4.1--4.3]{diouane-habiboullah-orban-2024b} immediately, and \citep[Problem Assumption~4.1]{diouane-habiboullah-orban-2024b} reduces to \(f \in C^1\).
Because the Cauchy step simplifies to \(\skcp(\nu_k) = -\nu_k \nabla f(x_k)\), the stationarity measure of \citet{diouane-habiboullah-orban-2024b} reduces exactly to the gradient norm:
\begin{equation}
  \label{eq:xi-smooth-equiv}
  \nu_k^{-1/2} \xikcp^{1/2}
  = \nu_k^{-1/2} \bigl(\nu_k \|\nabla f(x_k)\|^2\bigr)^{1/2}
  = \|\nabla f(x_k)\|.
\end{equation}

The only remaining assumption for global convergence concerns the growth of \(\{B_k\}\), and corresponds to \citep[Model Assumption~4.4]{diouane-habiboullah-orban-2024b}.

\begin{modelassumption}%
  \label{asm:growth-Bk-conv}
  The sequence \(\{B_k\}\) satisfies
  \[
    \sum_{k=0}^{\infty} \frac{1}{1 + \max_{0 \leq j \leq k} \|B_j\|}
    = +\infty.
  \]
\end{modelassumption}

\Cref{asm:growth-Bk-conv} holds whenever \(\{B_k\}\) is bounded, or grows sufficiently slowly with \(k\).
We can now state the main convergence result.

\begin{theorem}[{\protect \citealp[Theorems~5.6 and~5.9]{diouane-habiboullah-orban-2024b}}]%
  \label{thm:global-convergence-mo}
  Let \(f \in C^1\) be bounded below, and \(\{B_k\}\) satisfy \Cref{asm:growth-Bk-conv}.
  Assume that $s_k$ satisfies~\eqref{eq:def-cauchy-decrease} with a constant $\tau > 0$ independent of $k$.
  Then, the sequence of iterates generated by \Cref{alg:r2n} satisfies
  \begin{equation}
    \liminf_{k \to \infty} \|\nabla f(x_k)\| = 0.
  \end{equation}
  If \Cref{alg:r2n} generates only finitely many successful iterations, the boundedness assumption on \(f\) can be removed, and there exists \(x^\star \in \R^n\) such that \(x_k = x^\star\) for all large enough \(k\), and \(\nabla f(x^\star) = 0\).
\end{theorem}

\begin{proof}
  By~\eqref{eq:xi-smooth-equiv}, the stationarity measure of \citet{diouane-habiboullah-orban-2024b} reduces to \(\|\nabla f(x_k)\|\).
  By assumption, $s_k$ satisfies~\eqref{eq:def-cauchy-decrease}, which coincides with the sufficient decrease condition of \citep[Proposition~3.7]{diouane-habiboullah-orban-2024b}.
  Since \Cref{asm:growth-Bk-conv} coincides with \citep[Model Assumption~4.4]{diouane-habiboullah-orban-2024b}, the results follow from \citep[Theorems~5.6 and~5.9]{diouane-habiboullah-orban-2024b}.
\end{proof}

The Cauchy decrease condition~\eqref{eq:def-cauchy-decrease} is satisfied by each subproblem solver under the assumptions stated in \Cref{sec:step}.

To establish worst-case evaluation complexity bounds, \citet{diouane-habiboullah-orban-2024b} additionally make a Lipschitz-continuity assumption on \(\nabla f\).
We expand on these bounds in \Cref{sec:complexity}, distinguishing progress along positive- and nonpositive-curvature directions.
\section{Structured Worst-Case Complexity}%
\label{sec:complexity}

The algorithm of \citet{diouane-habiboullah-orban-2024b} results in a worst-case iteration complexity in \(O(\epsilon^{-2})\) under the assumption that \(\nabla f\) is Lipschitz continuous and \(\{B_k\}\) remains bounded.
However, close examination of their analysis reveals that Lipschitz continuity is not necessary, but that \Cref{asm:grad-remainder} below is sufficient.
Thus, the same bound applies to \Cref{alg:r2n} under the same assumptions.

As we now show, the iteration complexity bound can be refined to distinguish between iterations in which a nonpositive curvature direction was identified and those in which a positive curvature direction was identified.
We show that the former contribute \(O(\epsilon^{-3/2})\) terms to the overall complexity bound, and that the asymptotic order of those terms deteriorates with the quality of the Hessian curvature approximation along the step.
We do so by correcting the analysis of \citet{liu-roosta-2022b}, which contains an error.
Because positive curvature directions contribute \(O(\epsilon^{-2})\) terms to the complexity bound, the tightness of the bound established in \citep{diouane-habiboullah-orban-2024b} shows that there exists \(f\) satisfiying the assumptions below, and for which \(O(\epsilon^{-2})\) iterations are necessary to achieve \(\|\nabla f(x_k)\| \leq \epsilon\).
Thus, contrary to claims in \citep{liu-roosta-2022b}, the overall complexity bound cannot be \(O(\epsilon^{-3/2})\).
However, we show below that the overall complexity bound can be decomposed according to the contributions of positive- and nonpositive-curvature accepted steps.
Moreover, we exhibit an example for which the bound \(O(\epsilon^{-2})\) is attained, and one for which the bound \(O(\epsilon^{-3/2})\) is attained.

We denote by \(\calS\) the set of successful iterations, and partition it according to the curvature along the step:
\[
  \calS \eqdef \{k \in \N \mid \rho_k \geq \eta_1\},
  \mathhfill
  \calS_{> 0} \eqdef \{k \in \calS \mid s_k^T (B_k + \sigma_k I) s_k > 0\},
  \mathhfill
  \calS_{\leq 0} \eqdef \calS \setminus \calS_{> 0}.
\]

\begin{problemassumption}%
  \label{asm:grad-remainder}
  There exists \(L_g \geq 0\) such that, for all \(k \in \N\),
  \[
    |f(x_{k+1}) - (f(x_k) + \nabla f(x_k)^T s_k)| \leq \tfrac{1}{2} L_g \|s_k\|^2,
  \]
  where \(x_k\) and \(s_k\) are the iterate and step generated at iteration \(k\) of \Cref{alg:r2n}.
\end{problemassumption}

For completeness, we begin by restating the decrease along directions of positive curvature as established by \citet{diouane-habiboullah-orban-2024b}.
Their proof does not depend on boundedness of \(\{B_k\}\).
We only require the result below on successful iterations where a positive-curvature step is computed.

We implicitly assume that \(|\calS|\) is infinite for otherwise, the final iterate can be shown to be stationary \citep[Theorem~\(6.3\)]{diouane-habiboullah-orban-2024b}.
The proof of the following result is classic and constitutes the first part of the proof of \citet[Theorem~\(6.4\)]{diouane-habiboullah-orban-2024b}.

\begin{proposition}%
  \label{prop:decrease-pc}
  For \(k \in \calS_{> 0}\),
  \[
    f(x_k) - f(x_{k+1}) \geq \eta_1 \tau \nu_k \|\nabla f(x_k)\|^2.
  \]
\end{proposition}

\begin{proof}
  By definition of \(\rho_k\), the identity \(x_{k+1} = x_k + s_k\),~\eqref{eq:def-xi-cp} and~\eqref{eq:def-cauchy-decrease} yield
  \begin{align*}
    f(x_k) - f(x_{k+1}) \geq \eta_1 \xi(x_k, \sigma_k)
                        \geq \eta_1 \tau \xicp(x_k, \nu_k)
                        \geq \eta_1 \tau \nu_k \|\nabla f(x_k)\|^2. \tag*{\qedhere}
  \end{align*}
\end{proof}

\Cref{prop:decrease-pc} contradicts \citep[Lemma~\(8\)]{liu-roosta-2022b}, where the authors use MINRES to compute \(s_k\).
Indeed, the decrease achieved along positive curvature directions is entirely independent of MINRES's nonpositive-curvature detection abilities.
Below, we show an example where the bound given by \Cref{prop:decrease-pc} is attained, even though~\eqref{eq:kkt-model} is solved exactly, which suggests an error in \citep[Lemma~\(8\)]{liu-roosta-2022b}.

To quantify the objective decrease along nonpositive-curvature directions, we require the following assumption.

\begin{problemassumption}%
  \label{asm:hess-remainder}
  \(f\) is \(C^2\) and there exists \(L_H \geq 0\) such that, for all \(k \in \calS_{\leq 0}\),
  \[
    |f(x_{k+1}) - (f(x_k) + \nabla f(x_k)^T s_k + \tfrac{1}{2} s_k^T \nabla^2 f(x_k) s_k)| \leq \tfrac{1}{3!} L_H \|s_k\|^3,
  \]
  where \(x_k\) and \(s_k\) are the iterate and step generated at iteration \(k\) of \Cref{alg:r2n}.
\end{problemassumption}

Although \Cref{asm:grad-remainder,asm:hess-remainder} are not assumptions on the problem alone, their interest resides in the facts that they are weaker than local or global Lipschitz continuity, and that they allow for the construction of examples that attain the worst-case complexity bounds.
We develop such examples in \Cref{sec:slow,sec:conv-npc}.

In addition, we require model Hessians to match the objective Hessian sufficiently along the step using the same assumption as in \citep{dussault-migot-orban-2024}.
Note that the latter is significantly weaker than that used by \citet{lim-roosta-2025}, who impose a uniform bound of the form \(\|\nabla^2 f(x_k) - B_k\| \leq \delta\) where \(\delta > 0\).

\begin{modelassumption}%
  \label{asm:hess-error}
  There exist \(C_H > 0\) and \(0 < \zeta \leq 1\) such that, for all \(k \in \calS_{\leq 0}\), \(\|(\nabla^2 f(x_k) - B_k) s_k\| \leq C_H \|s_k\|^{1 + \zeta}\), where \(x_k\) and \(s_k\) are the iterate and step generated at iteration \(k\) of \Cref{alg:r2n}.
\end{modelassumption}

The following result corrects \citep[Lemma~\(9\)]{liu-roosta-2022b} and generalizes it to the case of inexact Hessian.

\begin{theorem}%
  \label{thm:decrease-npc}
  Let \Cref{asm:hess-remainder} and \Cref{asm:hess-error} be satisfied.
  Let \(k \in \calS_{\leq 0}\) and assume that \(\rho_U < \kappa_1 / \kappa_2\).
  Then,
  \[
    f(x_k) - f(x_{k+1}) \geq \rho_L \kappa_1 \min \left( c_1 \|\nabla f(x_k)\|^{3/2}, \, c_2 \|\nabla f(x_k)\|^{\frac{2 + \zeta}{1 + \zeta}} \right),
  \]
  where
  \[
    c_1 \eqdef
    \sqrt{
      \frac{
        3 (\kappa_1 - \kappa_2 \rho_U)
      }{
        L_H \kappa_2^3
      }
    },
    \quad \text{and} \quad
    c_2 \eqdef
    \left( \frac{
    \kappa_1 - \kappa_2 \rho_U
    }{
    C_H \kappa_2^{2 + \zeta}
    } \right)^{\frac{1}{1+\zeta}}.
  \]
\end{theorem}

\begin{proof}
  Because \(s_k\) is accepted, \(x_{k+1} = x_k + s_k\).
  By \Cref{def:gradient-related-dir,def:npc-dir}, \Cref{asm:hess-remainder}, \Cref{asm:hess-error}, and \(s_k = \alpha_k d_k\),
  \begin{multline*}
    f(x_{k+1}) - f(x_k) \leq \nabla f(x_k)^T s_k + \tfrac{1}{2} s_k^T B_k s_k + \tfrac{1}{2} s_k^T (\nabla^2 f(x_k) - B_k) s_k + \tfrac{1}{3!} L_H \|s_k\|^3
    \\
               \leq -\kappa_1 \alpha_k \|\nabla f(x_k)\|^2 + \tfrac{1}{2} C_H \kappa_2^{2 + \zeta} \alpha_k^{2 + \zeta} \|\nabla f(x_k)\|^{2 + \zeta} + \tfrac{1}{3!} L_H \kappa_2^3 \alpha_k^3 \|\nabla f(x_k)\|^3.
  \end{multline*}
  The left-hand side inequality of the Armijo-Goldstein conditions combines with~\eqref{eq:slope-lower-bound} to yield \(- \kappa_2 \rho_U \alpha_k \|\nabla f(x_k)\|^2 \leq \rho_U \alpha_k \nabla f(x_k)^T d_k \leq f(x_{k+1}) - f(x_k)\).
  We combine the latter inequality with the above to obtain
  \begin{multline*}
    - \kappa_2 \rho_U \alpha_k \|\nabla f(x_k)\|^2 \leq \\
    -\kappa_1 \alpha_k \|\nabla f(x_k)\|^2 + \tfrac{1}{2} C_H \kappa_2^{2 + \zeta} \alpha_k^{2 + \zeta} \|\nabla f(x_k)\|^{2 + \zeta} + \tfrac{1}{3!} L_H \kappa_2^3 \alpha_k^3 \|\nabla f(x_k)\|^3,
  \end{multline*}
  which simplifies, after dividing by \(\alpha_k \|\nabla f(x_k)\|^2\), to
  \begin{equation}%
    \label{eq:npc-quadratic-alpha}%
    (\kappa_1 - \kappa_2 \rho_U)
    \leq
    \tfrac{1}{3!} L_H \kappa_2^3 \alpha_k^2 \|\nabla f(x_k)\| +
    \tfrac{1}{2} C_H \kappa_2^{2 + \zeta} \alpha_k^{1 + \zeta} \|\nabla f(x_k)\|^{\zeta}.
  \end{equation}
  By assumption, \(\kappa_1 - \kappa_2 \rho_U > 0\).
  There are two possible regimes.

  Firstly, if \(\tfrac{1}{3!} L_H \kappa_2^3 \alpha_k^2 \|\nabla f(x_k)\| \geq \tfrac{1}{2} C_H \kappa_2^{2 + \zeta} \alpha_k^{1 + \zeta} \|\nabla f(x_k)\|^\zeta\),~\eqref{eq:npc-quadratic-alpha} yields
  \[
    (\kappa_1 - \kappa_2 \rho_U) \leq
    \tfrac{1}{3} L_H \kappa_2^3 \alpha_k^2 \|\nabla f(x_k)\|
    \quad \Longleftrightarrow \quad
    \alpha_k \geq c_1 \, \|\nabla f(x_k)\|^{-1/2}.
  \]
  Secondly, if \(\tfrac{1}{3!} L_H \kappa_2^3 \alpha_k^2 \|\nabla f(x_k)\| < \tfrac{1}{2} C_H \kappa_2^{2 + \zeta} \alpha_k^{1 + \zeta} \|\nabla f(x_k)\|^\zeta\),~\eqref{eq:npc-quadratic-alpha} yields
  \[
    (\kappa_1 - \kappa_2 \rho_U) \leq
    C_H \kappa_2^{2 + \zeta} \alpha_k^{1 + \zeta} \|\nabla f(x_k)\|^{\zeta}
    \quad \Longleftrightarrow \quad
    \alpha_k \geq c_2 \, \|\nabla f(x_k)\|^{-\frac{\zeta}{1+\zeta}}.
  \]
  In both cases,
  \[
    \alpha_k \geq
    \min \left( c_1 \, \|\nabla f(x_k)\|^{-1/2} , \ c_2 \, \|\nabla f(x_k)\|^{-\frac{\zeta}{1+\zeta}} \right).
  \]
  We substitute the latter bound on \(\alpha_k\) into the Armijo condition, i.e., the right-hand side inequality of~\eqref{eq:ag}, and use \Cref{def:gradient-related-dir} to obtain
  \[
    f(x_{k+1}) - f(x_k) \leq \rho_L \alpha_k \nabla f(x_k)^T d_k \leq -\rho_L \kappa_1 \alpha_k \|\nabla f(x_k)\|^2.
  \]
  The identity \(2 - \zeta / (1 + \zeta) = (2 + \zeta) / (1 + \zeta)\) yields the desired result.
\end{proof}

\Cref{thm:decrease-npc} implicitly assumes that the step is not reset to the Cauchy step on \Cref{alg:r2n:step-comparison} of \Cref{alg:r2n}, which is a safe assumption if \(\theta_2\) is chosen large enough, given the second part of \Cref{def:gradient-related-dir}.

The proof of \Cref{thm:decrease-npc} simplifies when the curvature along the step is exact and gives a tighter bound.

\begin{corollary}%
  \label{cor:decrease-npc}
  Let \Cref{asm:hess-remainder} hold and assume that for all \(k\), \(s_k^T B_k s_k = s_k^T \nabla^2 f(x_k) s_k\), as would happen, e.g., if \(B_k = \nabla^2 f(x_k)\). 
  Assume that at iteration \(k\), a direction \(d_k\) of nonpositive curvature is identified at \Cref{alg:r2n:step-computation} of \Cref{alg:r2n}, and that \(\rho_U < \kappa_1 / \kappa_2\).
  Finally, assume that the corresponding step \(s_k = \alpha_k d_k\) is accepted.
  Then,
  \[
    f(x_k) - f(x_{k+1}) \geq \rho_L \kappa_1 c_3 \, \|\nabla f(x_k)\|^{3/2},
    \qquad
    c_3 \eqdef
    \sqrt{
      \frac{
        3!
        (\kappa_1 - \kappa_2 \rho_U)
      }{
        L_H \kappa_2^3
      }
    }.
  \]
\end{corollary}

\begin{proof}
  We follow the logic of the proof of \Cref{thm:decrease-npc}, but the terms involving \(\zeta\) are absent.
  We obtain directly the bound
  \[
    \alpha_k \geq c_3 \, \|\nabla f(x_k)\|^{-1/2}.
  \]
  The same logic then yields the result.
\end{proof}

For any $0 < \zeta < 1$, $\frac{2+\zeta}{1+\zeta} > \frac{3}{2}$.
Hence, when $\|\nabla f(x_k)\|$ is sufficiently small, the term in $\|\nabla f(x_k)\|^{\frac{2+\zeta}{1+\zeta}}$ achieves the minimum in \Cref{thm:decrease-npc}.
As a consequence, any degree of inexactness in the curvature approximation deteriorates the worst-case complexity bound relative to the case where the curvature is exact.
The optimal $O(\epsilon^{-3/2})$ complexity along directions of nonpositive curvature is recovered only when $\zeta = 1$ (as shown in \Cref{cor:decrease-npc}) or when the curvature approximation is exact.

Let \(\epsilon > 0\).
We are interested in finding the smallest \(k_\epsilon\) such that \(\|\nabla f(x_{k_\epsilon - 1})\| > \epsilon\) and \(\|\nabla f(x_{k_\epsilon})\| \leq \epsilon\).
Because \Cref{alg:r2n} generates iterates by searching along directions of nonpositive and positive curvature, we define the sets
\begin{subequations}
  \begin{align}
    \calSpos(\epsilon) & \eqdef \{ k < k_\epsilon \mid \rho_k \geq \eta_1 \text{ and } s_k^T (B_k + \sigma_k I) s_k > 0 \}    \\
    \calSneg(\epsilon) & \eqdef \{ k < k_\epsilon \mid \rho_k \geq \eta_1 \text{ and } s_k^T (B_k + \sigma_k I) s_k \leq 0 \} \\
    \calU(\epsilon)    & \eqdef \{ k < k_\epsilon \mid \rho_k < \eta_1 \},
  \end{align}
\end{subequations}
respectively the sets of successful iterations on which positive and nonpositive curvature was detected, and the set of unsuccessful iterations.
We also define
\begin{equation}%
  \label{eq:successful}
  \calS(\epsilon) \eqdef \calSpos(\epsilon) \cupdot \calSneg(\epsilon),
\end{equation}
where \(\cupdot\) denotes the disjoint union.
Clearly, \(|\calS(\epsilon)| = |\calSpos(\epsilon)| + |\calSneg(\epsilon)|\).

\Cref{thm:decrease-npc} readily lets us bound \(|\calS_{\leq 0}(\epsilon)|\) as follows.

\begin{theorem}%
  \label{thm:succ-npc}
  Let \Cref{asm:hess-remainder} and \Cref{asm:hess-error} be satisfied.
  Assume further that \(f\) is bounded below by \(f_{\text{low}}\).
  Then,
  \[
    |\calS_{\leq 0}(\epsilon)| \leq
    \frac{ f(x_0) - f_{\text{low}} }{ \rho_L \kappa_1 \, \min \left( c_1 \epsilon^{3/2}, \, c_2 \epsilon^{\frac{2 + \zeta}{1 + \zeta}} \right) },
  \]
  where \(c_1\) and \(c_2\) are as in \Cref{thm:decrease-npc}.
\end{theorem}

\begin{proof}
  Let \(k \in \calS_{\leq 0}(\epsilon)\).
  By \Cref{thm:decrease-npc}, \(f(x_k) - f(x_{k+1}) \geq \rho_L \kappa_1 \min ( c_1 \epsilon^{3/2}, \, c_2 \epsilon^{\frac{2 + \zeta}{1 + \zeta}} )\).
  We sum those inequalities over \(\calS_{\leq 0}(\epsilon)\) and obtain
  \begin{align*}
    f(x_0) - f_{\text{low}} & \geq
    \sum_{k \in \calS_{\leq 0}(\epsilon)} (f(x_k) - f(x_{k+1})) \\
                            & \geq
    |\calS_{\leq 0}(\epsilon)| \rho_L \kappa_1 \min ( c_1 \epsilon^{3/2}, \, c_2 \epsilon^{\frac{2 + \zeta}{1 + \zeta}} ).
    \tag*{\qedhere}
  \end{align*}
\end{proof}

Under the stronger assumptions of \Cref{cor:decrease-npc}, the bound in \Cref{thm:succ-npc} becomes
\[
  |\calS_{\leq 0}(\epsilon)| \leq
  \frac{ f(x_0) - f_{\text{low}} }{ \rho_L \kappa_1 c_3 } \, \epsilon^{-3/2}.
\]

\citet{diouane-habiboullah-orban-2024b} make two additional assumptions, which we only require on iterations in \(\calS_{> 0}\).
The first assumption concerns adequacy between \(f(x_k)\) and \(\varphi(\cdot; x_k)\), and states that there exists \(\kappa_m > 0\) such that for all \(k \in \calS_{> 0}\),
\[
  |f(x_k + s_k) - \varphi(s_k; x_k)| \leq \kappa_m (1 + \|B_k\|) \, \|s_k\|^2.
\]
By \Cref{asm:grad-remainder} and the definition of \(\varphi(\cdot; x_k)\), the previous inequality is always satisfied with \(\kappa_m = \tfrac{1}{2} \max(1, \, L_g)\).

The second concerns the growth rate of \(\{B_k\}\) and is necessary to quantify the decay of \(\nu_k\) in the bound of \Cref{prop:decrease-pc}.

\begin{modelassumption}%
  \label{asm:growth-Bk}
  There exist \(\mu > 0\) and \(0 \leq p \leq 1\) such that, for all \(k \in \N\),
  \[
    \|B_k\| \leq \mu (1 + |\calS_{> 0, k-1}|^p),
  \]
  where \(\calS_{> 0, k} = \{ 0 \leq j \leq k \mid j \in \calS_{> 0}\}\) and \(\calS_{> 0, -1} = \varnothing\).
\end{modelassumption}

The classic assumption, which requires \(\{B_k\}\) to be bounded, corresponds to \(p = 0\) in \Cref{asm:growth-Bk}.
However, any value \(0 \leq p \leq 1\) is allowed to capture the growth of \(\{B_k\}\).
A similar assumption was first formulated by \citet{leconte-orban-2026} in the context of trust-region methods.
\Cref{asm:growth-Bk} is compatible with \Cref{asm:hess-remainder,asm:hess-error}.
However, in our complexity analysis, they are used in different situations; the first applies along positive-curvature steps while the other two apply along nonpositive-curvature steps.

We are now in position to establish the following tight bound.

\begin{proposition}[{\protect \citealp[Theorem~\(6.4\)]{diouane-habiboullah-orban-2024b}}]%
  \label{prop:succ-pc}
  Let \Cref{asm:grad-remainder} and \Cref{asm:growth-Bk} be satisfied.
  Assume further that \(f\) is bounded below by \(f_{\text{low}}\).
  Then,
  \[
    |\calS_{> 0}(\epsilon)| \leq
    \begin{cases}
      ((1 - p) \kappa_s \epsilon^{-2} + 1)^{1 / (1 - p)} - 1 & \text{ if } \ 0 \leq p < 1 \\
      \exp(\kappa_s \epsilon^{-2}) - 1                       & \text{ if } \ p = 1,
    \end{cases}
  \]
  where
  \begin{align*}
    \kappa_s & \eqdef \frac{ (f(x_0) - f_{\text{low}}) (b_{\max} + 2 \mu (1 + b_{\max})) }{ \eta_1 \theta_1 \tau },
    \\
    b_{\max} & \eqdef \max(\sigma_0 (1 + \|B_0\|)^{-1}, \, 2 \kappa_m \gamma_2 / (1 - \eta_2)).
  \end{align*}
\end{proposition}

Finally, \citet[Theorem~\(6.5\)]{diouane-habiboullah-orban-2024b} establish that
\begin{equation}%
  \label{eq:unsuccessful}
  |\calU(\epsilon)| \leq \log_{\gamma_1}(b_{\max} / \sigma_0) + |\log_{\gamma_1}(\gamma_3)| \, |\calS(\epsilon)| + \log_{\gamma_1}(1 + \mu (1 + |\calS(\epsilon)|^p)).
\end{equation}

\Cref{thm:succ-npc} and \Cref{prop:succ-pc} show that the overall complexity bound contains terms in \(\epsilon^{-2}\) driven by successful iterations along positive-curvature steps, and terms in \(\epsilon^{-3/2} + \epsilon^{-(2 + \zeta)/(1 + \zeta)}\) driven by successful iterations along nonpositive-curvature steps.
Therefore, the worst-case complexity of \Cref{alg:r2n} remains in \(\epsilon^{-2}\), whether or not MINRES is used to compute steps, because no special property of MINRES helps improve that bound along positive-curvature steps.
In that sense, the bound reported in \citep{liu-roosta-2022b} is misleading and erroneous.
Indeed, as shown in the next section, there exists \(f\) satisfying our assumptions and a choice of \(B_k\) such that every step in \Cref{alg:r2n} is a positive-curvature step, every step is computed exactly, every iteration is very successful, and therefore, such that the total number of iterations attains \(\epsilon^{-2}\).

\section{Slow Convergence of \texorpdfstring{\Cref{alg:r2n}}{Algorithm R2N} in the Absence of Nonpositive Curvature}%
\label{sec:slow}

We exhibit an example on which \Cref{alg:r2n} performs exactly a number of successful iterations of the asymptotic order given in \Cref{prop:succ-pc} before finding an approximately stationary point.
The example is directly inspired from similar examples in \citep{leconte-orban-2026,diouane-habiboullah-orban-2026}, and covers the cases of bounded and unbounded model Hessians.
Because the construction of such an example uses well-established results, we refer the reader to \citep{leconte-orban-2026,diouane-habiboullah-orban-2026} for the details.

Let \(0 < \epsilon \leq \tfrac{1}{2}\).
We set \(0 \leq p < 1\) and \(k_\epsilon \eqdef \lceil \epsilon^{-2 / (1 - p)} \rceil\).
The case \(p = 1\) can be handled using an example similar to that in \citep{diouane-habiboullah-orban-2026}.
For \(k = 0, \ldots, k_\epsilon\), we set
\[
  w_k \eqdef (k_{\epsilon} - k) / k_\epsilon, \qquad
  g_k \eqdef - \epsilon (1 + w_k), \qquad
  B_k = (k + 1)^p.
\]
It is clear that \(|g_k| > \epsilon\) for \(k = 0, \ldots, k_\epsilon - 1\), that \(|g_{k_\epsilon}| = \epsilon\), and that \(B_k + \sigma_k \succ 0\) for all \(k = 0, \ldots, k_\epsilon\).
No nonpositive-curvature is ever detected.

The step is computed exactly as
\[
  s_k = -(B_k + \sigma_k)^{-1} g_k = \frac{ \epsilon (1 + w_k) }{ (k+1)^p + \sigma_k },
  \qquad
  k = 0, \ldots, k_\epsilon - 1.
\]
The iterates and function values are generated as \(x_{k+1} \eqdef x_{k} + s_{k}\), and
\begin{equation}%
  \label{eq:def-f}
  f_0 \eqdef \frac{ 4 (1 + \epsilon^{\frac{ 2 }{ 1 - p }}) }{ 1 - p }, \mathhfill
  f_{k+1} \eqdef f_{k} + g_{k} s_{k} = f_{k} - \frac{g_{k}^2}{(k+1)^p + \sigma_k}, \mathhfill
  k = 0, \ldots, k_\epsilon - 1.
\end{equation}

\begin{lemma}%
  \label{lem:f-g-bnd}
  The definitions above ensure that there exists \(\kappa_f > 0\) such that
  \begin{enumerate}
    \item $f_k \ge 0$ for all $k \in \{0, \ldots, k_{\epsilon}\}$, and $\{f_k\}$ is strictly decreasing;
    \item for $k \in \{0, \ldots, k_{\epsilon}-1\}$, $|g_{k+1} - g_k| \le \kappa_f s_k$;
    \item for all $k \in \{0, \ldots, k_{\epsilon}-1\}$, $|f_{k+1} - (f_k + g_k s_k)| \le s_k^2$.
  \end{enumerate}
\end{lemma}

\begin{proof}
  It follows from~\eqref{eq:def-f} that \(\{f_k\}\) is strictly decreasing.
  To show $f_k \geq 0$, we sum the total descent from $f_0$:
  \[
    f_0 - f_{k} = \sum_{i=0}^{k-1} \frac{g_i^2}{(i+1)^p + \sigma_i} \leq \sum_{i=0}^{k-1} \frac{\epsilon^2 (1 + w_i)^2}{(i+1)^p}.
  \]
  Because $(1 + w_i)^2 \le 4$,
  \[
    f_0 - f_{k} \le 4 \epsilon^2 \sum_{i=0}^{k-1} (i+1)^{-p} \le 4 \epsilon^2 \int_0^k x^{-p} dx = \frac{4 \epsilon^2}{1 - p} k^{1 - p}.
  \]
  Since $k \le k_{\epsilon} = \lceil \epsilon^{-2 / (1 - p)} \rceil \leq \epsilon^{-2 / (1 - p)} + 1$, $k^{1 - p} \le (\epsilon^{-2 / (1 - p)} + 1)^{1 - p}$.
  Thus, \(f_0 - f_k \leq 4 (1 + \epsilon^{2 / (1 - p)}) / (1 - p) = f_0\).
  Hence, $f_k \ge 0$.

  For the second part, $|g_{k+1} - g_k| = \epsilon / k_{\epsilon}$.
  Because the step size $0 < s_k < 2 \epsilon / (k + 1)^p$, there exists $\kappa_f$ such that the interpolation bound is satisfied.

  Finally,~\eqref{eq:def-f} yields $|f_{k+1} - f_k - g_k s_k| = 0$.
\end{proof}

\begin{theorem}%
  \label{thm:slow-conv}
  The definitions above ensure existence of \(f: \R \to \R\) continuously differentiable for which exactly $k_\epsilon = \lceil \epsilon^{-2 / (1 - p)} \rceil$ iterations are necessary to produce $x_{k_{\epsilon}}$ such that $\|\nabla f(x_{k_{\epsilon}})\| \le \epsilon$.
  Furthermore, along the sequence of iterates generated, we only encounter strictly positive curvature, solve the subproblem exactly, and each iteration is very successful.
\end{theorem}

\begin{proof}
  Existence of \(f\) that interpolates the values \(f_k\) such that \(\nabla f\) interpolates the values \(g_k\) is guaranteed by \Cref{lem:f-g-bnd} and \citep[Proposition~\(6\)]{leconte-orban-2026}.

  We now show that every iteration is very successful.
  The predicted reduction is
  \begin{align*}
    \varphi(0; x_k) - \varphi(s_k; x_k) & = -g_k s_k - \tfrac{1}{2} s_k^2 B_k                                                                       \\
                                        & = \frac{g_k^2}{(k+1)^p + \sigma_k} - \tfrac{1}{2} \left(\frac{g_k}{(k+1)^p + \sigma_k}\right)^2 (k+1)^p   \\
                                        & = \frac{g_k^2}{(k+1)^p + \sigma_k} \left( 1 - \tfrac{1}{2} \frac{ (k+1)^p }{ (k+1)^p + \sigma_k } \right) \\
                                        & \le \frac{g_k^2}{(k+1)^p + \sigma_k} = f_k - f_{k+1},
  \end{align*}
  where the final value is precisely the actual reduction.
  Therefore, \(\rho_k \geq 1 > \eta_2\).
\end{proof}

\begin{figure}[ht]
  \centering
  \includetikzgraphics[width=\linewidth]{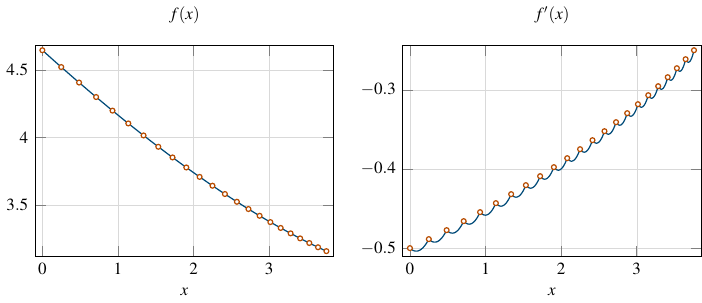}
  \caption{%
    \label{fig:nonpc}
    Graph of \(f\) and \(f'\) in the example of \Cref{sec:slow} with $\epsilon=\tfrac{1}{4}$, and \(p = \tfrac{1}{10}\), which results in \(\lceil \epsilon^{-2 / (1 - p)} \rceil = 22\) iterations.
    Circles indicate the \(23\) iterates \(x_k\).
  }
\end{figure}

\Cref{fig:nonpc} graphs \(f\) and \(f'\) over the interval covered by the iterates \(x_k\).
We set \(\epsilon = \tfrac{1}{4}\) and \(p = \tfrac{1}{10}\).

\section{Slow Convergence of \texorpdfstring{\Cref{alg:r2n}}{Algorithm R2N} in the Presence of Nonpositive Curvature}%
\label{sec:conv-npc}

We now exhibit an example on which \Cref{alg:r2n} performs exactly a number of successful iterations of the asymptotic order given in \Cref{thm:succ-npc} before finding an approximate stationary point.
The example is also inspired by similar examples in \citep{leconte-orban-2026,diouane-habiboullah-orban-2026}, but adapted to account for the Armijo-Goldstein line search along nonpositive-curvature directions.
Let \(0 < \epsilon \leq \tfrac{1}{2}\).
We set \(k_\epsilon \eqdef \lceil \epsilon^{-3/2} \rceil\).
For \(k = 0, \ldots, k_\epsilon\), we set
\[
  w_k \eqdef (k_{\epsilon} - k) / k_\epsilon, \qquad
  g_k \eqdef - \epsilon (1 + w_k), \qquad
  B_k = -\sigma_0 < 0.
\]
It is clear that \(|g_k| > \epsilon\) for \(k = 0, \ldots, k_\epsilon - 1\), that \(|g_{k_\epsilon}| = \epsilon\), and that \(B_k + \sigma_0 = 0\) for all \(k = 0, \ldots, k_\epsilon\).
Because the example will ensure that \(\sigma_k = \sigma_0\) for all \(k\), nonpositive curvature is detected at every iteration.
Because the model curvature is zero, the model \(\varphi(\cdot; x_k)\) is linear, so the algorithm takes the descent direction \(d_k = -g_k\), which is simultaneously the Cauchy direction and a direction of nonpositive curvature.
Our construction depends on two constants \(c_0 > 0\), and \(\bar{\rho} \in (0, 1)\).
Define
\[
  \alpha_k \eqdef c_0 |g_k|^{-1/2}, \qquad
  s_k = \alpha_k d_k = c_0 |g_k|^{1/2},
  \qquad
  k = 0, \ldots, k_\epsilon - 1.
\]
The iterates are generated as \(x_{k+1} \eqdef x_{k} + s_{k}\).
We set
\begin{equation}%
  \label{eq:def-f-npc}
  f_0 \eqdef (2^{3/2} + 1) \bar{\rho} c_0, \quad
  f_{k+1} \eqdef f_k - \bar{\rho} c_0 |g_k|^{3/2}, \quad
  k = 0, \ldots, k_\epsilon - 1.
\end{equation}
Finally, we arrange so the Hessian of the interpolant is \(H_k \eqdef -\sigma_0\), and we set \(B_k = H_k\) in the model.

\begin{lemma}%
  \label{lem:f-g-bnd-npc}
  The definitions of \(g_k\), \(s_k\), \(H_k\), and~\eqref{eq:def-f-npc} ensure that there is \(\kappa_f > 0\) such that
  \begin{enumerate}
    \item%
      \label{itm:fpos}%
      \(\{f_k\}\) is strictly decreasing and \(f_k \ge 0\) for all \(k \in \{0, \ldots, k_{\epsilon}\}\);
    \item%
      \label{itm:fs3}%
      \(|f_{k+1} - (f_k + g_k s_k + \frac{1}{2}H_k s_k^2)| \le \kappa_f s_k^3\) for all \(k \in \{0, \ldots, k_{\epsilon}-1\}\);
    \item%
      \label{itm:gs2}%
      \(|g_{k+1} - (g_k + H_k s_k)| \le \kappa_f s_k^2\) for all \(k \in \{0, \ldots, k_{\epsilon}-1\}\);
    \item%
      \label{itm:Hs}%
      \(|H_{k+1} - H_k| \le \kappa_f s_k\) for all \(k \in \{0, \ldots, k_{\epsilon}-1\}\).
  \end{enumerate}
\end{lemma}

\begin{proof}
  It follows from~\eqref{eq:def-f-npc} that \(\{f_k\}\) is strictly decreasing.
  By definition, \(f_0 > 0\).
  Let \(0 < k \leq k_\epsilon\).
  For \(i = 0, \ldots, k-1\), \(|g_i| \le 2 \epsilon\).
  Hence, because
  \begin{align*}
    f_0 - f_k =
    \sum_{i=0}^{k-1} (f_i - f_{i+1}) =
    \bar{\rho} c_0 \sum_{i=0}^{k-1} |g_i|^{3/2} \leq
    \bar{\rho} c_0 (2 \epsilon)^{3/2} k_\epsilon & \leq
    \bar{\rho} c_0 2^{3/2} (1 + \epsilon^{3/2})
    \\
                                                 & \leq (2^{3/2} + 1) \bar{\rho} c_0 = f_0,
  \end{align*}
  where we used the inequalities \(k_\epsilon \le \epsilon^{-3/2} + 1\), and \(\epsilon^{3/2} \leq 2^{-3/2}\).
  Thus, \(f_k \ge 0\), which establishes~\ref{itm:fpos}.

  For~\ref{itm:fs3}, let \(0 \leq k < k_\epsilon\).
  Because \(g_k < 0\), \(s_k = c_0 |g_k|^{1/2}\), and \(H_k = -\sigma_0\),~\eqref{eq:def-f-npc} yields
  \begin{align*}
    |f_{k+1} - (f_k + g_k s_k + \tfrac{1}{2}H_k s_k^2)| & = |-\bar{\rho} c_0 |g_k|^{3/2} + c_0 |g_k|^{3/2} + \tfrac{1}{2}\sigma_0 s_k^2| \\
                                                        & = (1 - \bar{\rho}) c_0 |g_k|^{3/2} + \tfrac{1}{2}\sigma_0 s_k^2                \\
                                                        & = \frac{1 - \bar{\rho}}{c_0^2} s_k^3 + \tfrac{1}{2}\sigma_0 s_k^2.
  \end{align*}
  Because \(|g_k| > \epsilon\), \(s_k > c_0 \epsilon^{1/2}\), i.e., \(1 < s_k / (c_0 \sqrt{\epsilon})\).
  We multiply both sides of the latter inequality by \(\tfrac{1}{2} \sigma_0 s_k^2\), and insert into the final equality above to obtain
  \[
    |f_{k+1} - (f_k + g_k s_k + \tfrac{1}{2}H_k s_k^2)| <
    \left( \frac{1 - \bar{\rho}}{c_0^2} + \frac{\sigma_0}{2 c_0 \sqrt{\epsilon}} \right) s_k^3.
  \]

  We now turn to~\ref{itm:gs2}.
  By construction, \(|g_{k+1} - g_k| = \epsilon / k_\epsilon\).
  The definitions of \(s_k\) and \(H_k\), and the inequality \(k_\epsilon \ge \epsilon^{-3/2}\) yield
  \[
    |g_{k+1} - (g_k + H_k s_k)| =
    |g_{k+1} - g_k + \sigma_0 s_k| =
    \frac{\epsilon}{k_\epsilon} + \sigma_0 s_k \leq
    \epsilon \, \epsilon^{3/2} + \sigma_0 s_k.
  \]
  Since \(|g_k| \ge \epsilon\), and \(s_k \ge c_0 \sqrt{\epsilon}\), i.e., \(\epsilon \leq s_k^2 / c_0^2\).
  In addition, using an observation similar to that for part~\ref{itm:fs3}, what precedes can be written \(1 \leq s_k / (c_0 \sqrt{\epsilon})\), and we multiply both sides of the latter inequality by \(\sigma_0 s_k\) to obtain \(\sigma_0 s_k \leq \sigma_0 s_k^2 / (c_0 \sqrt{\epsilon})\).
  Thus,
  \[
    |g_{k+1} - (g_k + H_k s_k)|
    \le \left( \frac{\epsilon^{3/2}}{c_0^2} + \frac{\sigma_0}{c_0 \sqrt{\epsilon}} \right) s_k^2.
  \]

  Finally, part~\ref{itm:Hs} holds because \(H_k = H_{k+1}\) for all \(k\).

  The result follows with \(\kappa_f \eqdef \max\left( \frac{\epsilon^{3/2}}{c_0^2} + \frac{\sigma_0}{c_0 \sqrt{\epsilon}}, \, \frac{1 - \bar{\rho}}{c_0^2} + \frac{\sigma_0}{2 c_0 \sqrt{\epsilon}} \right)\).
\end{proof}

Among others,  \Cref{lem:f-g-bnd-npc} shows that \Cref{asm:hess-remainder} holds.
The next result establishes that our construction is consistent with \Cref{alg:r2n}.

\begin{theorem}%
  \label{thm:conv-npc-example}
  Assume \(\rho_U > \max(\rho_L, \eta_2)\), \(\gamma_3 = 1\), and choose \(\bar{\rho} \in (\max(\rho_L, \eta_2), \rho_U)\) in~\eqref{eq:def-f-npc}.%
  Assume further that \(\sigma_0\) and \(c_0\) are chosen such that
  \[
    \sigma_0 \le \min\left( \frac{\theta_1 \theta_2}{2}, \, 2\left(\frac{\bar{\rho}}{\eta_2} - 1\right) \right) \frac{\sqrt{\epsilon}}{c_0}.
  \]
  Then, the definitions of \(f_k\), \(g_k\), \(H_k\), \(B_k\), and \(s_k\) above ensure existence of \(f: \R \to \R\) with \(f \in C^2\) whose gradient and Hessian satisfy \Cref{asm:grad-remainder} and \Cref{asm:hess-remainder}, respectively, and for which exactly \(k_\epsilon = \lceil \epsilon^{-3/2} \rceil\) iterations of \Cref{alg:r2n} are necessary to produce \(x_{k_\epsilon}\) such that \(|\nabla f(x_{k_\epsilon})| \le \epsilon\).
  Along the sequence of iterates, nonpositive curvature is detected at every iteration, \(s_k\) is never reset to the Cauchy step, and every iteration is very successful.
\end{theorem}

\begin{proof}
  Existence of \(f \in C^2\) that interpolates the values \(f_k\) with \(\nabla f(x_k) = g_k\) and \(\nabla^2 f(x_k) = H_k\) follows from \Cref{lem:f-g-bnd-npc} and \citep[Theorem A.9.2]{diouane-habiboullah-orban-2026}.
  Because \(|g_k| > \epsilon\) for \(k < k_\epsilon\) and \(|g_{k_\epsilon}| = \epsilon\), exactly \(k_\epsilon\) iterations are required to reach the stopping criterion.

  Since \(B_k + \sigma_0 = 0\), the model \(\varphi(\cdot; x_k)\) is linear and any descent direction is a direction of nonpositive curvature.
  The algorithm therefore takes \(d_k \eqdef -g_k\), which is simultaneously the Cauchy direction and gradient-related with \(\kappa_1 = \kappa_2 = 1\), and performs the Armijo-Goldstein line search along \(d_k\).
  Using \(g_k d_k = -|g_k|^2\) and \(\alpha_k |g_k|^2 = c_0 |g_k|^{3/2}\), the Armijo-Goldstein condition for \(\alpha_k \eqdef c_0 |g_k|^{-1/2}\) reads
  \[
    f_k - \rho_U c_0 |g_k|^{3/2}
    \le f_{k+1}
    \le f_k - \rho_L c_0 |g_k|^{3/2}.
  \]
  By~\eqref{eq:def-f-npc} and the fact that \(\bar{\rho} \in [\rho_L, \rho_U]\), both inequalities hold and \(\alpha_k\) is admissible.

  Since \(g_k^T (B_k + \sigma_0 I) g_k = 0\), \Cref{alg:cauchy-step:nuk} of \Cref{alg:cauchy-step} sets \(\alphakcp = \theta_1 / (\|B_k\| + \sigma_0) = \theta_1 / (2 \sigma_0)\), so \(\|\skcp\| = \theta_1 |g_k| / (2 \sigma_0)\).
  The reset at \Cref{alg:r2n:step-comparison} of \Cref{alg:r2n} is triggered if \(\|s_k\| > \theta_2 \|\skcp\|\), i.e., \(c_0 |g_k|^{1/2} > \theta_1 \theta_2 |g_k| / (2 \sigma_0)\), i.e., \(2 c_0 \sigma_0 > \theta_1 \theta_2 |g_k|^{1/2}\).
  Because \(|g_k|^{1/2} \ge \sqrt{\epsilon}\), the initialization of \(\sigma_0\) directly ensures \(2 c_0 \sigma_0 \le \theta_1 \theta_2 \sqrt{\epsilon} \le \theta_1 \theta_2 |g_k|^{1/2}\) for all \(k\), so no reset occurs.

  Because \(B_k = -\sigma_0\) and \(s_k^2 = c_0^2 |g_k|\), the predicted and actual reductions are
  \begin{align*}
    \varphi(0; x_k) - \varphi(s_k; x_k) & = -g_k s_k - \tfrac{1}{2} H_k s_k^2 = c_0 |g_k|^{3/2} + \tfrac{1}{2} \sigma_0 c_0^2 |g_k|, \\
    f_k - f_{k+1}                       & = \bar{\rho} c_0 |g_k|^{3/2}.
  \end{align*}
  Because \(|g_k|^{1/2} \ge \sqrt{\epsilon}\), \(\rho_k = \frac{\bar{\rho}}{1 + \tfrac{1}{2} c_0 \sigma_0 |g_k|^{-1/2}} \ge \frac{\bar{\rho}}{1 + \tfrac{1}{2} c_0 \sigma_0 \epsilon^{-1/2}}\).
  Our assumption on \(\sigma_0\) guarantees that \(\tfrac{1}{2} c_0 \sigma_0 \epsilon^{-1/2} \le \bar{\rho}/\eta_2 - 1\), which gives \(\rho_k \ge \eta_2\) at every iteration. 

  Because we set \(\gamma_3 = 1\), the update rules of \Cref{alg:r2n} set \(\sigma_{k+1} = \sigma_k = \sigma_0\) for all \(k\), which confirms that the assumption \(\sigma_k = \sigma_0\) is consistent.
\end{proof}

\begin{figure}[ht]
  \centering
  \includetikzgraphics[width=\linewidth]{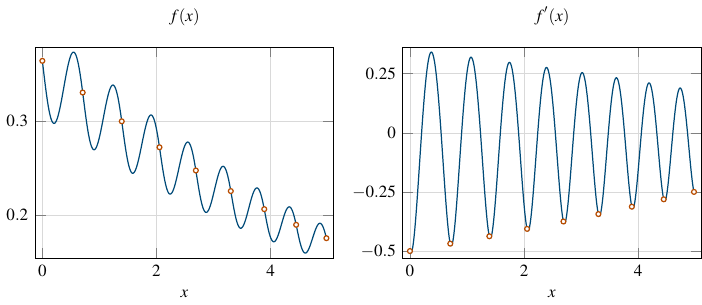}
  \caption{%
    \label{fig:npc}
    Graph of \(f\) and \(f'\) in the example of \Cref{sec:conv-npc} with $\epsilon=\tfrac{1}{4}$, which results in \(\epsilon^{-3/2} = 8\) iterations.
    We use \(\sigma_0 = c_0 = 1\), and \(\bar{\rho} = 0.095 = 2 \eta_2\).
    Circles indicate the \(9\) iterates \(x_k\).
  }
\end{figure}

\Cref{fig:npc} graphs \(f\) and \(f'\) over the interval covered by the iterates \(x_k\).
We use the values \(\epsilon = \tfrac{1}{4}\), \(c_0 = 1\), \(\sigma_0 = 1\).
The assumptions of \Cref{thm:conv-npc-example} are satisfied by choosing any \(0 < \theta_1 < 1\), \(\theta_2 = 4 / \theta_1\), \(0 < \eta_1 \leq \eta_2 < \tfrac{1}{2}\), \(\bar{\rho} = 2 \eta_2\), \(0 < \rho_L \leq \eta_2\), and \(\bar{\rho} < \rho_U < 1\).

On \Cref{fig:npc}, observe that the steps are chosen so that the algorithm makes minimal progress in \(f(x)\), and that negative curvature is observed at each iteration.
The iterates are located just to the left of the local minima of \(f'\), where \(f'\) is decreasing.

The interpolant is constructed using a piecewise Hermite spline of degree \(5\), which is $C^2$.
The domain is partitioned into the intervals $[x_k, x_{k+1}]$, $k = 0, \ldots, k_\epsilon - 1$, and the extra intervals $[x_0 - 1, x_0]$ and $[x_{k_\epsilon}, x_{k_\epsilon} + 1]$.
On each interval except the leftmost and rightmost intervals, the function value, first derivative, and second derivative match $f_k$, $g_k$, $H_k$ at both endpoints.
On the leftmost and rightmost intervals, these three quantities match the sequence at the interior endpoint, while the first and second derivatives vanish and the function value equals the adjacent $f_k$ at the far endpoint.
Outside all intervals, $f$ equals $f_0$ on $(-\infty, x_0 - 1)$ and $f_{k_\epsilon}$ on $(x_{k_\epsilon} + 1, +\infty)$.
The next result gives the explicit Hermite interpolant.

\begin{lemma}[Hermite Interpolation]%
  \label{lem:Hermite-spline}
  For each interval $[x_k, x_{k+1}]$, the unique degree-5 polynomial
  $\pi_k(x) = \sum_{i=0}^5 p_i (x - x_k)^i$ satisfying
  \begin{align*}
    \pi_k(x_k)     & = f_k,     & \pi'_k(x_k)     & = g_k,     & \pi''_k(x_k)     & = -\sigma_0, \\
    \pi_k(x_{k+1}) & = f_{k+1}, & \pi'_k(x_{k+1}) & = g_{k+1}, & \pi''_k(x_{k+1}) & = -\sigma_0,
  \end{align*}
  has coefficients $p_0 = f_k$, $p_1 = g_k$, $p_2 = -\tfrac{1}{2}\sigma_0$, and
  \begin{align*}
    p_3 & = \frac{10E_{f,k} - 4s_k E_{g,k}}{s_k^3}, \quad
    p_4 = \frac{-15E_{f,k} + 7s_k E_{g,k}}{s_k^4}, \quad
    p_5 = \frac{6E_{f,k} - 3s_k E_{g,k}}{s_k^5},
  \end{align*}
  where $E_{f,k} \eqdef f_{k+1} - (f_k + g_k s_k - \tfrac{1}{2}\sigma_0 s_k^2)$
  and $E_{g,k} \eqdef g_{k+1} - (g_k - \sigma_0 s_k)$.
\end{lemma}

\begin{proof}
  Set $t = x - x_k \in [0, s_k]$.
  At $t = 0$, we obtain $p_0 = f_k$, $p_1 = g_k$, and $p_2 = -\sigma_0/2$.
  Setting $u_j \eqdef p_{j+3} s_k^{j+3}$ for $j = 0, 1, 2$, the boundary conditions at $t = s_k$ reduce to
  \[
    \begin{pmatrix}%
      1 & 1  & 1  \\
      3 & 4  & 5  \\
      6 & 12 & 20
    \end{pmatrix}
    \begin{pmatrix}%
      u_0 \\
      u_1 \\
      u_2
    \end{pmatrix}
    =
    \begin{pmatrix}%
      E_{f,k}     \\
      s_k E_{g,k} \\
      0
    \end{pmatrix}
    ,
  \]
  where the third right-hand side vanishes because $H_k = -\sigma_0$ is constant.
  The coefficients given in the statement of the lemma solve the latter linear system.
\end{proof}

\section{A Variant with Fast Local Convergence}%
\label{sec:local}

In the vicinity of an isolated minimizer, we show below that all iterations are eventually very successful, i.e., \(\rho_k \geq \eta_2\).
As stated, however, \Cref{alg:r2n} does not permit fast local convergence, because it prevents \(\sigma_k\) from decaying sufficiently fast.
Variants have been proposed in the literature that enjoy fast local convergence and a favorable worst-case complexity bound.
For instance, instead of updating \(\sigma_k\) as in \Cref{alg:r2n}, \citet{bergou-diouane-kungurtsev-2020} set \(\sigma_k = \mu_k \|\nabla f(x_k)\|^2\) and update \(\mu_k\) as in \Cref{alg:r2n:step-update}.
\citet{bergou-diouane-kungurtsev-royer-2022} proceed similarly with \(\sigma_k = \mu_k \|\nabla f(x_k)\|\).

We now examine a \emph{minimal} change to \Cref{alg:r2n} that still allows us to obtain a worst-case complexity bound and permits fast local convergence.
To that end, we modify the update of \(\sigma_k\) when \(\rho_k \geq \eta_2\) as in \Cref{alg:r2n-fast}.

\begin{algorithm}[ht]%
  \caption[caption]{%
    \label{alg:r2n-fast}
    R2N with Fast Local Convergence.
  }
  \begin{algorithmic}[1]%
    \State Apply \Cref{alg:r2n} with the alternative update
    \[
      \sigma_{k+1} = \gamma_3 \min(\sigma_k, \, \|\nabla f(x_{k+1})\|) \quad \text{if} \quad \rho_k \geq \eta_2
    \]
    in \Cref{alg:r2n:step-update}.
  \end{algorithmic}
\end{algorithm}

\Cref{prop:succ-pc} is unaffected by the change in \Cref{alg:r2n-fast}.
Indeed, \citep[Theorem~\(6.2\)]{diouane-habiboullah-orban-2024b} continues to hold, and is the only prerequisite to \citep[Theorem~\(6.4\)]{diouane-habiboullah-orban-2024b}.
Similarly, \Cref{thm:succ-npc} is also unaffected.
In particular, \(\liminf \|\nabla f(x_k)\| = 0\) continues to hold.
The role of the update of \(\sigma_k\) when \(\rho_k \geq \eta_2\) in the complexity analysis is confined to~\eqref{eq:unsuccessful}.
We now seek a bound on \(|\calU(\epsilon)|\) for \Cref{alg:r2n-fast}.
Unfortunately, the proof technique of \citep[Theorem~\(6.5\)]{diouane-habiboullah-orban-2024b} no longer applies to the modified update.

We first state the following straightforward variation of \citep[Lemma~\(6.1\)]{diouane-habiboullah-orban-2024b}.

\begin{lemma}%
  \label{lem:bsucc}
  Let the assumptions of \Cref{prop:succ-pc} and \Cref{thm:succ-npc} be satisfied.
  Define \(\tilde{b}_{\textup{succ}} \eqdef 2 \kappa_m / (1 - \eta_1)\).
  If \(x_k\) is not stationary and \(\sigma_k (1 + \max_{0 \leq j \leq k} \|B_j\|)^{-1} \geq \tilde{b}_{\textup{succ}}\), iteration \(k\) is successful.
\end{lemma}

The threshold \(\tilde{b}_{\text{succ}}\) of \Cref{lem:bsucc}
is used below
to bound the total number of unsuccessful iterations.

Observe that necessarily, \(k_\epsilon \in \calS(\epsilon)\) because \(\|\nabla f(x_{k_\epsilon})\| \leq \epsilon < \|\nabla f(x_{k_\epsilon - 1})\|\), and hence, \(x_{k_\epsilon - 1} \neq x_{k_\epsilon}\).
Let
\begin{equation}%
  \label{eq:succ-enum}
  \calS(\epsilon) = \{k_0, k_1, \ldots, k_t = k_\epsilon\}
\end{equation}
be an enumeration of the successful iterations, where \(0 \leq k_0 < k_1 < \cdots < k_t\).
The proof of the following result is inspired from that of \citep[Lemma~\(3.3\)]{bergou-diouane-kungurtsev-2020}.

\begin{theorem}%
  \label{thm:unsuccessful-fast}
  Under the assumptions of \Cref{prop:succ-pc} and \Cref{thm:succ-npc},
  \begin{multline*}
    |\calU(\epsilon)| \leq
    1 + \log_{\gamma_1} \left(
    \frac{
      \tilde{b}_{\textup{succ}} (1 + \mu)
    }{
      \sigma_0
    }
    \right) +
    \sum_{i = 0}^{t - 1} \log_{\gamma_1} \left[
      \frac{
        1
      }{
        \gamma_3 \min(\sigma_{k_i}, \epsilon)
      }
      \right]
    + \\
    (|\calS(\epsilon)| - 1) \log_{\gamma_1} [\tilde{b}_{\textup{succ}} (1 + \mu (1 + (|\calS(\epsilon)| - 1)^p))],
  \end{multline*}
  where \(\tilde{b}_{\text{succ}}\) is defined in \Cref{lem:bsucc}, and \(t\) is defined in~\eqref{eq:succ-enum}.
\end{theorem}

\begin{proof}
  Let \(i \in \{0, \ldots, t-1\}\), where \(t\) is defined in~\eqref{eq:succ-enum}.
  We proceed in two steps: firstly, we bound the number \(k_{i+1} - k_i - 1\) of unsuccessful iterations between two consecutive successful iterations, and secondly, we bound the number \(k_0\) of unsuccessful iterations before the first successful iteration.

  For any unsuccessful iteration \(k \in \{k_i + 1, \ldots, k_{i+1} - 1\}\), \Cref{lem:bsucc} implies
  \[
    \sigma_k (1 + \max_{0 \leq j \leq k} \|B_j\|)^{-1} < \tilde{b}_{\text{succ}}.
  \]
  \Cref{alg:r2n-fast} ensures that \(\sigma_{k+1} \geq \gamma_1 \sigma_k\) when \(\rho_k < \eta_1\), \(\sigma_{k+1} \geq \sigma_k\) when \(\eta_1 \leq \rho_k < \eta_2\), and \(\sigma_{k+1} = \gamma_3 \min(\sigma_k, \|\nabla f(x_{k+1})\|)\) when \(\rho_k \geq \eta_2\).
  Because \(\|\nabla f(x_{k_i + 1})\| > \epsilon\),
  \[
    \sigma_k \geq
    \gamma_1^{k - k_i - 1} \sigma_{k_i + 1} \geq
    \gamma_1^{k - k_i - 1} \gamma_3 \min(\sigma_{k_i}, \, \|\nabla f(x_{k_i + 1})\|) \geq
    \gamma_1^{k - k_i - 1} \gamma_3 \min(\sigma_{k_i}, \, \epsilon).
  \]
  We combine the two observations above with \Cref{asm:growth-Bk} and obtain
  \[
    \gamma_1^{k - k_i - 1} \gamma_3 \min(\sigma_{k_i}, \, \epsilon) \leq
    \tilde{b}_{\text{succ}} (1 + \max_{0 \leq j \leq k} \|B_j\|) \leq
    \tilde{b}_{\text{succ}} (1 + \mu (1 + |\calS_{> 0, k}|^p)),
  \]
  where \(\calS_{> 0, k}\) is the set of successful iterations corresponding to positive-curvature steps up to iteration \(k\).
  By~\eqref{eq:succ-enum}, \(|\calS_{> 0, k}| = i + 1\).
  Thus,
  \[
    k_{i+1} - k_i - 1 \leq \log_{\gamma_1}
    \left[
      \frac{
        \tilde{b}_{\text{succ}} (1 + \mu (1 + (i + 1)^p))
      }{
        \gamma_3 \min(\sigma_{k_i}, \epsilon)
      }
      \right].
  \]
  We sum the above inequalities to bound the total number of unsuccessful iterations between iterations \(k_0\) and \(k_t\):
  \[
    \sum_{i = 0}^{t - 1} k_{i+1} - k_i - 1 \leq
    \sum_{i = 0}^{t - 1} \log_{\gamma_1} \left[
      \frac{
        1
      }{
        \gamma_3 \min(\sigma_{k_i}, \epsilon)
      }
      \right]
    + \sum_{i = 0}^{t - 1} \log_{\gamma_1} [\tilde{b}_{\text{succ}} (1 + \mu (1 + (i + 1)^p))].
  \]
  In the last term on the right-hand side, \(i + 1 \leq t = |\calS(\epsilon)| - 1\).
  Thus,
  \begin{multline*}
    \sum_{i = 0}^{t - 1} k_{i+1} - k_i - 1 \leq
    \sum_{i = 0}^{t - 1} \log_{\gamma_1} \left[
      \frac{
        1
      }{
        \gamma_3 \min(\sigma_{k_i}, \epsilon)
      }
      \right]
    + \\
    (|\calS(\epsilon)| - 1) \log_{\gamma_1} [\tilde{b}_{\text{succ}} (1 + \mu (1 + (|\calS(\epsilon)| - 1)^p))].
  \end{multline*}

  There remains to bound the number of unsuccessful iterations before \(k_0\).
  Let \(k \in \{0, \ldots, k_0 - 1\}\).
  The update of \(\sigma_k\) on unsuccessful iterations now yields \(\sigma_k \geq \sigma_0 \gamma_1^k\).
  Again, \Cref{lem:bsucc} and \Cref{asm:growth-Bk} yield
  \[
    \sigma_k \leq \tilde{b}_{\text{succ}} (1 + \mu (1 + |\calS_{> 0, k}|^p)) = \tilde{b}_{\text{succ}} (1 + \mu),
  \]
  because \(\calS_{> 0, k} = \varnothing\).
  Thus, with \(k = k_0 - 1\), we obtain
  \[
    k_0 \leq 1 + \log_{\gamma_1} \left(
    \frac{
      \tilde{b}_{\text{succ}} (1 + \mu)
    }{
      \sigma_0
    }
    \right).
  \]
  Combining the two bounds gives the desired result.
\end{proof}

We may expect that, for \(\epsilon\) small enough, \(\min(\sigma_{k_i}, \epsilon) = \epsilon\) for \(i = 0, \ldots, t-1\).
Hence, the bound given by \Cref{thm:unsuccessful-fast} is worse than that of \citep[Theorem~\(6.5\)]{diouane-habiboullah-orban-2024b} by an asymptotic factor of the form \(|\calS(\epsilon)| \log_{\gamma_1}(1 / (\gamma_3 \epsilon))\), which, although it deteriorates the worst case, only affects unsuccessful iterations, and remains modest in practice.
We stress again that our interest here is to minimally alter \Cref{alg:r2n} to preserve favorable complexity bounds and guarantee fast local convergence.

We now examine the behavior of \Cref{alg:r2n-fast} near an isolated minimizer.

\begin{assumption}%
  \label{asm:local}
  \(f\) is \(\mathcal{C}^2\) and \(\{x_k\} \to x^\star\) where \(\nabla f(x_\star) = 0\) and \(\nabla^2 f(x_\star) \succ 0\).
\end{assumption}

\Cref{asm:local} and continuity ensure that \(s_k \to 0\), \(\nabla f(x_k) \to 0\), \(\nabla^2 f(x_k) \succ 0\) for all large enough \(k\), and that there exists \(M > 0\) such that \(\|\nabla^2 f(x_k)\| \leq M\) for all \(k\).

To quantify the local convergence rate of \Cref{alg:r2n-fast}, we assume Lipschitz continuity of \(\nabla^2 f\) instead of \Cref{asm:hess-remainder}, and that the bound of \Cref{asm:hess-error} holds at every iteration instead of only for \(k \in \calS_{\leq 0}\).

\begin{problemassumption}%
  \label{asm:hess-lipschitz}
  \(\nabla^2 f\) is globally Lipschitz-continuous, i.e., there exists \(L_H \geq 0\) such that \(\|\nabla^2 f(x) - \nabla^2 f(y)\| \leq L_H \|x - y\|\) for all \(x\), \(y \in \R^n\).
\end{problemassumption}

\begin{modelassumption}%
  \label{asm:hess-error-all}
  There exist \(C_H > 0\) and \(0 < \zeta \leq 1\) such that, for all \(k \in \N\), \(\|(\nabla^2 f(x_k) - B_k) s_k\| \leq C_H \|s_k\|^{1 + \zeta}\), where \(x_k\) and \(s_k\) are the iterate and the trial step generated at iteration \(k\) of \Cref{alg:r2n-fast}.
\end{modelassumption}

To quantify the local convergence rate, we must specify how accurately~\eqref{eq:kkt-model} is solved at each iteration.
Define
\begin{equation}%
  \label{eq:def-r}
  r_k \eqdef (B_k + \sigma_k I) s_k + \nabla f(x_k),
\end{equation}
the residual at iteration \(k\).

\begin{lemma}
  Let \Cref{asm:hess-error-all} and \Cref{asm:local} be satisfied.
  If \(\{r_k\} \to 0\), \(\sigma_k \|s_k\| = O(\|r_k\| + \|\nabla f(x_k)\| + \|s_k\|)\).
  In particular, \(\{\sigma_k s_k\} \to 0\).
\end{lemma}

\begin{proof}
  The triangle inequality,~\eqref{eq:def-r} and \Cref{asm:hess-error-all} yield
  \begin{align*}
    \|\sigma_k s_k\| & = \|r_k - \nabla f(x_k) - B_k s_k\|
    \\
                     & = \|r_k - \nabla f(x_k) - \nabla^2 f(x_k) s_k - (B_k - \nabla^2 f(x_k)) s_k\|
    \\
                     & \leq \|r_k\| + \|\nabla f(x_k)\| + M \|s_k\| + C_H \|s_k\|^{1 + \zeta}.
  \end{align*}
  For all large enough \(k\), \(C_H \|s_k\|^{1 + \zeta} \leq M \|s_k\|\), which concludes the proof.
\end{proof}

\begin{lemma}%
  \label{lem:g=Ors}
  Let \Cref{asm:hess-error-all} and \Cref{asm:local} be satisfied.
  If \(\{r_k\} \to 0\), \(\|\nabla f(x_k)\| = O(\|r_k\| + \sigma_k \|s_k\| + \|s_k\|)\).
\end{lemma}

\begin{proof}
  For all large enough \(k\), \(C_H \|s_k\|^{1 + \zeta} \leq M \|s_k\|\).
  Thus, the triangle inequality,~\eqref{eq:def-r} and \Cref{asm:hess-error-all} ensure that for all large enough \(k\),
  \begin{align*}
    \|r_k - \nabla f(x_k)\| & = \|\nabla^2 f(x_k) s_k + (B_k - \nabla^2 f(x_k)) s_k + \sigma_k s_k\|
    \\
                            & \leq M \|s_k\| + C_H \|s_k\|^{1 + \zeta} + \sigma_k \|s_k\|
    \\
                            & \leq 2 M \|s_k\| + \sigma_k \|s_k\|.
  \end{align*}

  The triangle inequality ensures that \(\|r_k - \nabla f(x_k)\| \geq \|\nabla f(x_k)\| - \|r_k\|\).
  Therefore, \(\|\nabla f(x_k)\| \leq \|r_k\| + 2 M \|s_k\| + \sigma_k \|s_k\|\).
\end{proof}

\begin{lemma}%
  \label{lem:grad-next}
  Let \Cref{asm:hess-lipschitz,asm:hess-error-all} and \Cref{asm:local} be satisfied.
  If \(\{r_k\} \to 0\), \(\|\nabla f(x_k + s_k)\| = O(\|r_k\| + \|s_k\|^{1 + \zeta} + \sigma_k \|s_k\|)\).
\end{lemma}

\begin{proof}
  \Cref{asm:hess-lipschitz,asm:hess-error-all} and~\eqref{eq:def-r} yield
  \begin{align*}
    \nabla f(x_k + s_k) & = \nabla f(x_k) + \int_0^1 \nabla^2 f(x_k + \tau s_k) s_k \, \mathrm{d}\tau
    \\
                        & = \int_0^1 \nabla^2 f(x_k + \tau s_k) s_k \, \mathrm{d}\tau + r_k - \nabla^2 f(x_k) s_k - (B_k - \nabla^2 f(x_k)) s_k - \sigma_k s_k
    \\
                        & = \int_0^1 (\nabla^2 f(x_k + \tau s_k) - \nabla^2 f(x_k)) s_k \, \mathrm{d}\tau + r_k - (B_k - \nabla^2 f(x_k)) s_k - \sigma_k s_k.
  \end{align*}
  The triangle inequality and \Cref{asm:hess-lipschitz} give
  \begin{align*}
    \|\nabla f(x_k + s_k)\| & \leq \int_0^1 L_H \tau \|s_k\|^2 \, \mathrm{d}\tau + \|r_k\| + C_H \|s_k\|^{1 + \zeta} + \sigma_k \|s_k\|
    \\
                            & = \tfrac{1}{2} L_H \|s_k\|^2 + \|r_k\| + C_H \|s_k\|^{1 + \zeta} + \sigma_k \|s_k\|.
  \end{align*}
  Because \(0 < \zeta \leq 1\) and \(s_k \to 0\), \(\tfrac{1}{2} L_H \|s_k\|^2 = O(\|s_k\|^{1 + \zeta})\).
  The result follows.
\end{proof}

\begin{lemma}%
  \label{lem:s=Og}
  Let \Cref{asm:hess-error-all} and \Cref{asm:local} be satisfied.
  Assume that \(\|r_k\| = O(\|\nabla f(x_k)\|)\).
  Then, \(\|s_k\| = O(\|\nabla f(x_k)\|)\) and \(\sigma_k \|s_k\| = O(\|\nabla f(x_k)\|)\).
\end{lemma}

\begin{proof}
  Let \(\lambda_{\min} > 0\) be the smallest eigenvalue of \(\nabla^2 f(x_\star)\).
  \Cref{asm:local} implies that for all large enough \(k\), the smallest eigenvalue of \(\nabla^2 f(x_k)\) is at least \(\tfrac{1}{2} \lambda_{\min} > 0\).
  Thus, for all large enough \(k\), \Cref{asm:hess-error-all},~\eqref{eq:def-r} and the triangle inequality yield
  \begin{align*}
    \|r_k - \nabla f(x_k)\| & = \|\nabla^2 f(x_k) s_k + (B_k - \nabla^2 f(x_k)) s_k + \sigma_k s_k\|
    \\
                            & \geq \|(\nabla^2 f(x_k) + \sigma_k I) s_k\| - \|(B_k - \nabla^2 f(x_k)) s_k\|
    \\
                            & \geq (\tfrac{1}{2} \lambda_{\min} + \sigma_k) \|s_k\| - C_H \|s_k\|^{1 + \zeta}.
  \end{align*}
  Because \(s_k \to 0\), for all large enough \(k\), \(C_H \|s_k\|^{1 + \zeta} \leq \tfrac{1}{4} \lambda_{\min} \|s_k\|\).
  Substituting this into the inequality yields:
  \begin{equation}%
    \label{eq:smg}
    \|r_k - \nabla f(x_k)\| \geq (\tfrac{1}{4} \lambda_{\min} + \sigma_k) \|s_k\|.
  \end{equation}
  The triangle inequality and our assumption imply \(\|r_k - \nabla f(x_k)\| \leq \|r_k\| + \|\nabla f(x_k)\|\).

  On the one hand,~\eqref{eq:smg} gives \(\tfrac{1}{4} \lambda_{\min} \|s_k\| \leq (\tfrac{1}{4} \lambda_{\min} + \sigma_k) \|s_k\| \leq \|r_k\| + \|\nabla f(x_k)\| = O(\|\nabla f(x_k)\|)\).
  Similarly,~\eqref{eq:smg} gives \(\sigma_k \|s_k\| \leq (\tfrac{1}{4} \lambda_{\min} + \sigma_k) \|s_k\| = O(\|\nabla f(x_k)\|)\).
\end{proof}
\begin{corollary}%
  \label{cor:gradfp}
  Let \Cref{asm:hess-lipschitz,asm:hess-error-all} and \Cref{asm:local} be satisfied.
  If \(\|r_k\| = O(\|\nabla f(x_k)\|^{1 + \zeta})\) and \(\sigma_k = O(\|\nabla f(x_k)\|^{\zeta})\), \(\|\nabla f(x_k + s_k)\| = O(\|\nabla f(x_k)\|^{1 + \zeta})\).
\end{corollary}

\begin{proof}
  The result follows directly from \Cref{lem:grad-next,lem:s=Og}.
\end{proof}

\begin{corollary}%
  \label{cor:g=Os}
  Let \Cref{asm:hess-error-all} and \Cref{asm:local} be satisfied.
  If \(\|r_k\| = O(\|s_k\|)\) and \(\sigma_k = O(1)\), \(\|\nabla f(x_k)\| = O(\|s_k\|)\).
  Thus, if \(\|r_k\| = O(\max(\|\nabla f(x_k)\|, \, \|s_k\|))\) and \(\sigma_k = O(1)\), \(\|\nabla f(x_k)\| = \Theta(\|s_k\|)\).
\end{corollary}

\begin{proof}
  The first part follows directly from \Cref{lem:g=Ors}.
  The second part follows from the first part and \Cref{lem:s=Og}.
\end{proof}

\begin{lemma}%
  \label{lem:eta2}
  Let \Cref{asm:hess-lipschitz}, \Cref{asm:hess-error-all} and \Cref{asm:local} be satisfied.
  If \(B_k s_k = \nabla^2 f(x_k) s_k\) for all \(k\), assume also that \(\|r_k\| = O(1)\).
  Otherwise, assume that \(\|r_k\| = O(\|s_k\|^{1 + \zeta})\) and \(\sigma_k = O(\|s_k\|^{\zeta})\).
  Then, \(\rho_k \geq \eta_2\) for all large enough \(k\).
\end{lemma}

\begin{proof}
  \Cref{asm:hess-lipschitz} ensures that
  \[
    |f(x_k + s_k) - f(x_k) - \nabla f(x_k)^T s_k - \tfrac{1}{2} s_k^T \nabla^2 f(x_k) s_k| \leq \tfrac{1}{3!} L_H \|s_k\|^3.
  \]
  Thus, \Cref{asm:hess-error-all} yields
  \begin{align*}
    f(x_k + s_k) - \varphi(s_k; x_k) & \leq \tfrac{1}{3!} L_H \|s_k\|^3 - \tfrac{1}{2} s_k (B_k - \nabla^2 f(x_k)) s_k
    \\
                                     & \leq \tfrac{1}{3!} L_H \|s_k\|^3 + \tfrac{1}{2} C_H \|s_k\|^{2 + \zeta}.
  \end{align*}

  Let \(\lambda_{\min} > 0\) be the smallest eigenvalue of \(\nabla^2 f(x_\star)\).
  \Cref{asm:local} ensures that \(\nabla^2 f(x_k) \succ 0\) for all large enough \(k\), and, in particular, that its smallest eigenvalue is larger than \(\tfrac{1}{2} \lambda_{\min}\).
  We take the inner product of~\eqref{eq:def-r} with \(s_k\) and obtain
  \begin{align*}
    \varphi(0; x_k) - \varphi(s_k; x_k) & = -\nabla f(x_k)^T s_k - \tfrac{1}{2} s_k^T B_k s_k
    \\
                                        & = \tfrac{1}{2} s_k^T B_k s_k + \sigma_k \|s_k\|^2 - r_k^T s_k
    \\
                                        & = \tfrac{1}{2} s_k^T \nabla^2 f(x_k) s_k + \tfrac{1}{2} s_k^T (B_k - \nabla^2 f(x_k)) s_k + \sigma_k \|s_k\|^2 - r_k^T s_k
    \\
                                        & \geq \tfrac{1}{4} \lambda_{\min} \|s_k\|^2 -\tfrac{1}{2} \|s_k\| \|(B_k - \nabla^2 f(x_k)) s_k\| + \sigma_k \|s_k\|^2 - \|r_k\| \|s_k\|
    \\
                                        & \geq \tfrac{1}{4} \lambda_{\min} \|s_k\|^2 -\tfrac{1}{2} C_H \|s_k\|^{2 + \zeta} + \sigma_k \|s_k\|^2 - \|r_k\| \, \|s_k\|.
  \end{align*}
  The above combine to yield
  \begin{align*}
    \rho_k & = 1 + \frac{ \varphi(s_k; x_k) - f(x_k + s_k) }{ \varphi(0; x_k) - \varphi(s_k; x_k) }
    \\
           & \geq 1 -  \frac{ \tfrac{1}{3!} L_H \|s_k\|^3 + \tfrac{1}{2} C_H \|s_k\|^{2 + \zeta} }{ \tfrac{1}{4} \lambda_{\min} \|s_k\|^2 - \tfrac{1}{2} C_H \|s_k\|^{2 + \zeta} - \|r_k\| \, \|s_k\| + \sigma_k \|s_k\|^2 }.
  \end{align*}
  Our assumptions ensure that the numerator above is \(\Theta(\|s_k\|^{2 + \zeta})\) while the denominator is \(\Theta(\|s_k\|^2)\).
  Therefore, \(\rho_k \geq \eta_2\) for all large enough \(k\).

  If \(B_k s_k = \nabla^2 f(x_k) s_k\) for all \(k\), the terms involving \(\zeta\) above are simply absent, i.e.,
  \(
  f(x_k + s_k) - \varphi(s_k; x_k) \leq \tfrac{1}{3!} L_H \|s_k\|^3,
  \)
  and
  \(
  \varphi(0; x_k) - \varphi(s_k; x_k) > (\tfrac{1}{4} \lambda_{\min} + \sigma_k) \|s_k\|^2 - \|r_k\| \, \|s_k\|.
  \)
  Hence, \(\rho_k \to 1\), and the result continues to hold.
\end{proof}

We compile the conclusions of this section into the following result.

\begin{theorem}%
  \label{thm:r2n-fast}
  Let \Cref{asm:hess-lipschitz}, \Cref{asm:hess-error-all} and \Cref{asm:local} be satisfied.
  For all large enough \(k\), \(\rho_k \geq \eta_2\) and every iteration is very successful.
  If \(\|r_k\| = O(\|\nabla f(x_k)\|^{1 + \zeta}\|)\) and \(\sigma_k = O(\|\nabla f(x_k)\|^{\zeta})\), \(\|\nabla f(x_{k+1})\| = O(\|\nabla f(x_k)\|^{1 + \zeta})\) and \(\|x_{k+1} - x_\star\| = O(\|x_k - x_\star\|^{1 + \zeta})\).
  In other words, the sequence \(\{x_k\}\) generated by \Cref{alg:r2n-fast} converges locally superlinearly with rate \(1 + \zeta\) to an isolated minimizer.
  If \(B_k s_k = \nabla^2 f(x_k) s_k\) for all \(k\), we may choose \(\zeta = 1\), and the rate is quadratic.
\end{theorem}

\begin{proof}
  By \Cref{cor:g=Os}, \(\|\nabla f(x_k)\| = \Theta(\|s_k\|)\), so that our assumptions on \(r_k\) and \(\sigma_k\) imply those of \Cref{lem:eta2}, which yields the first assertion.
  The second follows from \Cref{cor:gradfp}.
  A Taylor expansion about \(x_\star\) yields, for all large enough \(k\),
  \[
    \|\nabla f(x_k)\| = \|\nabla^2 f(x_\star) (x_k - x_\star) + o(\|x_k - x_\star\|)\| = \Theta(\|x_k - x_\star\|).
  \]
  Thus, \(\|x_{k+1} - x_\star\| = O(\|x_k - x_\star\|^{1 + \zeta})\).
\end{proof}

In practice, the value of \(\zeta\) is often unknown.
However, it should be clear from the analysis above that if one sets \(\|r_k\| = O(\|\nabla f(x_k)\|^{1 + \mu})\) and \(\sigma_k = O(\|\nabla f(x_k)\|^{\omega})\) with \(0 < \mu\), \(\omega \leq 1\), the overall convergence rate will be \(1 + \min(\mu, \omega, \zeta)\).
Alternatively, one may simply set \(\|r_k\| = O(\|\nabla f(x_k)\|^2)\) and \(\sigma_k = O(\|\nabla f(x_k)\|)\), as in \Cref{alg:r2n-fast}, and the overall convergence rate will be \(1 + \zeta\).
Although the latter choice may require solving~\eqref{eq:kkt-model} to higher accuracy than strictly necessary, it dispenses with trying to estimate \(\zeta\).

\section{A Variant for Nonlinear Least-Squares}%
\label{sec:nlls}

\Cref{alg:r2n} specializes naturally to the smooth nonlinear least-squares problem
\begin{equation}%
  \label{eq:nlls}
  \minimize{x \in \R^n} \ f(x), \qquad f(x) \eqdef \tfrac{1}{2} \|F(x)\|^2,
\end{equation}
where $F : \R^n \to \R^m$ is continuously differentiable.
The gradient is given by $\nabla f(x) = J(x)^T F(x)$, where $J(x)$ is the Jacobian of $F$.
For simplicity, we write $F_k \eqdef F(x_k)$, $J_k \eqdef J(x_k)$, and $g_k \eqdef \nabla f(x_k) = J_k^T F_k$.

Following the classic framework of \citet{levenberg-1944} and \citet{marquardt-1963}, we use the Gauss-Newton model by selecting \(B_k \eqdef J_k^T J_k\).
Because $\sigma_k > 0$, $B_k + \sigma_k I \succ 0$, the subproblem is strongly convex, and hence admits a unique (global) minimizer.
In addition, nonpositive-curvature directions never occur, and neither does the Armijo-Goldstein line search.
The step $s_k$ can be viewed as an approximate solution of the linear least-squares model
\begin{equation}%
  \label{eq:m-r2nls}
  m(s; x_k, \sigma_k) \eqdef \tfrac{1}{2} \|J_k s + F_k\|^2 + \tfrac{1}{2} \sigma_k \|s\|^2 =
  \tfrac{1}{2}
  \left\|
  \begin{bmatrix}
    J_k \\
    \sqrt{\sigma_k} I
  \end{bmatrix}
  s +
  \begin{bmatrix}
    F_k \\
    0
  \end{bmatrix}
  \right\|^2.
\end{equation}

In this context,~\eqref{eq:kkt-model} is the system of normal equations associated with~\eqref{eq:m-r2nls}.
The solvers used to compute a step parallel those used in \Cref{alg:r2n}, except that there is no quasi-Newton solver, but they exploit the least-squares structure of \(m(\cdot; x_k, \sigma_k)\):
\begin{itemize}
  \item CGLS\@: a variant of CG described by \citet{hestenes-stiefel-1952} in which products with \(J_k\) and \(J_k^T\) are decoupled, and that is equivalent to CG applied to~\eqref{eq:kkt-model} in exact arithmetic, but that is typically more robust in practice.
  \item LSQR \citep{paige-saunders-1982}: a Krylov method that is equivalent to CG applied to~\eqref{eq:kkt-model} in exact arithmetic, but that is typically more robust in practice.
  \item LSMR \citep{fong-saunders-2011}: a Krylov method that is equivalent to MINRES applied to~\eqref{eq:kkt-model} in exact arithmetic, but that is typically more robust in practice.
  \item QRMumps \citep{agullo-buttari-guermouche-lopez-2016}: a multi-core sparse QR factorization used to factorize the matrix of the linear residual in~\eqref{eq:m-r2nls}.
    Indeed the \(R\) factor of the latter is the transposed Cholesky factor of \(B_k + \sigma_k I\).
    This approach parallels the sparse factorization of \(B_k + \sigma_k I\) in~\eqref{eq:kkt-model}.
\end{itemize}

Our implementation uses the Julia implementation of CGLS, LSQR and LSMR found in \citep{montoison-orban-2023} and the QRMumps interface found in \citep{buttari-montoison-qrmumps-2021}.
By extension of the subproblem solvers for \Cref{alg:r2n}, the solvers above guarantee Cauchy decrease under the same assumptions.

The resulting procedure is a globally convergent Levenberg-Marquardt method.
In~\eqref{eq:nlls}, denote \(F(x) = (f_1(x), \ldots, f_m(x))\).
If each \(f_i\) is twice differentiable, \Cref{asm:hess-error-all} can be stated as
\[
  \|(\sum_{i=1}^m f_i(x_k) \nabla^2 f_i(x_k)) s_k\| \leq C_H \|s_k\|^{1 + \zeta},
\]
and holds, among others, when \(F\) is linear near a solution, and when \(x_k\) approaches a zero-residual solution, but may potentially hold in other situations.
All global convergence, worst-case complexity bounds, and local convergence of \Cref{sec:convergence,sec:complexity,sec:local} carry over without modification under the assumption above.

\section{Implementation and Numerical Experiments}%
\label{sec:implementation}

We describe our implementation of R2N and R2NLS and assess it on problems from the CUTEst collection.
In \Cref{sec:benchmark-setup,sec:parameters,sec:metrics}, we describe the computational environment, test set, algorithmic parameters and performance measures.
Our first objective is to identify a best default configuration for R2N by evaluating options for the treatment of nonpositive curvature, the Krylov subproblem solver, the update of the regularization parameter, and direct vs.\ iterative subproblem solves.
We use exact second derivatives in \Cref{sec:r2n-exact-benchmarks}.
We end the section with a comparison with a trust-region method, which has natural treatment of nonpositive curvature directions.
We repeat similar experiments using limited-memory BFGS Hessian approximations in \Cref{sec:r2n-lbfgs-benchmarks}.
Finally, we consider nonlinear least-squares problems in \Cref{sec:r2nls-benchmarks}, for which nonpositive curvature cannot arise.

\subsection{Implementation and Benchmark Setup}%
\label{sec:benchmark-setup}

All experiments were conducted on a Linux system equipped with two AMD EPYC 7402 24-core processors and 1\,TiB of RAM, though only 8\,GiB are available to any individual user, running AlmaLinux 9.8 (Olive Jaguar).
The solvers are implemented in Julia 1.13, and every run uses a single CPU thread without GPU acceleration.
Julia relies on a BLAS library, while MA57 and QRMumps both rely on another, and both are set to use a single thread.
Problems are solved sequentially, so the reported times are free of contention between runs.
All computations are carried out in Float64.

The R2N and R2NLS solvers are part of the JSOSolvers.jl package \citep{migot-monnet-orban-siqueira-2026} version 0.14.9, a component of the JuliaSmoothOptimizers (JSO) organization.
Their implementation is based on several JSO libraries: NLPModels.jl \citep{orban-siqueira-nlpmodels-2023} (v0.21.12) defines abstract types to represent problems, NLPModelsModifiers.jl \citep{migot-orban-siqueira-nlpmodelsmodifiers-2023} (v0.8.0) is used to replace a model's exact second derivatives with quasi-Newton approximations, LinearOperators.jl \citep{leconte-orban-siqueira-linearoperators-2023} (v2.14.1) is used to model exact Hessians and limited-memory quasi-Newton approximations as abstract linear operators, Krylov.jl \citep{montoison-orban-2023} (v0.10.6) has implementations of the various Krylov methods used in R2N and R2NLS, HSL.jl \citep{fowkes-lister-montoison-orban-2024} (v0.5.2, LibHSL 2025.7.21) provides an interface to MA57, QRMumps.jl \citep{buttari-montoison-qrmumps-2021} (v0.3.2, QRMumps 3.1.1) interfaces QRMumps, and SolverCore.jl \citep{migot-orban-siqueira-solvercore-2023} (v0.3.10) provides solver building blocks.
Finally, CUTEst.jl (v1.3.7, CUTEst 2.6.0, SIFDecode 3.1.0) provides access to problems from the CUTEst collection using the abstractions of NLPModels.jl, and SolverBenchmark.jl \citep{orban-siqueira-solverbenchmark-2023} (v0.6.5) standardizes benchmarking activities.
Problems obtained from CUTEst.jl allow access to explicit Hessians as sparse matrices, or as abstract linear operators via operator-vector products.

While evaluating the model curvature along the gradient, \(\gamma_k\), requires only one operator-vector product and a few dot products, computing the exact Cauchy step \(\skcp\) (\Cref{alg:cauchy-step}) is not necessarily cheap.
When \(\gamma_k \leq 0\), it requires estimating the spectral norm \(\|B_k\|\) at \Cref{alg:cauchy-step:nuk}.
We use the estimator implemented by default in LinearOperators.jl, which computes the norm via an eigendecomposition for \(n \leq 5\) and otherwise calls the implicitly restarted Lanczos iteration of ARPACK \citep{lehoucq-sorensen-yang-1998} via the Arpack.jl (v0.5.4) package.
Thus, falling back to \(\skcp\) is a safe global convergence mechanism, but it is not free.
When a direction of nonpositive curvature is detected and a step length \(\alpha_k\) satisfying the Armijo-Goldstein conditions has been found, the resulting step is subjected to the usual ratio test by default.
\Cref{prop:ag-accept} shows that \(\rho_k \geq \tilde{\eta}_1 > 0\) for such steps, so the test may optionally be bypassed; we evaluate that variant in \Cref{sec:case-variants}.
Before factorizing, the direct subproblem solvers apply a density test: MA57 is not used if \(n > 12{,}000\) and \(\mathrm{nnz}(H) > n (n+1) / 4\), i.e., the Hessian is more than half full.
Similarly, QRMumps is not used if \(m n > 12{,}000^2\) and \(\mathrm{nnz}(J) > \tfrac{1}{2} m n\).
Below the size thresholds, even a fully dense matrix is inexpensive to factorize, and density alone is no reason to avoid a direct solve.
When the test fails, the solver declines to factorize and recommends a Krylov subproblem solver instead.

\paragraph{Test problems.}
Our benchmark uses problems drawn from the CUTEst test set \citep{gould-orban-toint-2015}.
For the general unconstrained experiments of \Cref{sec:r2n-exact-benchmarks,sec:r2n-lbfgs-benchmarks}, we retain the smooth unconstrained problems, excluding those of nonlinear least-squares type, which are treated separately in \Cref{sec:r2nls-benchmarks}.
No bound is imposed on the number of variables \(n\), and the suite spans problems from a few variables to large-scale instances.

Problems FLETCBV2, FLETCBV3, FLETCHBV, LRCOVTYPE, and NONMSQRT were excluded because they exceed the time limit of 600s in all solver configurations.
In addition, PENALTY1 and VAREIGVL were excluded because they fail the density test for MA57.
The remaining suite contains \(87\) problems.

In the nonlinear least-squares experiments of \Cref{sec:r2nls-benchmarks}, we use the CUTEst least-squares problems with at most \(5,000\) variables.
The seven problems ARGLCLE, DIAMON3D, DMN15103, DMN37142, DMN37143, KSS, and PATTERNNE fail the density test, and are excluded.
The remaining suite contains \(177\) problems.

\subsection{Algorithm Parameters and Tolerances}%
\label{sec:parameters}

All global tolerances, experimental settings, and algorithmic constants for the step-acceptance logic and regularization multipliers are summarized in \Cref{tab:constants}.
Termination is controlled by the wall-clock time limit and by the gradient-norm criterion \(\|\nabla f(x_k)\| \le \epsilon_a + \epsilon_r \|\nabla f(x_0)\|\).
For least-squares problems, the solvers also stop when \(\|F(x_k)\| \le \epsilon_{Fa} + \epsilon_{Fr} \|F(x_0)\|\).
The same criteria, tolerances and time limit are used for every solver compared in \Cref{sec:implementation}, including the trust-region solver Trunk of \Cref{sec:case-cg-vs-trunk,sec:case-trunk-vs-lbfgs-cg,sec:r2nls-benchmarks}.
Apart from these settings, Trunk is used as distributed in JSOSolvers.jl with its default parameters; in particular, it backtracks along rejected steps and selects its own initial trust-region radius.

We set the initial regularization parameter to
\begin{equation}
  \label{eq:sigma0}
  \sigma_0 \eqdef \max\!\left( \sigma_{\min}, \; \frac{2^{\mathrm{round}(\log_2(\|\nabla f(x_0)\| + 1))}}{\|\nabla f(x_0)\|} \right),
\end{equation}
and \(\sigma_0 \eqdef \max(\sigma_{\min}, 1)\) in the case \(\nabla f(x_0) = 0\).
Rounding to the nearest power of two perturbs the numerator by at most a factor \(\sqrt{2}\), and therefore \(\sigma_0 \geq (1 + \|\nabla f(x_0)\|^{-1}) / \sqrt{2} > 2^{-1/2}\).

Hence, the lower bound in~\eqref{eq:sigma0} is typically never active and \(\sigma_{\min}\) serves only as a floor on subsequent decreases of \(\sigma_k\).
This choice requires no per-problem tuning, and it keeps the problem-dependent constants \(b_{\max}\) of \Cref{prop:succ-pc} and the leading term of \Cref{thm:unsuccessful-fast}, both of which involve \(\sigma_0\), of moderate size.
We set \(\sigma_{\min} = \epsilon_M \approx 2.2 \times 10^{-16}\).
A strictly positive \(\sigma_{\min}\) satisfies assumption~(i) of \Cref{prop:npc-dir-cauchy}; at this magnitude it does not interfere numerically with the local convergence analysis of \Cref{sec:local}, which requires \(\sigma_k \to 0\) near an isolated minimizer.

We set \(\theta_2 = \epsilon_M^{-1}\), so that the safeguard on \Cref{alg:r2n:step-comparison} of \Cref{alg:r2n} is effectively disabled, and the assumption following \Cref{thm:decrease-npc} holds.

When a nonpositive curvature direction is detected, the default handler is the Armijo-Goldstein line search with the constants \(\rho_L\) and \(\rho_U\) of \Cref{tab:constants}.
It starts from \(\alpha = 1\), multiplies the trial step length by \(1.5\) while the Goldstein condition fails, and halves it while the Armijo condition fails, with at most \(10\) expansions and \(10\) reductions.
If either budget is exhausted, or if the direction returned by the subproblem solver is not a descent direction, the step is rejected and \(\sigma_k\) is increased.
A step whose predicted reduction does not exceed \(\epsilon_M \max(1, |f(x_k)|)\) is also rejected, without evaluating \(f\).
When using the MA57 direct solver, nonpositive curvature is detected as a by-product of the factorization itself, which reports the inertia.
If detected, the step is rejected and \(\sigma_k\) is increased.
Indeed, the Armijo-Goldstein line search is not applicable because an exact solution of~\eqref{eq:kkt-model} that is a nonpositive-curvature direction cannot be a descent direction.
For Krylov-based inner solvers, the absolute stopping tolerance matches \(\epsilon_a\), while the relative tolerance \(\epsilon_{rk}\) follows the forcing sequence
\begin{equation}
  \label{eq:forcing}
  \epsilon_{rk} = \max\!\left(\epsilon_r,\,
  \min\!\left(0.1,\, \sqrt{\|\nabla f(x_k)\|},\,
  0.9\,\epsilon_{r,k-1}\right)\right),
\end{equation}
updated after each accepted step.
The Krylov subproblem solvers of Trunk use the same rule, but update it at every iteration, including after a rejected step, in which case their tolerance decreases by a factor \(0.9\) while that of R2N is unchanged.
In the limited-memory experiments of \Cref{sec:r2n-lbfgs-benchmarks}, \(B_k\) is the L-BFGS approximation implemented in LinearOperators.jl with memory \(q = 5\), initialized with \(B_{k,0} = (y_k^T y_k / s_k^T y_k) I\) from the most recent stored pair, and an update is performed whenever \(s_k^T y_k \geq \epsilon_M\).
The same model and memory are used on both sides of the comparison of \Cref{sec:case-trunk-vs-lbfgs-cg}.
R2NLS uses the update of \(\sigma_k\) given in \Cref{alg:r2n}; the local analysis of \Cref{sec:local} does not apply to least-squares problems in general.
An objective is declared unbounded when \(f(x_k) < \min(-1, f(x_0)) / \epsilon_M\).

\begin{table}[ht]
  \centering
  \caption{Experimental settings and algorithmic constants.}%
  \label{tab:constants}
  \begin{tabular}{llcl}
    \hline\noalign{\smallskip}
    \textbf{Parameter} & \textbf{Symbol} & \textbf{Value} & \textbf{Context / Condition} \\
    \noalign{\smallskip}\hline\noalign{\smallskip}
    \multicolumn{4}{l}{\textit{Global tolerances and settings}} \\
    Wall-clock time limit & \(t_{\max}\) & \(600\) s & All solvers \\
    Relative stopping tolerance & \(\epsilon_r\) & [\(10^{-6}\)] & All solvers \\
    Absolute stopping tolerance & \(\epsilon_a\) & [\(10^{-6}\)] & All solvers \\
    Residual stopping tolerances & \(\epsilon_{Fa}, \epsilon_{Fr}\) &
    [\(\sqrt{\epsilon_M}, \epsilon_M\)] & Least-squares solvers \\
    Density test size threshold & --- & \(n > 12{,}000\) & \cref{sec:benchmark-setup} \\
    (MA57) & & & \\
    Density test size threshold & --- & \(m n > 12{,}000^2\) & \cref{sec:benchmark-setup} \\
    (QRMumps) & & & \\
    \noalign{\smallskip}
    \multicolumn{4}{l}{\textit{Default constants (\Cref{alg:r2n})}} \\
    Cauchy step parameter & \(\theta_1\) & \(0.9\) &\eqref{eq:def-nuk} \\
    Maximum step ratio bound & \(\theta_2\) & \(\epsilon_M^{-1}\) &
    \Cref{alg:r2n:step-comparison} \\
    Accept step reduction ratio & \(\eta_1\) & \(\epsilon_M^{1/4}\) &
    \Cref{alg:r2n:step-accept} \\
    Very successful step ratio & \(\eta_2\) & \(0.95\) &
    \Cref{alg:r2n:step-update} \\
    Reg.\ increase (successful step) & \(\gamma_1\) & \(1.5\) &
    \Cref{alg:r2n:step-update} \\
    Reg.\ increase (rejected step) & \(\gamma_2\) & \(2.5\) &
    \Cref{alg:r2n:step-update} \\
    Reg.\ decrease (very successful step) & \(\gamma_3\) & \(0.5\) &
    \Cref{alg:r2n:step-update} \\
    Initial regularization & \(\sigma_0\) & scale-adaptive &\eqref{eq:sigma0} \\
    Minimum regularization & \(\sigma_{\min}\) & \(\epsilon_M\) &
    Floor on \(\sigma_k\) \\
    Nonpositive-curvature handler & --- & Armijo-Goldstein & \\
    Armijo parameter & \(\rho_L\) & \(10^{-4}\) &\eqref{eq:ag} \\
    Goldstein parameter & \(\rho_U\) & \(1 - 10^{-4}\) &\eqref{eq:ag} \\
    Initial step length & --- & \(1\) &\eqref{eq:ag} \\
    Expansion factor & --- & \(1.5\) & Max.\ \(10\) expansions \\
    Reduction factor & --- & \(0.5\) & Max.\ \(10\) reductions \\
    \noalign{\smallskip}
    \multicolumn{4}{l}{\textit{Variant parameters (\Cref{alg:r2n-fast}, \cref{sec:case-variants})}} \\
    Fast update & --- & on / off & Replaces the \(\sigma_k\) update \\
    & & & when \(\rho_k \geq \eta_2\) \\
    Skip \(\rho\) & --- & on / off & Off by default; \(\sigma_k\) \\
    & & & unchanged after AG step \\
    L-BFGS memory & \(q\) & \(5\) & \Cref{sec:r2n-lbfgs-benchmarks} \\
    \noalign{\smallskip}\hline
  \end{tabular}
\end{table}

\subsection{Performance Metrics}%
\label{sec:metrics}

We compare solvers by means of \citet{dolan-more-2002} performance profiles, produced by SolverBenchmark.jl \citep{orban-siqueira-solverbenchmark-2023}.
For each problem \(p\) and solver \(s\), let \(t_{p,s} > 0\) denote the measure of interest, with \(t_{p,s} = \infty\) if \(s\) fails on \(p\).
The performance ratio is \(r_{p,s} \eqdef t_{p,s} / \min_{s'} t_{p,s'}\), and the profile of solver \(s\) gives, as a function of \(\tau \geq 1\), the fraction of problems for which \(r_{p,s} \leq \tau\).
The measures are wall-clock time, objective evaluations, gradient evaluations and equivalent Hessian evaluations \(n_H + n_{Hv} / n\), which counts \(n\) Hessian-vector products as one Hessian evaluation.
Because the gradient is evaluated once per accepted step, the gradient profile also measures the number of successful iterations.
For least-squares problems, we report residual evaluations, iterations and equivalent Jacobian evaluations \(n_J + n_{Jv} / n + n_{J^T v} / m\).
\Cref{tab:solved} reports the number of problems solved by each configuration and how the failed runs terminated.
A run is counted under Exception when the solver stops on an error before meeting a stopping criterion.
For instance, on INDEF, whose objective is unbounded below, such a termination is expected.
With more than two solvers, a profile identifies the solver most likely to be within a factor \(\tau\) of the best, but does not compare them two-by-two \citep{gould-scott-2016}.
Below, we do not report results with R2NLS in full detail because it is an implementation of the well-known \citet{levenberg-1944} and \citet{marquardt-1963} method, as is its trust-region variant, first described by \citet{more-1978}.
We use our own implementation of the latter, named TrunkNLS, with LSMR as subproblem solver.

\begin{table}[htbp]
  \centering
  \caption{Problems solved, and reason for termination otherwise, for each configuration appearing in \Cref{sec:r2n-exact-benchmarks,sec:r2n-lbfgs-benchmarks,sec:r2nls-benchmarks}.
  The unconstrained suite contains \(87\) problems and the least-squares suite \(177\).
  Exception (non-descent): termination on a search direction that is not a descent direction.}%
  \label{tab:solved}
  \footnotesize
  \setlength{\tabcolsep}{2.5pt}%
  \begin{tabular}{llcccc}
    \toprule
    & \textbf{Configuration} & \textbf{Solved} & \textbf{Timed out} & \textbf{Unbounded} & \textbf{Exception} \\
    \midrule
    \multicolumn{6}{l}{\textit{Exact Hessian, 87 problems}} \\
    \cref{sec:case-npc}
      & R2N / MINRES-QLP / AG                       & 85 & 1 & 1 & 0 \\
      & R2N / MINRES-QLP / Prev                     & 85 & 1 & 1 & 0 \\
      & R2N / MINRES-QLP / Sigma                    & 85 & 1 & 1 & 0 \\
      & R2N / MINRES-QLP / CP                       & 80 & 6 & 1 & 0 \\
    \cref{sec:case-cg-vs-minresqlp}
      & R2N / CG / AG                               & 84 & 2 & 1 & 0 \\
    \cref{sec:case-variants}
      & R2N / MINRES-QLP / AG / Fast                & 85 & 1 & 1 & 0 \\
      & R2N / MINRES-QLP / AG / Skip \(\rho\)        & 84 & 2 & 1 & 0 \\
      & R2N / MINRES-QLP / AG / Fast / Skip \(\rho\) & 84 & 2 & 1 & 0 \\
    \cref{sec:case-ma57}
      & R2N / MA57 / Sigma                          & 85 & 0 & 2 & 0 \\
      & R2N / MA57 / Sigma / Fast                   & 85 & 0 & 2 & 0 \\
    \cref{sec:case-cg-vs-trunk}
      & R2N / CG / AG / Fast                        & 85 & 1 & 1 & 0 \\
      & Trunk / CG                                  & 82 & 4 & 1 & 0 \\
    \midrule
    \multicolumn{6}{l}{\textit{L-BFGS model, 87 problems}} \\
    \cref{sec:case-lbfgs-hessian}
      & R2N / Shifted / Fast                        & 76 & 8 & 3 & 0 \\
      & R2N / MINRES-QLP / Fast                     & 78 & 6 & 3 & 0 \\
    \cref{sec:case-trunk-vs-lbfgs-cg}
      & R2N / CG / Fast                             & 78 & 6 & 3 & 0 \\
      & Trunk / CG                                  & 78 & 7 & 1 & 1 \\
    \midrule
    \multicolumn{6}{l}{\textit{Nonlinear least squares, 177 problems}} \\
    \cref{sec:r2nls-benchmarks}
      & R2NLS / LSMR                                & 169 & 8 & --- & 0 \\
      & R2NLS / QRMumps                           & 170 & 7 & --- & 0 \\
    \cref{sec:r2nls-benchmarks}
      & TrunkNLS / LSMR                             & 174 & 2 & --- & 1 \\
    \bottomrule
  \end{tabular}
\end{table}

\subsection{R2N with Exact Second Derivatives}%
\label{sec:r2n-exact-benchmarks}

We first evaluate R2N when exact second derivatives are available, either through Hessian-vector products or as an assembled sparse Hessian.
We progressively compare options to build an efficient and robust default configuration.
We examine in turn the treatment of nonpositive-curvature directions (\Cref{sec:case-npc}), the Krylov subproblem solver (\Cref{sec:case-cg-vs-minresqlp}), the regularization parameter update, and compare options combined with an iterative, and with a direct subproblem solver (\Cref{sec:case-variants,sec:case-ma57}).
We then compare direct vs.\ iterative solves (\Cref{sec:case-exact-vs-krylov}), and finally compare our best configuration of R2N with a trust-region solver (\Cref{sec:case-cg-vs-trunk}).

\subsubsection{Nonpositive-Curvature Handlers}%
\label{sec:case-npc}

Our first experiment concerns the treatment of nonpositive-curvature (NPC) directions.
We fix MINRES-QLP as subproblem solver and vary only the handler; the model, the regularization update and the stopping tolerances are identical across the four runs.
The four handlers are those discussed in \Cref{sec:npc-treatment}:
\begin{itemize}
  \item \textbf{AG:} an Armijo-Goldstein line search along the NPC direction, i.e.,~\eqref{eq:ag}; 
  \item \textbf{Prev:} the most recent Krylov iterate preceding the detection of the NPC direction, which continues to satisfy sufficient decrease; when NPC appears at the first inner iteration, no previous iterate exists, and the subproblem solver returns the right-hand side \(-\nabla f(x_k)\), which Prev uses as the step;
  \item \textbf{Sigma:} an increase of \(\sigma_k\) followed by a restart of the subproblem solver, which promotes convexity of the model at the price of the work already performed;
  \item \textbf{CP:} a fallback on the Cauchy step, which ignores curvature information and requires an estimate of \(\|B_k\|\) whenever \(\gamma_k \leq 0\) (\Cref{alg:cauchy-step:nuk} of \Cref{alg:cauchy-step}).
\end{itemize}

\begin{figure}[ht]
  \centering
  \includetikzgraphics[width=\linewidth]{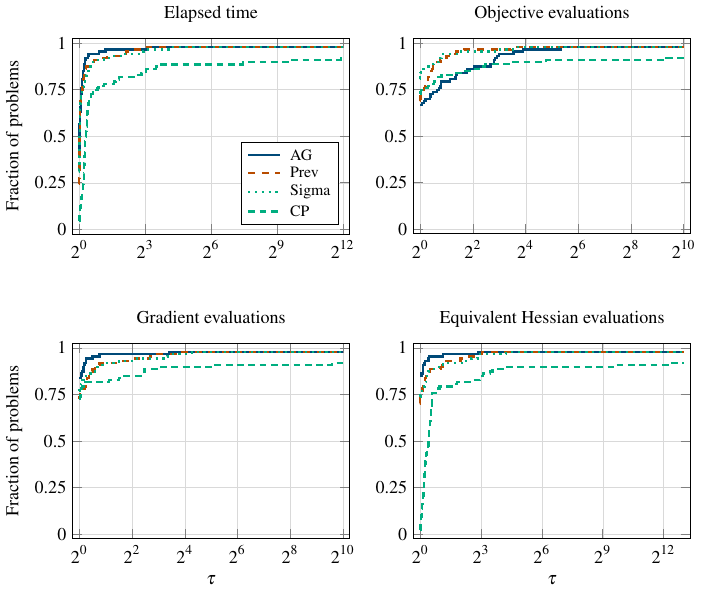}
  \caption{%
    Performance profiles of R2N with exact second derivatives, the MINRES-QLP subproblem solver and the \(\sigma_k\) update of \Cref{alg:r2n}.
    Comparison between nonpositive-curvature (NPC) handlers: Armijo-Goldstein line search (AG), previous MINRES-QLP iterate (Prev), increase of \(\sigma_k\) (Sigma), and fallback on the Cauchy step (CP).
  }%
  \label{fig:r2n-npc}
\end{figure}

\Cref{fig:r2n-npc} shows that CP is dominated by the three other handlers on elapsed time and on equivalent Hessian evaluations.
It is also the least robust: it solves five fewer problems than AG, Prev and Sigma, which solve the same \(85\) problems (\Cref{tab:solved}).
The first-order Cauchy step increases the number of outer iterations, and on those iterations where \(\gamma_k \leq 0\) it also requires an estimate of \(\|B_k\|\), obtained as described in \Cref{sec:benchmark-setup}; both effects are likely to contribute to its cost.
Falling back on the Cauchy step is a sound global convergence mechanism, but an expensive default.

The three remaining handlers are close, and their profiles largely overlap.
AG has the highest profiles for gradient evaluations and equivalent Hessian evaluations, but the lowest for objective evaluations at small \(\tau\), because each Armijo-Goldstein line search evaluates \(f\) several times; the elapsed-time profiles overlap.
AG thus trades objective evaluations for gradient evaluations and Hessian-vector products, and is preferable when the latter dominate the cost of an iteration, whereas Prev or Sigma may be preferable when objective evaluations are the more expensive.
We also ran AG with other values of \(\rho_L\) and \(\rho_U\) and of the expansion and reduction factors of the line search, and observed no noticeable change in the results, which we do not report.
We adopt AG as the default for the remainder of our experiments.

\subsubsection{Truncated CG versus MINRES-QLP}%
\label{sec:case-cg-vs-minresqlp}

With the handler fixed to AG, we compare the two Krylov methods available for the subproblem.
Truncated CG \citep{steihaug-1983} detects nonpositive curvature based on its search directions.
MINRES-QLP \citep{choi-paige-saunders-2011} coincides with MINRES in exact arithmetic on nonsingular systems, so that the guarantees of \Cref{sec:modified-minres} apply to it, and remains well defined on singular systems.
Equipped with the mechanism of \citet{liu-roosta-2022}, MINRES-QLP detects nonpositive curvature directly from quantities generated by the Lanczos process.
The outer loop, the handler and the stopping tolerances are identical.

\begin{figure}[ht]
  \centering
  \includetikzgraphics[width=\linewidth]{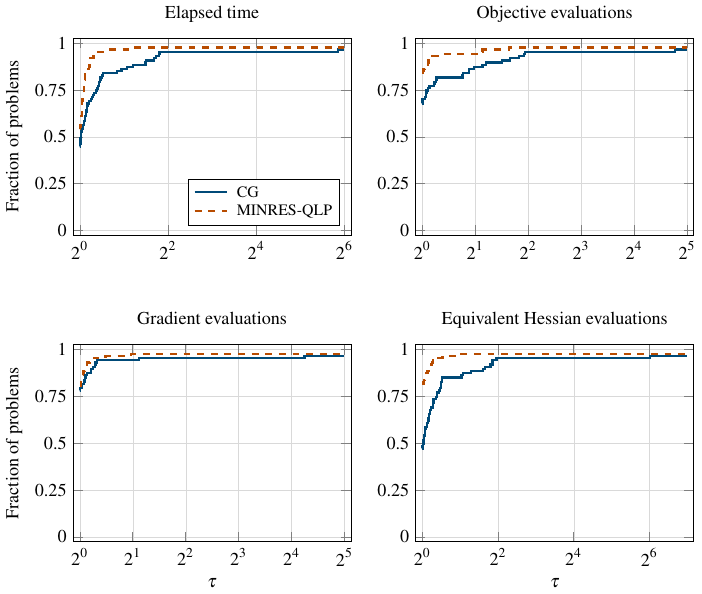}
  \caption{%
    Performance profiles of R2N with exact second derivatives, the AG handler and the base update of \(\sigma_k\).
    Comparison between subproblem solvers: truncated CG (CG) and MINRES-QLP (MINRES-QLP).
  }%
  \label{fig:r2n-cg-vs-minresqlp}
\end{figure}

The profiles of MINRES-QLP lie above those of truncated CG for all four measures of \Cref{fig:r2n-cg-vs-minresqlp}.
MINRES-QLP solves one more problem, DJTL, on which truncated CG makes little progress until the time limit (\Cref{tab:solved}).

We use MINRES-QLP as the default subproblem solver in the sequel.
Truncated CG remains a reasonable widely-available alternative.

\subsubsection{Outer-Loop Variants with MINRES-QLP}%
\label{sec:case-variants}

We now study the two outer-loop switches of R2N, holding MINRES-QLP and the AG handler fixed.
The Base configuration is that of \Cref{alg:r2n} with its update of \(\sigma_k\) and the parameters of \Cref{tab:constants}.
The Fast variant replaces the \(\sigma_k\) update on very successful iterations by that of \Cref{alg:r2n-fast}, which drives \(\sigma_k\) to zero.
The Skip \(\rho\) variant accepts a step obtained from an Armijo-Goldstein line search without evaluating \(\rho_k\), as motivated by \Cref{prop:ag-accept}.
The fourth configuration enables both.

\begin{figure}[ht]
  \centering
  \includetikzgraphics[width=\linewidth]{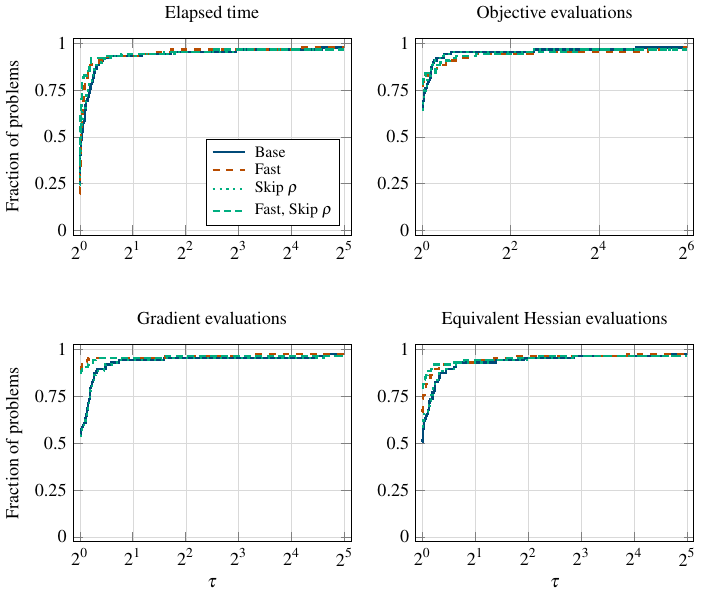}
  \caption{%
    Performance profiles of R2N with exact second derivatives, the MINRES-QLP subproblem solver and the AG handler.
    Comparison between variants: update of \(\sigma_k\) of \Cref{alg:r2n} (Base), update of \Cref{alg:r2n-fast} (Fast), acceptance of AG steps without the ratio test (Skip \(\rho\)), and both changes (Fast, Skip \(\rho\)).
  }%
  \label{fig:r2n-variants}
\end{figure}

The Fast update produces the visible difference in \Cref{fig:r2n-variants}.
Its profiles lie above those of the base update at small \(\tau\) for gradient evaluations, i.e., successful iterations, and for equivalent Hessian evaluations, while the elapsed-time and objective evaluations profiles largely overlap and robustness is unchanged (\Cref{tab:solved}).
This suggests that the base update of \Cref{alg:r2n}, which decreases \(\sigma_k\) by \(\gamma_3\) per very successful iteration, keeps \(\sigma_k\) larger than necessary as the iterates approach a minimizer.
In \Cref{sec:case-ma57}, exact subproblem solves illustrate the same effect.

The profiles with and without the Skip \(\rho\) switch are superimposed, as \Cref{prop:ag-accept} leads us to expect: an Armijo-Goldstein step satisfies \(\rho_k \geq \tilde{\eta}_1 > 0\), so that, with \(\eta_1 = \epsilon_M^{1/4}\), the ratio test can reject such a step only if \(\tilde{\eta}_1 < \eta_1\).
Because \(f(x_k + s_k)\) has already been computed by the line search, skipping the test saves no evaluation; its only effect is on the regularization parameter, which is left unchanged after every accepted Armijo-Goldstein step, whereas with the test in place a step with \(\rho_k \geq \eta_2\) decreases it.
The switch costs one problem, SCOSINE\@.
When nonpositive curvature is detected on many consecutive iterations, freezing \(\sigma_k\) may leave it at a value fixed early in the run and over-regularize the steps that follow, which is a plausible explanation of that loss.
Because skipping the test saves no evaluation and costs one problem, we keep the ratio test.

The AG handler combined with the Fast update is the configuration we carry into the comparisons that follow.

\subsubsection{Outer-Loop Variants with MA57}%
\label{sec:case-ma57}

\begin{figure}[ht]
  \centering
  \includetikzgraphics[width=\linewidth]{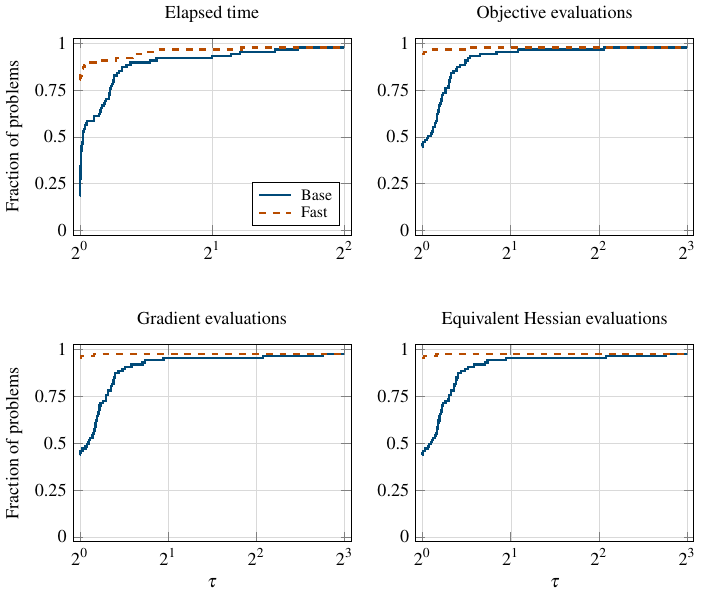}
  \caption{%
    Performance profiles of R2N with exact second derivatives and the MA57 subproblem solver.
    Comparison between the update of \(\sigma_k\) of \Cref{alg:r2n} (Base) and that of \Cref{alg:r2n-fast} (Fast).
    Nonpositive curvature is handled by rejecting the step and increasing \(\sigma_k\).
  }%
  \label{fig:r2n-ma57}
\end{figure}

We repeat the comparison of \Cref{sec:case-variants} with the direct solver MA57 \citep{duff-2004} in place of MINRES-QLP, retaining only the Fast switch.
The subproblem is then solved exactly. %
\Cref{fig:r2n-ma57} shows that the Fast update is above the base update on all four measures, including elapsed time. %

The update rule does not affect robustness.
Both configurations solve the same \(85\) problems (\Cref{tab:solved}), and both declare INDEF and HIELOW unbounded below.
Neither reaches the time limit. %
With exact subproblem solves, the Fast update should therefore be preferred; the price of such a solve, which is to assemble and factorize the Hessian at every iteration, is the subject of \Cref{sec:case-exact-vs-krylov}.

\subsubsection{MA57 versus MINRES-QLP}%
\label{sec:case-exact-vs-krylov}

\begin{figure}[ht]
  \centering
  \includetikzgraphics[width=\linewidth]{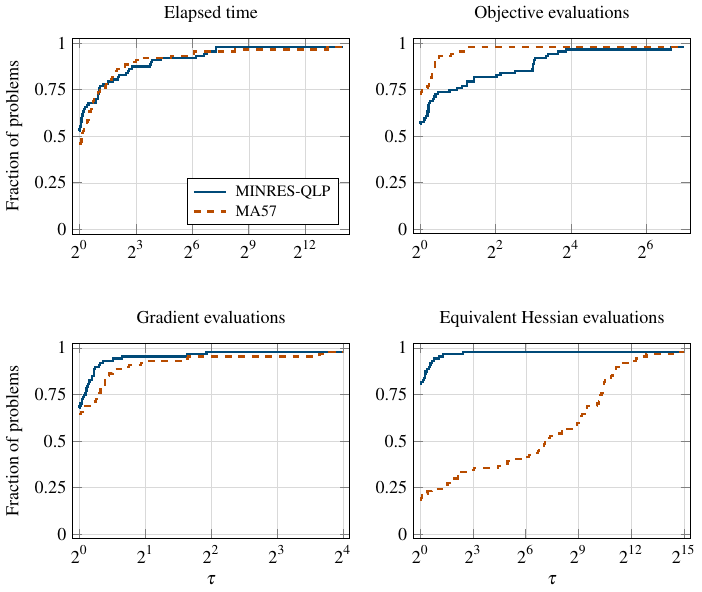}
  \caption{%
    Performance profiles of R2N with exact second derivatives and the Fast update of \(\sigma_k\).
    Comparison between subproblem solvers: MINRES-QLP with the AG handler and MA57.
    Equivalent Hessian evaluations count one assembled Hessian as \(n\) Hessian-vector products and do not reflect the cost of factorization.
  }%
  \label{fig:r2n-exact-vs-krylov}
\end{figure}

We compare the two ways of solving the subproblem with exact second derivatives, both under the Fast update of \Cref{alg:r2n-fast}: MINRES-QLP with the AG handler, and MA57 \citep{duff-2004}, which increases \(\sigma_k\) on NPC directions, and hence, requires a fresh factorization.
The equivalent Hessian panel of \Cref{fig:r2n-exact-vs-krylov} counts one assembled Hessian as \(n\) products, and therefore does not measure the relative cost of the two solvers.

The two configurations are closely matched on this test set, and their elapsed-time profiles nearly coincide.
MA57 requires fewer objective evaluations on more problems, which is consistent with an exact step being accepted more often and with the absence of a line search in the Sigma strategy, while the gradient-evaluation profiles are close.
Both solve \(85\) problems (\Cref{tab:solved}).
Besides INDEF, MINRES-QLP reaches the time limit on SCOSINE, and MA57 declares HIELOW unbounded.

We retain MINRES-QLP as the default subproblem solver, since it needs only operator-vector products and its memory requirements do not depend on fill-in.
MA57 is the better choice for small to moderate sizes when the Hessian is sparse, inexpensive to assemble and has sparse factors.

\subsubsection{R2N versus Trunk}%
\label{sec:case-cg-vs-trunk}

We now compare R2N, with truncated CG, the AG handler and the Fast update, with Trunk, the Steihaug-Toint trust-region method \citep{steihaug-1983}.
The two solvers share the model, the Krylov method, the forcing formula~\eqref{eq:forcing} and the stopping test, but they are complete solvers that differ in more than the step control.
Trunk applies a backtracking line search along rejected steps, uses its own initial radius and update constants, tightens the inner tolerance after a rejected step, and steps to the trust-region boundary along a nonpositive-curvature direction.
The comparison therefore assesses R2N against an established solver, not one ingredient in isolation.

\begin{figure}[ht]
  \centering
  \includetikzgraphics[width=\linewidth]{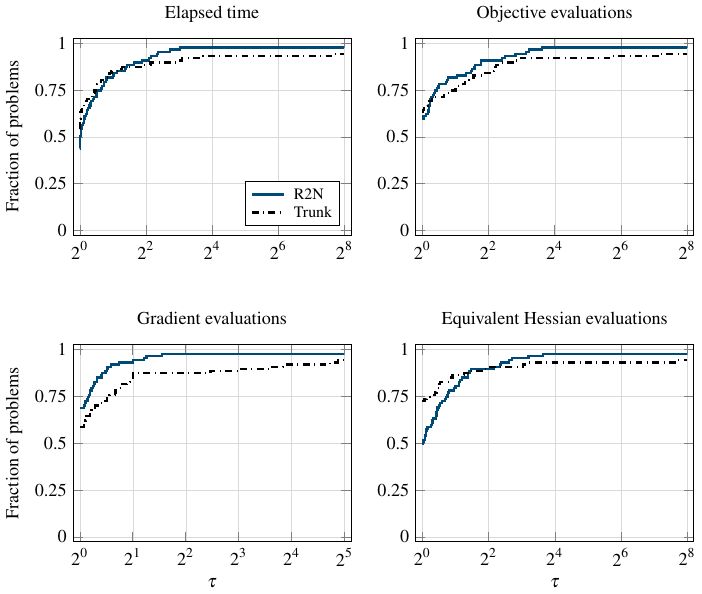}
  \caption{%
    Performance profiles with exact second derivatives and truncated CG as subproblem solver on both sides.
    Comparison between R2N with the AG handler and the Fast update of \(\sigma_k\) (R2N) and the trust-region method Trunk (Trunk).
  }%
  \label{fig:r2n-cg-vs-trunk}
\end{figure}

Neither solver dominates in \Cref{fig:r2n-cg-vs-trunk}.
Trunk is the fastest on more problems, but the elapsed-time profiles cross at small \(\tau\), and R2N is slightly more robust.
R2N requires fewer objective and gradient evaluations, hence fewer iterations.
Trunk requires fewer Hessian-vector products on more problems, but again, the profiles cross for small \(\tau\).

R2N solves three more problems than Trunk, and every problem solved by Trunk is also solved by R2N (\Cref{tab:solved}).
Of the three, Trunk reaches the time limit on SCOSINE and STRATEC, and declares HIELOW unbounded.
On INDEF, whose objective is unbounded below, R2N detects unboundedness, whereas Trunk reaches the time limit.

R2N is thus competitive with an established trust-region solver.
With exact second derivatives, R2N is preferable when gradient evaluations dominate the cost of an iteration, and Trunk when Hessian-vector products do.
\Cref{sec:case-trunk-vs-lbfgs-cg,sec:r2nls-benchmarks} repeat the comparison with models for which nonpositive curvature cannot arise.

\subsection{R2N with a Limited-Memory Model}%
\label{sec:r2n-lbfgs-benchmarks}

We now replace the exact Hessian with an L-BFGS approximation with memory \(q = 5\) described in \Cref{sec:parameters}, which satisfies \(B_k \succ 0\) at every iteration.
The shifted operator \(B_k + \sigma_k I\) is then positive definite; no direction of nonpositive curvature arises, and the NPC handler is inactive.
We compare two subproblem solvers in \Cref{sec:case-lbfgs-hessian}, and repeat the comparison with Trunk in \Cref{sec:case-trunk-vs-lbfgs-cg}.
We report elapsed time, objective evaluations and gradient evaluations.

\begin{figure}[htp]
  \centering
  \includetikzgraphics[width=\linewidth]{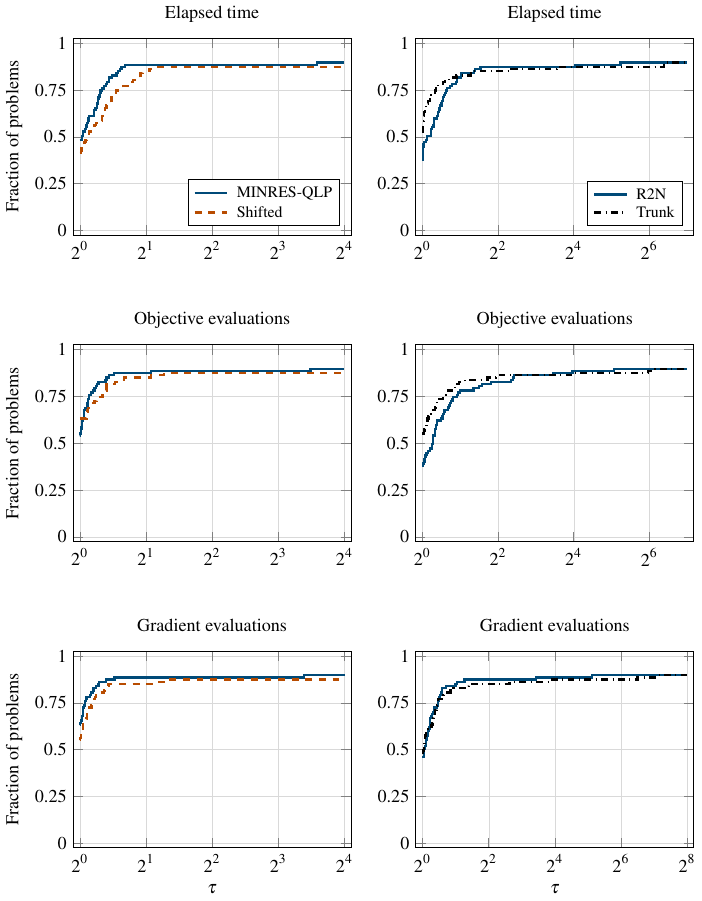}
  \caption{%
    Performance profiles with the L-BFGS model (\(q = 5\)) and the Fast update of \(\sigma_k\) in R2N\@; rows show elapsed time, objective evaluations and gradient evaluations.
    Left column: R2N with MINRES-QLP (MINRES-QLP) and with the shifted solve of \citet{erway-jain-marcia-2014} (Shifted).
    Right column: R2N and Trunk, both with truncated CG as subproblem solver.
  }%
  \label{fig:r2n-lbfgs}
\end{figure}

\subsubsection{Shifted L-BFGS Solver versus MINRES-QLP}%
\label{sec:case-lbfgs-hessian}

We compare the shifted solver of \citet{erway-jain-marcia-2014}, which exploits the representation of the L-BFGS approximation to compute the exact minimizer of~\eqref{eq:kkt-model} at a cost of \(O(q^2 n)\) per solve, with a matrix-free MINRES-QLP solve of the same model, whose cost grows with the number of inner iterations.
Both configurations use the Fast update of \Cref{alg:r2n-fast}, and the subproblem solver is the only difference.

The two are practically equivalent on this test set, and their profiles in the left column of \Cref{fig:r2n-lbfgs} are close on all three measures, with MINRES-QLP slightly above.
On this test set, solving the subproblem exactly did not reduce the number of iterations. %
MINRES-QLP solves two more problems (\Cref{tab:solved}): the shifted solver additionally reaches the time limit on JIMACK and on DJTL, where it makes slow progress.

The solver of \citet{erway-jain-marcia-2014} remains the natural choice when an exact solve of the subproblem is wanted, for instance to separate the effect of the model from that of the subproblem accuracy.

\subsubsection{R2N versus Trunk with L-BFGS}%
\label{sec:case-trunk-vs-lbfgs-cg}

We repeat the comparison of \Cref{sec:case-cg-vs-trunk} with the L-BFGS model, using CG as subproblem solver on both sides. %

The two solvers are closely matched.
Trunk is slightly faster and requires slightly fewer objective evaluations on more problems in the right column of \Cref{fig:r2n-lbfgs}, and the advantage in gradient evaluations that R2N showed with exact second derivatives is now absent.
Both solve \(78\) problems (\Cref{tab:solved}).
R2N additionally declares PARKCH unbounded, and Trunk fails with an exception on DJTL\@.
Compared to the results of \Cref{sec:case-cg-vs-trunk}, R2N loses seven problems and Trunk four, and the difference in robustness between the two disappears.

\subsection{R2NLS on Nonlinear Least-Squares Problems}%
\label{sec:r2nls-benchmarks}

We evaluate the nonlinear least-squares specialization of \Cref{alg:r2n} on the \(177\) least-squares problems described in \Cref{sec:benchmark-setup}.

We first compare R2NLS with the subproblem solvers LSMR \citep{fong-saunders-2011}, which requires only products with \(J_k\) and \(J_k^T\), and a sparse QR factorization of the matrix in~\eqref{eq:m-r2nls}, computed with QRMumps \citep{agullo-buttari-guermouche-lopez-2016}, which requires \(J_k\) explicitly.
The rest of the algorithm is identical on both sides.

In terms of iterations and residual evaluations, QRMumps is the most efficient on more problems, yet LSMR is the fastest on more problems, and the time profile of QRMumps catches up with that of LSMR only for larger \(\tau\). 
The two are equally robust, with \(169\) and \(170\) problems solved, and all failures are exhaustions of the time limit (\Cref{tab:solved}).
The solved sets differ by a few problems in each direction.
As in \Cref{sec:case-exact-vs-krylov}, the test set favors the direct solver.  %
Of course, R2NLS with LSMR requires neither the assembled Jacobian nor storage for its factors, and it remains the configuration we recommend for large problems.
For small to moderate problems, QRMumps is a highly competitive alternative. %

We repeat the comparison of \Cref{sec:case-cg-vs-trunk} on least-squares problems with TrunkNLS, the least-squares variant of Trunk in JSOSolvers.jl.
Both share the Gauss-Newton model, LSMR as subproblem solver and the stopping tests.
TrunkNLS is ahead on every measure. %
In elapsed time and in iterations, it is the most efficient solver on more than twice as many problems as R2NLS, and the gap is smaller in residual evaluations.
The gap narrows as \(\tau\) grows, but TrunkNLS solves five more problems, and every problem solved by R2NLS is also solved by TrunkNLS (\Cref{tab:solved}).
The five are exhaustions of the time limit, and the difference in robustness is small.
These results are consistent with the authors' experience with the two implementations of the Levenberg-Marquardt method.

\section{Discussion}%
\label{sec:discussion}

As \Cref{prop:succ-pc,thm:succ-npc} show, in a regularized Newton method, the worst-case complexity is dictated by positive-curvature steps, not by the subproblem solver.
Detecting nonpositive curvature, with MINRES or otherwise, does not improve the worst-case complexity of \Cref{alg:r2n}---see \Cref{thm:slow-conv}.

In practice, however, exploiting nonpositive curvature matters for performance.
R2N only has an edge over a standard trust-region method when nonpositive curvature can arise---see \Cref{fig:r2n-cg-vs-trunk,fig:r2n-lbfgs}---, which suggests that its treatment deserves as much care in regularization methods as it traditionally receives in trust-region methods.

When subproblems are convex, we observe no advantage of regularization over a trust region.
The main argument in favor of regularization is then that subproblems are unconstrained, so that direct solvers and specialized solvers such as that of \citet{erway-jain-marcia-2014} apply without modification.

\citet{kanzow-steck-2023} use limited-memory quasi-Newton models and solve the subproblem exactly with variants of the direct method of \citet{erway-jain-marcia-2014}.
They assume that \(\{B_k\}\) is bounded, and increase the regularization parameter when nonpositive curvature occurs.
While it is known that several limited-memory quasi-Newton approximations remain bounded under standard assumptions \citep{aravkin-baraldi-orban-2021,burdakov-gong-zikrin-yuan-2017,diouane-habiboullah-khouna-orban-2026}, boundedness remains an open question for full quasi-Newton approximations.
Our analysis covers inexact solves and potentially unbounded \(\{B_k\}\).
In addition, our results suggest that exact solves with limited-memory quasi-Newton models yield no discernible advantage.
In \cref{sec:case-lbfgs-hessian}, the direct solver did not reduce the number of iterations compared to MINRES-QLP, and solved two fewer problems.
Our approach and that of \citet{kanzow-steck-2023} differ most with indefinite models such as L-SR1, where R2N exploits nonpositive curvature instead of increasing \(\sigma_k\).
Experiments with L-SR1 models, for R2N and Trunk alike, are the natural next step.

Several improvements to the solver are possible.
Our density test is a proxy for the cost of a factorization.
The symbolic analysis phase of a sparse factorization predicts the size of the factors, and could be used to choose between a direct and an iterative solver automatically.
Shift invariance of the Krylov subspace could be exploited after a rejected step, either by reusing the Krylov subspace generated for the previous value of \(\sigma_k\), or by solving for several values of \(\sigma_k\) at once with a shifted Krylov method such as the shifted variant of CG-Lanczos \citep{frommer-maass-1999}.
A similar idea is exploited in \citep{dussault-migot-orban-2024} in the context of a cubic regularization method.
Preconditioning the Krylov solvers, e.g., with an incomplete factorization, is a natural way to combine both approaches, provided the preconditioner is positive definite so that nonpositive-curvature detection remains meaningful.
Other Krylov solvers are potential candidates to solve the subproblem.
A recent example is MinAres \citep{montoison-orban-saunders-2025}, which minimizes \(\|A r_k\|\) and is well suited to inconsistent systems, which could occur when \(B_k + \sigma_k I\) is singular.
Using it in R2N requires nonpositive-curvature detection and a Cauchy decrease guarantee similar to those of \cref{sec:truncated-cr,sec:modified-minres}.

Several questions remain open.
\citet{gratton-jerad-toint-2025} regularize the Hessian with a multiple of \(\|\nabla f(x_k)\|^{1/2}\) and obtain a bound in \(O(|\log \epsilon| \, \epsilon^{-3/2})\).
Whether such a rule can be combined with the update of \Cref{alg:r2n}, and how it affects the local rate, remains to be investigated.
Equipping Trunk with MINRES-QLP would tell whether the advantage of MINRES-QLP over truncated CG carries over to trust-region methods.
Finally, we do not know of an example that attains the bound of \Cref{thm:succ-npc} when \(0 < \zeta < 1\), or whether the bound of \Cref{thm:unsuccessful-fast} is tight.

\subsection*{Use of AI Tools}%

We used the large language model Claude Opus 5.5 (Anthropic) during the preparation of this paper.
Below, we list each section of the paper and state whether it was written entirely by the human authors or whether AI contributed, and how.
In every section, AI was used to edit the English and make the text more readable.
All AI contributions were checked by the human authors, who also cross-checked each other's work, and the human authors take full responsibility for the content of this paper.

\begin{description}
  \item[Abstract, Introduction, \nameref{sec:background}, \nameref{sec:models}.] Written by the human authors.
  \item[\nameref{sec:algorithm}, \nameref{sec:step}.] The statement and proof of all results were developed and verified by the human authors, and AI was used to re-verify the proofs.
  \item[\nameref{sec:convergence}.] Written by the human authors.
  \item[\nameref{sec:complexity}.] The statements and proofs of all results were developed and verified by the human authors, and AI was used to re-verify the proofs.
  \item[\nameref{sec:slow}.] The example, together with its lemma and proof, was constructed and verified by the human authors, re-verified with AI, and finally checked numerically to confirm that it satisfies the assumptions of the corresponding results.
  \item[\nameref{sec:conv-npc}.] The design of the sequences \(\{f_k\}\), \(\{g_k\}\), \(\{H_k\}\) and \(\{s_k\}\) was done jointly with AI\@.
  The assumptions of \Cref{thm:conv-npc-example} were formulated by the human authors with the help of AI\@.
  \Cref{lem:f-g-bnd-npc,lem:Hermite-spline,thm:conv-npc-example} and their proofs were developed by both the human authors and AI\@: AI was used to re-evaluate them and to check that the assumptions hold, and the human authors extended the lemmas and their proofs.
  The final construction was re-evaluated by the human authors and checked numerically to confirm that it satisfies the assumptions of the corresponding results.
  \item[\nameref{sec:local}.] The statements and proofs of all results were developed and verified by the human authors, and AI was used to re-verify the proofs.
  \item[\nameref{sec:nlls}.] Written by the human authors.
  \item[\nameref{sec:implementation}.] The text was written by the human authors. Most of the code was written by F.~Rahbarnia, and an AI coding agent was used to clean and review the code as needed. The implementation and the numerical results were validated by the human authors.
  \item[\nameref{sec:discussion}.] The human authors outlined the content of this section. AI was used to turn the outline into sentences, which the human authors edited for readability.  
  \item[References.] AI was used to check bibliographic data and to format entries, all of which were checked by the human authors.
\end{description}

\small
\bibliographystyle{abbrvnat}
\bibliography{abbrv,r2n-smooth}
\normalsize

\section*{Declarations}

\paragraph{Availability of data and materials}
All test problems are from the CUTEst collection, which is publicly available.
The results of all experiments are available at \url{https://github.com/farhadrclass/r2n-smooth-benchmarks}.

\paragraph{Code availability}
R2N and R2NLS are part of JSOSolvers.jl v0.14.9, available at \url{https://github.com/JuliaSmoothOptimizers/JSOSolvers.jl}.
The scripts that reproduce all numerical results are available at \url{https://github.com/farhadrclass/r2n-smooth-benchmarks}.

\newpage

\hypertarget{contents}{}  %
\scriptsize
\tableofcontents
\normalsize

\end{document}